\documentclass[twoside,11pt,reqno]{amsart}
\usepackage{divres}

\begin{document}

\title[Ext algebras for standard modules]{Some extension algebras for standard modules over KLR algebras of type $A$}

\author{\sc Doeke Buursma}
\address{Department of Mathematics\\University of Oregon\\Eugene\\OR 97403, USA}
\email{dbuursma@uoregon.edu}

\author{\sc Alexander Kleshchev}
\address{Department of Mathematics\\University of Oregon\\Eugene\\OR 97403, USA}
\email{klesh@uoregon.edu}

\author{\sc David J. Steinberg}
\address{Department of Mathematics\\University of Oregon\\Eugene\\OR 97403, USA}
\email{dsteinbe@uoregon.edu}


\thanks{The second author was supported by the NSF grants DMS-1161094, DMS-1700905, the Max-Planck-Institut, the Fulbright Foundation, and the DFG Mercator program through the University of Stuttgart.}

\begin{abstract}
	Khovanov-Lauda-Rouquier algebras $R_\theta$ of finite Lie type are affine quasihereditary with standard modules $\Delta(\pi)$ labeled by Kostant partitions of $\theta$.
	Let $\Delta$ be the direct sum of all standard modules.
	It is known that the Yoneda algebra $\CE_\theta:=\Ext_{R_\theta}^*(\Delta, \Delta)$ carries a structure of an $A_\infty$-algebra which can be used to reconstruct the category of standardly filtered $R_\theta$-modules.
	In this paper, we explicitly describe $\CE_\theta$ in two special cases: (1) when $\theta$ is a positive root in type $\ttA$, and (2) when $\theta$ is of Lie type $\ttA[2]$.
	In these cases, $\CE_\theta$ turns out to be torsion free and intrinsically formal.
	We provide an example to show that the $A_\infty$-algebra $\CE_\theta$ is non-formal in general.
\end{abstract}

\maketitle

\section{Introduction}\label{SIntro}

Let $R_{\theta, \F}$ be a Khovanov-Lauda-Rouquier (KLR) algebra of finite Lie type over a field $\F$ corresponding to $\theta \in Q_+$~\cite{KL,R}.
It is known that $R_{\theta, \F}$ is affine quasihereditary~\cite{BKM,Kat,Kdonkin,KlL}, and in particular it comes with a family of \emph{standard modules} $\{ \Delta(\pi)_\F \mid \pi\in \KP(\theta)\}$, where $\KP(\theta)$ is the set of Kostant partitions of $\theta$.
KLR algebras are defined over $\Z$, so we have a $\Z$-algebra $R_{\theta}$ with $R_{\theta, \F} \cong R_{\theta} \otimes_\Z \F$.
The standard modules have natural integral forms $\Delta(\pi)$ with $\Delta(\pi)_\F \cong \Delta(\pi) \otimes_\Z \F$.
All modules and algebras are explicitly graded, and we refer to these gradings as KLR gradings.

Let $\Delta := \bigoplus_{\pi\in \KP(\theta)}\Delta(\pi)$.
The Yoneda algebra $\CE_\theta:=\Ext_{R_\theta}^*(\Delta, \Delta)$ carries a structure of an $A_\infty$-algebra~\cite{Kad}, which can be used to reconstruct the category $\mathcal{F}(\Delta)$ of modules which admit a finite filtration by standard modules~\cite{Ke1,Ke2,KKO}.
In view of~\cite[Corollary 3.14]{BKM} and~\cite[Theorem 4.28]{KS}, understanding $\mathcal{F}(\Delta)$ is relevant for computing formal characters of the simple modules $L(\pi)_\F$ of $R_{\theta, \F}$.
The smallest known example where the formal characters depend on the characteristic of $\F$ occurs in Lie type $\ttA[5]$, see~\cite{Wi} (cf.~\cite[\S2.6]{BKM}).

We now assume that the Lie type is $\ttA[\infty]$ with simple roots $\{\alpha_i\mid i\in\Z\}$ so that the set of positive roots $\Phi_+$ is $\{\alpha_i+\alpha_{i+1}+\dots+\alpha_j\mid i\leq j\}$.
There is a natural lexicographic total order $>$ on $\Phi_+$.
Let $Q_+$ be the positive root lattice, and fix $\theta\in Q_+$.
If $\theta=\sum k_i\alpha_i$, we define the \emph{height} of $\theta$ as $\height(\theta):=\sum k_i$.
A \emph{Kostant partition} of $\theta$ is a sequence $\pi = (\beta_1^{m_1}, \ldots, \beta_t^{m_t})$ where $m_1,\dots, m_t\in \Z_{>0}$, $\beta_1 > \cdots > \beta_t$ are positive roots, and $m_1\beta_1 + \cdots + m_t\beta_t = \theta$.

We consider the Yoneda algebra $\CE_\theta$ as the $\Z$-linear category whose objects are $\KP(\theta)$, and the set of morphisms from $\rho\in \KP(\theta)$ to $\sigma\in \KP(\theta)$ is
$$\CE_\theta(\rho, \sigma) := \Ext^*_{R_\theta}(\Delta(\rho), \Delta(\sigma)).$$
The composition $gf$ of $g\in \CE_\theta(\sigma, \tau)$ and $f\in \CE_\theta(\rho, \sigma)$ is obtained using the composition of lifts of $g$ in $\Hom_{R_\theta}(P_\bullet^\sigma, P_\bullet^\tau)$ and $f$ in $\Hom_{R_\theta}(P_\bullet^\rho, P_\bullet^\sigma)$, where $P_\bullet^\pi$ is a projective resolution of $\Delta(\pi)$ for $\pi \in \KP(\theta)$.
The category $\CE_\theta$ has a \emph{homological grading} for which the homogeneous components are $\CE_\theta^m(\rho, \sigma) := \Ext_{R_\theta}^m(\Delta(\rho), \Delta(\sigma))$, and a \emph{KLR grading} which is inherited from the KLR grading on the standard modules.
We use $q$ to denote the KLR degree shift functor.
Theorems~\ref{TA} and~\ref{TB} describe the category $\CE_\theta$ (as a bigraded category) in two special cases: (1) when $\theta$ is an arbitrary positive root, and (2) when $\theta$ is of type $\ttA[2]$, i.e. $\theta$ is of the form $c_1\alpha_1 + c_2\alpha_2$.

\subsection{The case where $\theta$ is a positive root}

Let $\theta = \alpha_a + \alpha_{a+1} + \cdots + \alpha_{b+1} \in \Phi_+$.
Set $l := b+2-a=\height(\theta)$ and consider the polynomial algebra $\CX := \Z[x_1, \ldots, x_l]$.
We consider $\CX$ to be graded with $\deg x_r = 2$.
Note that $\KP(\theta)$ is in bijection with the set of subsets of $[1,l-1]$: the subset associated to $\rho = (\beta_1, \ldots, \beta_u) \in \KP(\theta)$ is $D_\rho := \{ d_1, \ldots, d_{u-1}\}$ where $d_t := \height(\beta_1)+\cdots+\height(\beta_t)$.
For such $D_\rho$, set $d_0 := 0$ and $d_u := l$, and let $J^\rho$ be the ideal of $\CX$ generated by all $x_r-x_s$ such that there is $1\leq t\leq u$ with $d_{t-1} < r,s \leq d_t.$
Define $\CX^\rho := \CX/J^\rho$.
If $D_\rho \subseteq D_\sigma$, then $J^\sigma\subseteq J^\rho$ so we have a natural projection $\ttp_\rho^\sigma : \CX^\sigma \onto \CX^\rho$.
We use the notation $C\subseteq_m D$ to indicate that $C\subseteq D$ with $|D\setminus C| = m$.

\begin{MainTheorem}\label{TA}
	Let $\theta = \alpha_a + \alpha_{a+1} + \cdots + \alpha_{b+1} \in \Phi_+$ be a positive root.
	We have
	$$\CE_\theta^m(\rho,\sigma) \cong
	\begin{cases}
		q^{-m}\CX^\rho & \text{if } D_\rho \subseteq_m D_\sigma, \\
		0 & \text{otherwise}.
	\end{cases}$$
	If $D_\rho \subseteq_m D_\sigma \subseteq_n D_\tau$ with $f\in q^{-m}\CX^\rho \cong \CE_\theta(\rho, \sigma)$ and $g\in q^{-n}\CX^\sigma \cong \CE_\theta(\sigma, \tau)$, then the composition of $g$ with $f$ is given by $\ttp_\rho^\sigma(g)f \in q^{-(m+n)}\CX^\rho \cong \CE_\theta(\rho, \tau)$.
\end{MainTheorem}

\subsection{The $\ttA[2]$ case}

For a nonnegative integer $k$, let $\Lambda_k$ be the algebra of symmetric polynomials in $k$ variables.
We impose a grading on $\Lambda_k$ where linear symmetric polynomials have degree $2$.
The space $\Lambda_k$ is a free $\Z$-module with basis $\{ s_\lambda \mid \lambda \in \Par(k)\}$, where $\Par(k)$ is the set of partitions with at most $k$ parts, and $s_\lambda$ is the \emph{Schur polynomial} corresponding to $\lambda$~\cite[\S I.3]{Mac}.
Letting $V$ be the free graded $\Z$-module with basis $\{ v_0, v_1, v_2, \ldots\}$ such that $\deg v_i := 2i$, there is an isomorphism of graded $\Z$-modules
\begin{equation}\label{EGamma}
	\gamma_k: \Lambda_k \iso q^{-k(k-1)}\Wed^k V, \ s_{(\lambda_1, \ldots, \lambda_k)} \mapsto v_{\lambda_k}\wedge v_{\lambda_{k-1}+1}\wedge \cdots \wedge v_{\lambda_1+k-1}
\end{equation}
where $\Wed^k V$ is the $k\mathrm{th}$ exterior power of $V$.
Define
\begin{equation}\label{EStar}
	\blank \star \blank : \Lambda_a \otimes \Lambda_b \to q^{2ab} \Lambda_{a+b}, \ f\otimes g \mapsto \gamma_{a+b}^{-1}(\gamma_a(f)\wedge \gamma_b(g)).
\end{equation}

Considering $\Lambda_{a+b}$ to be a subalgebra of $\Lambda_a\otimes \Lambda_b$ in the obvious way, we have that $\Lambda_a \otimes \Lambda_b$ is free as a $\Lambda_{a+b}$-module with basis $\{ s_\lambda \otimes 1 \mid \lambda \in \Par(a,b)\}$, where $\Par(a,b)$ is the set of partitions with at most $a$ nonzero parts, the first part being at most $b$, see~\cite[PARTL.1.5]{LLPT} and~\cite[Proposition 2.6.8]{Man}.
Moreover,~\cite[SCHUB.1.7]{LLPT} provides an explicit algorithm for writing any element of $\Lambda_a\otimes \Lambda_b$ as a $\Lambda_{a+b}$-linear combination of the basis elements $s_\lambda \otimes 1$.

Let $c_1, c_2 \in \Z_{\geq 0}$ and $\theta = c_1 \alpha_1 + c_2 \alpha_2$.
Note that there is a bijection
$$[0, \min\{ c_1, c_2\} ] \overset{\sim}{\longleftrightarrow} \KP(\theta), \ r\mapsto (\alpha_2^{c_2-r}, (\alpha_1+\alpha_2)^r, \alpha_1^{c_1-r}), \ r_\rho \mapsfrom \rho.$$
For $\rho, \sigma \in \KP(\theta)$ with $r_\rho \geq r_\sigma$, let
\begin{align*}
	\omega(\rho, \sigma) &:= -(r_\rho-r_\sigma)(1+(c_1-r_\rho)+(c_2-r_\rho)), \\
	\Lambda(\rho,\sigma) &:= q^{\omega(\rho,\sigma)}\Lambda_{c_2-r_\rho}\otimes \Lambda_{r_\rho-r_\sigma}\otimes \Lambda_{r_\sigma}\otimes \Lambda_{c_1-r_\rho}, \\
	\Par_{\rho,\sigma}&:=\Par(r_\rho-r_\sigma,r_\sigma).
\end{align*}
If $f\in \Lambda_{r_\rho-r_\sigma}$, we write
$$f^{\rho,\sigma}:=1_{\Lambda_{c_2-r_\rho}}\otimes f\otimes 1_{\Lambda_{r_\sigma}}\otimes 1_{\Lambda_{c_1-r_\rho}}\in \Lambda(\rho,\sigma).$$
Then note that $\Lambda(\rho,\sigma)$ is a free right $\Lambda(\rho,\rho)$-module with basis $\{ s_\lambda^{\rho,\sigma} \mid \lambda \in \Par_{\rho,\sigma}\}$.
We make $\Lambda(\rho, \sigma)$ into a left $\Lambda(\sigma, \sigma)$-module via the composition of algebra homomorphisms:
$$\xi: \Lambda(\sigma,\sigma) \into \Lambda_{c_2-r_\rho}\otimes\Lambda_{r_\rho-r_\sigma}\otimes\Lambda_{r_\sigma}\otimes\Lambda_{r_\rho-r_\sigma}\otimes\Lambda_{c_1-r_\rho} \onto q^{-\omega(\rho,\sigma)}\Lambda(\rho,\sigma);$$
the first map uses the embeddings $\Lambda_{c_2-r_\sigma}\into \Lambda_{c_2-r_\rho}\otimes\Lambda_{r_\rho-r_\sigma}$ and $\Lambda_{c_1-r_\sigma}\into \Lambda_{r_\rho-r_\sigma}\otimes\Lambda_{c_1-r_\rho}$, and the second map is $a\otimes b\otimes c\otimes d\otimes e\mapsto a\otimes bd\otimes c\otimes e$ (which we think of as identifying the two factors of $\Lambda_{r_\rho-r_\sigma}$).

If $\rho, \sigma, \tau \in \KP(\theta)$ with $r_\rho \geq r_\sigma \geq r_\tau$, the tensor product $\Lambda(\sigma,\tau)\otimes_{\Lambda(\sigma,\sigma)}\Lambda(\rho,\sigma)$ is now a free right $\Lambda(\rho,\rho)$-module with basis
$$\{ s_\mu^{\sigma,\tau} \otimes s_\lambda^{\rho,\sigma}\mid \mu\in \Par_{\sigma,\tau}, \lambda\in \Par_{\rho,\sigma}\}$$
and we define a map of right $\Lambda(\rho,\rho)$-modules
$$\Theta: \Lambda(\sigma,\tau)\otimes_{\Lambda(\sigma,\sigma)}\Lambda(\rho,\sigma) \to \Lambda(\rho,\tau), \ s_\mu^{\sigma,\tau}\otimes s_\lambda^{\rho,\sigma} \mapsto (s_\mu\star s_\lambda)^{\rho,\tau}.$$
Let
$$-\diamond- : \Lambda(\sigma,\tau)\otimes_\Z \Lambda(\rho,\sigma) \to \Lambda(\rho,\tau), \ g\otimes f \mapsto \Theta(g\otimes f).$$
Thus, to compute $g\diamond f$ for some $g\in \Lambda(\sigma, \tau)$ and $f \in \Lambda(\rho, \sigma)$, the following steps must be performed: (1) write $g = \sum_{\mu \in \Par_{\sigma,\tau}}  s_\mu^{\sigma,\tau}g_\mu$ with $g_\mu \in \Lambda(\sigma, \sigma)$, (2) for each $\mu\in \Par_{\sigma,\tau}$, write $\xi(g_\mu)f = \sum_{\lambda\in \Par_{\rho,\sigma}} s_\lambda^{\rho,\sigma}h_{\mu,\lambda}$ with $h_{\mu,\lambda} \in \Lambda(\rho, \rho)$, (3) we have
$$g\diamond f = \sum_{\mu\in \Par_{\sigma,\tau},\, \lambda \in \Par_{\rho,\sigma}}(s_\mu \star s_\lambda)^{\rho,\tau}h_{\mu, \lambda}.$$

\begin{MainTheorem}\label{TB}
	Let $c_1, c_2 \in \Z_{\geq 0}$ and $\theta = c_1 \alpha_1 + c_2 \alpha_2$.
	We have
	$$\CE_\theta^m(\rho, \sigma) =
	\begin{cases}
		\Lambda(\rho, \sigma) & \text{if } m = r_\rho-r_\sigma \geq 0, \\
		0 & \text{otherwise}.
	\end{cases}$$
	If $r_\rho \geq r_\sigma \geq r_\tau$ with $f\in \Lambda(\rho,\sigma) \cong \CE_\theta(\rho, \sigma)$ and $g\in \Lambda(\sigma,\tau) \cong \CE_\theta(\sigma, \tau)$, then the composition of $g$ with $f$ is given by $g\diamond f \in \Lambda(\rho,\tau) \cong \CE_\theta(\rho, \tau)$.
\end{MainTheorem}

\subsection{Formality}

For $\theta$ as in Theorem~\ref{TA} or~\ref{TB}, note that $\Ext_{R_\theta}(\Delta(\rho), \Delta(\sigma))$ is torsion-free as a $\Z$-module.
We do not know if this is true in general.

Also note that because $\Ext_{R_\theta}(\Delta(\rho), \Delta(\sigma))$ is concentrated in homological degree $|X_\sigma \setminus X_\rho|$ (in the case of Theorem~\ref{TA}) or $r_\rho-r_\sigma$ (in the case of Theorem~\ref{TB}), the $A_\infty$-category structure of $\CE_\theta$ must have $m_n = 0$ unless $n = 2$, so that $\CE_\theta$ is \emph{intrinsically formal}, see~\cite[\S3.3]{Ke2}.
In Section~\ref{SExample}, we show that intrinsic formality and even formality does not occur in general:

\begin{MainExample}\label{EC}
	If $\theta = \alpha_1 + 2\alpha_2 + \alpha_3$, then the $A_\infty$-category $\CE_\theta$ is non-formal.
\end{MainExample}

\subsection{The structure of the paper}

The proofs of Theorems~\ref{TA} and~\ref{TB} occupy Sections~\ref{SRoot} and~\ref{SA2}, respectively.
In the preliminary Section~\ref{SPrelim}, we review the definition of the KLR algebra $R_\theta$ and the standard modules $\Delta(\pi)$.

In \S\ref{SSRootRes} and \S\ref{SSA2Res}, we record the relevant special cases of the projective resolution $P_\bullet^\rho$ of $\Delta(\rho)$ constructed in~\cite{BKS}.
This resolution is finite and has the form $P_\bullet^\rho = \cdots \overset{d_1}{\to} P_1^\rho \overset{d_0}{\to} P_0^\rho \overset{\eps_\rho}{\to} \Delta(\rho)$ with $P_n^\rho = \bigoplus_{x\in X_n}q^{s_x}R_\theta 1_x$ for some explicit index set $X_n$, integers $s_x$, and idempotents $1_x$.
The map $d_n: P_{n+1}^\rho \to P_n^\rho$ can be described as right multiplication by an $X_{n+1}\times X_n$ matrix $(d_n^{y,x})$ for some $d_n^{y,x} \in 1_y R_\theta 1_x$.

In \S\ref{SSRootExtMod} and \S\ref{SSA2ExtMod}, we use the isomorphism $\Hom_{R_\theta}(R_\theta 1_x, \Delta(\sigma)) \cong 1_x \Delta(\sigma)$ to describe the complex $\Hom_{R_\theta}^\bullet(P_\bullet^\rho, \Delta(\sigma))$ in terms of familiar objects from commutative algebra; in the case of Theorem~\ref{TA}, these objects are polynomial rings and in the case of Theorem~\ref{TB}, they are rings of symmetric polynomials.
It turns out that in both cases, the complex $\Hom_{R_\theta}^\bullet(P_\bullet^\rho, \Delta(\sigma))$ is isomorphic to a Koszul complex corresponding to an explicit regular sequence, and we can therefore compute its homology $H(\Hom_{R_\theta}^\bullet(P_\bullet^\rho, \Delta(\sigma))) =: \CE_\theta(\rho, \sigma)$ as a bigraded $\Z$-module.

It remains to describe the composition in the category $\CE_\theta$.
This is done in \S\ref{SSRootCat} and \S\ref{SSA2Cat}, where we explicitly lift elements of $\Hom_{R_\theta}^\bullet(P_\bullet^\rho, \Delta(\sigma))$ to $\Hom_{R_\theta}^\bullet(P_\bullet^\rho, P_\bullet^\sigma)$.
The function composition map $\Hom_{R_\theta}^\bullet(P_\bullet^\sigma, P_\bullet^\tau) \otimes \Hom_{R_\theta}^\bullet(P_\bullet^\rho, P_\bullet^\sigma) \to \Hom_{R_\theta}^\bullet(P_\bullet^\rho, P_\bullet^\tau)$ induces a map on homology $\CE_\theta(\sigma, \tau) \otimes \CE_\theta(\rho, \sigma) \to \CE_\theta(\rho, \tau)$ which is the composition in the category $\CE_\theta$.

In Section~\ref{SExample}, we review the $A_\infty$-category structure on $\CE_\theta$ following~\cite{Kad,Ke1,Ke2} and provide details for Example~\ref{EC}.

\section{Preliminaries}\label{SPrelim}

\subsection{Basic notation}

Throughout, we work over an arbitrary principal ideal domain~$\k$ (since everything is defined over $\Z$, one could just consider the case $\k=\Z$).

For $r,s\in\Z$, we use the segment notation $[r,s]:=\{ t\in\Z\mid r\leq t\leq s\}$, $[r,s):=\{ t\in\Z\mid r\leq t< s\}$, etc.

Let $q$ be a variable, and $\Z((q))$ be the ring of Laurent series. For $n \in \Z_{\geq 0}$, we define
$$[n] := \frac{q^n-q^{-n}}{q-q^{-1}}, \quad [n]_\pm := q^{\pm(n-1)}[n], \quad [n]_{(\pm)}^!:=[1]_{(\pm)}[2]_{(\pm)}\cdots [n]_{(\pm)},$$
and if $0 \leq m \leq n$,
$$\qbinom{n}{m}_{(\pm)} := \frac{[n]_{(\pm)}^!}{[m]_{(\pm)}^![n-m]_{(\pm)}^!}.$$

We denote by $\Sgp_d$ the symmetric group on $d$ letters considered as a Coxeter group with generators $\{s_r:=(r,r+1)\mid 1\leq r<d\}$ and the corresponding length function $\ell$.
The longest element of $\Sgp_d$ is denoted $w_0$ or $w_{0,d}$.
By definition, $\Sgp_d$ acts on $[1,d]$ on the left.
For a set $I$ the $d$-tuples from $I^d$ are written as words $\bi=i_1\cdots i_d$.
The group $\Sgp_d$ acts on $I^d$ via place permutations: $w\cdot \bi=i_{w^{-1}(1)}\cdots i_{w^{-1}(d)}$.

Given a composition $\mu=(\mu_1, \ldots, \mu_k)$ of $d$, we have the corresponding standard parabolic subgroup $\Sgp_\mu:=\Sgp_{\mu_1} \times \cdots \times \Sgp_{\mu_k} \leq \Sgp_d$.
We denote by $\D^\mu$ the set of the shortest coset representatives for $\Sgp_d/\Sgp_\mu$.

\subsection{Symmetric polynomials and the nil-Hecke algebra}\label{SSSym}

We impose a grading on the polynomial algebra $\CX_d:=\k[x_1,\dots,x_d]$ such that $\deg(x_r) = 2$.
The symmetric group $\Sgp_d$ acts by automorphisms on $\CX_d$ via $(w\cdot f)(x_1, \ldots, x_d) := f(x_{w(1)}, \ldots, x_{w(d)})$.
The symmetric polynomial algebra $\Lambda_d := \CX_d^{\Sgp_d}$ has a basis consisting of Schur polynomials
$$\{ s_\lambda \mid \lambda \in \Par(d)\},$$
where $\Par(d)$ is the set of partitions with at most $d$ nonzero parts, see~\cite[\S I.3]{Mac}.

For a composition $\mu = (\mu_1, \ldots, \mu_k)$ of $d$, the algebra of $\mu$-partially symmetric polynomials is $\Lambda_\mu := \CX_d^{\Sgp_\mu}$.
We often write $\Lambda_{\mu_1,\dots,\mu_k}$ for $\Lambda_\mu$ and identify it with $\Lambda_{\mu_1} \otimes \cdots \otimes \Lambda_{\mu_k}$.
For $a,b \in \Z_{\geq 0}$, let $\Par(a, b)$ be the set of partitions with at most $a$ nonzero parts, each part being at most $b$.
The following is known:

\begin{Proposition}\label{PPartSymBasis}
	The algebra $\Lambda_\mu$ is free as a $\Lambda_d$-module, and in the case where $\mu = (a,b)$ with $a+b=d$, a basis is given by
	$$\{ s_\lambda \otimes 1 \mid \lambda\in \Par(a, b)\}.$$
\end{Proposition}
\begin{proof}
	The freeness assertion is~\cite[PARTL.1.5]{LLPT}, which, when combined with~\cite[Proposition 2.6.8]{Man}, gives the basis assertion.
\end{proof}

For an integer $r$ with $1\leq r < d$ and a reduced decomposition $w = s_{r_1}\cdots s_{r_k}\in \Sgp_d$, the Demazure operators on $\CX_d$ are defined as follows:
$$\partial_r := \frac{\id_{\CX_d} - s_r}{x_{r+1} - x_r} \quad \text{and} \quad \partial_w := \partial_{r_1} \cdots \partial_{r_k}.$$
Note that $\partial_w$ does not depend on the choice of reduced decomposition and is a degree $-2\ell(w)$ element of $\End_\k{\CX_d}$.

For integers $r,i,j \geq 0$ with $r+i+j \leq d$, define
\begin{equation}\label{EThickTransposition}
	U_{r;i,j}\in \Sgp_d
\end{equation}
to be the permutation which maps the interval $[r+1, r+i]$ increasingly onto $[r+1+j, r+i+j]$, and the interval $[r+i+1, r+i+j]$ increasingly onto $[r+1, r+j]$, and fixes all other elements of $[1,d]$.
For example, we have $U_{r;1,1} = s_{r+1} = (r+1 \ r+2)$.
Recalling (\ref{EStar}), we have:

\begin{Proposition}~\cite[Proposition 2.9]{KLMS}\label{PThickPartial}
	Let $a,b \in \Z_{\geq 0}$ with $f\in \Lambda_a$ and $g\in \Lambda_{b}$.
	Then $\partial_{U_{0;a,b}}(f\otimes g) = f\star g$.
\end{Proposition}
\begin{proof}
	Since we use different conventions from~\cite{KLMS}, we provide a  translation for the reader's convenience.
	For $w\in \Sgp_{a+b}$, let $\partial_w' := (-1)^{\ell(w)}\partial_w$.
	If $k\in \Z_{\geq 0}$, recalling (\ref{EGamma}), let $\gamma_k' := (-1)^{\binom{k}{2}}\gamma_k$ and note that $\gamma_k'(s_\lambda) = v_{\lambda_1+k-1}\wedge v_{\lambda_2+k-2}\wedge \cdots \wedge v_{\lambda_k} \in \Wed^k V$ for any $\lambda = (\lambda_1, \ldots, \lambda_k)\in \Par(k)$.
	By~\cite[Proposition 2.9]{KLMS}, we have $\partial_{U_{0;a,b}}'(f\otimes g) = (\gamma_{a+b}')^{-1}m(\gamma_a'\otimes \gamma_b')(f\otimes g)$, so
	\begin{align*}
		\partial_{U_{0;a,b}}(f\otimes g) &= (-1)^{\ell(w)}\partial_{U_{0;a,b}}'(f\otimes g) \\
		&= (-1)^{ab}(\gamma_{a+b}')^{-1}m(\gamma_a'\otimes \gamma_b')(f\otimes g) \\
		&= (-1)^{ab + \binom{a+b}{2} + \binom{a}{2} + \binom{b}{2}}(\gamma_{a+b})^{-1}m(\gamma_a\otimes \gamma_b)(f\otimes g) \\
		&= f \star g. \qedhere
	\end{align*}
\end{proof}

Define
$$x_0 = x_{0,d} := \prod_{r=1}^d x_r^{r-1} \in \CX_d.$$
The following is well-known and easy to check:

\begin{Lemma}\label{LLambda}
	We have $\partial_{w_0}(fx_0)=f$ for all $f\in\Lambda_d$.
\end{Lemma}

The \emph{nil-Hecke algebra} $\NH_d$ is given by generators $\tau_1,\dots,\tau_{d-1},\, x_1,\dots,x_d$ subject only to the relations
\begin{align*}
  x_r x_t=x_t x_r,\ \tau_r^2=0,\ \tau_r\tau_{r+1}\tau_r=\tau_{r+1}\tau_r\tau_{r+1},\ \tau_r\tau_s=\tau_s\tau_r \ (|r-s|>1), \\
  \tau_r x_r=x_{r+1}\tau_r-1,\ \tau_r x_{r+1}=x_r\tau_r+1,\ \tau_r x_s=x_s\tau_r \ (s\neq r,r+1).
\end{align*}
For a reduced decomposition $w = s_{r_1}\cdots s_{r_k}\in \Sgp_d$, we have a well-defined element $\tau_w := \tau_{r_1} \cdots \tau_{r_k}$.
It is well-known that $\{\tau_w x_1^{k_1}\dots x_d^{k_d}\mid w\in\Sgp_d,\, k_1,\dots,k_d\in\Z_{\geq 0}\}$ is a basis of $\NH_d$.
In particular, we identify $\CX_d$ as a subalgebra of $\NH_d$.

\begin{Theorem}~\cite[Theorem 2.9]{KL}\label{TNHCenter}
	The center of $\NH_d$ is equal to $\Lambda_d$.
\end{Theorem}

We consider $\CX_d$ as an $\NH_d$-module with $\tau_w$ acting by $\partial_w$ and $x_r$ acting by multiplication with $x_r$.
Then we have the following easy to check and well-known properties:

\begin{Lemma}\label{LPartialPoly}
	Let $d\in \Z_{\geq 0}$, $f\in\CX_d$, and $w\in \Sgp_d$.
	Then in $\NH_d$, we have
	$$\tau_{w_0} f \tau_w = \tau_{w_0}\partial_{w^{-1}}(f) \quad \text{and} \quad \tau_w f \tau_{w_0} = \partial_w(f)\tau_{w_0}.$$
\end{Lemma}

Define the following elements of $\NH_d$:
\begin{equation}\label{ENHIdem}
	e_d := x_{0,d}\tau_{w_{0,d}} \quad \text{and} \quad e_d' := \tau_{w_{0,d}}x_{0,d}.
\end{equation}
Lemmas~\ref{LLambda} and~\ref{LPartialPoly} yield:

\begin{Lemma}\label{LW0X0}
	In $\NH_d$, the elements $e_d$ and $e_d'$ are idempotents.
	Moreover,
	\begin{enumerate}
		\item\label{LW0X0i} $\tau_{w_0}f\tau_{w_0}=0$ for any $f \in \CX_d$ with $\deg f < d(d-1)$, and
		\item\label{LW0X0ii} $\tau_{w_0}x_0\tau_{w_0}=\tau_{w_0}$.
	\end{enumerate}
\end{Lemma}

\subsection{KLR Algebras}\label{SSKLR}

From now on, we set $I := \Z$.
If $i,j \in I$ with $|i-j|=1$ we set $\eps_{i,j}:=j-i\in\{ 1,-1\}$.
We identify $I$ with the set of vertices of the Dynkin diagram of type $\ttA[\infty]$ and denote by $(\ttc_{i,j})_{i,j\in I}$ the corresponding Cartan matrix so that $\ttc_{i,j}=2$ if $i=j$, $\ttc_{i,j}=-1$ if $|i-j| = 1$, and $\ttc_{i,j}=0$ otherwise.
We use the notation $Q_+$, $\Phi_+$, $\height$, etc. introduced in Section~\ref{SIntro}.

For $\theta\in Q_+$ of height $d$, we define $I^\theta:=\{\bi=i_1\cdots i_d\in I^d\mid \alpha_{i_1}+\dots+\alpha_{i_d}=\theta\}$.
Let $\Z((q))\cdot I^\theta:=\bigoplus_{\bi\in I^\theta}\Z((q))\cdot \bi$.
For $1 \leq k \leq t$, suppose $\theta_k\in Q_+$ are such that $\theta_1+\cdots + \theta_t = \theta$, and set $d_k := \height(\theta_k)$.
If $\bi^k \in I^{\theta_k}$, then the concatenation $\bi^1\cdots\bi^t $ is considered as an element of $I^\theta$. Set $\bi^1\cdots\bi^t =: i_1\cdots i_d$.
Then the \emph{quantum shuffle product} is
\begin{equation}\label{EQSh}
	\bi^1 \circ \cdots \circ \bi^t := \sum_{w \in \D^{(d_1, \ldots, d_t)}} q^{-e(w)}w\cdot(\bi^1\cdots\bi^t) \in \Z((q))\cdot I^\theta,
\end{equation}
where $e(w) := \sum_{n < m, \, w(n) > w(m)}\ttc_{i_n,i_m}$.
If $a_k\in \Z((q))\cdot I^{\theta_k}$, we define $a_1 \circ \cdots \circ a_t \in \Z((q))\cdot I^\theta$ extending (\ref{EQSh}) by linearity.

The \emph{KLR algebra}~\cite{KL,R} corresponding to $\theta$ as above is the unital $\k$-algebra $R_\theta$ (with identity denoted $1_\theta$) with generators
$$\{ 1_\bi \mid  \bi \in I^\theta \} \cup \{ y_1, \dots, y_d\} \cup \{\psi_1, \dots, \psi_{d-1}\}$$
and defining relations
\begin{align}
	& y_r y_s = y_s y_r;\label{Rpoly} \tag{R1} \\
	& 1_\bi 1_\bj=\delta_{\bi,\bj}1_\bi \quad \text{and} \quad \sum_{\bi\in I^\theta}1_\bi=1_\theta;\label{Ridem} \tag{R2}\\
	& y_r 1_\bi = 1_\bi y_r \quad \text{and} \quad \psi_r 1_\bi = 1_{s_r\cdot \bi}\psi_r;\label{Rcolor} \tag{R3} \\
	& (\psi_r y_t - y_{s_r(t)} \psi_r)1_\bi = \delta_{i_r,i_{r+1}}(\delta_{t,r+1}-\delta_{t,r})1_\bi;\label{Rdots} \tag{R4} \\
	& \psi_r^2 1_\bi =
	\begin{cases}
		0 & \text{if } i_r=i_{r+1}, \\
		\eps_{i_r,i_{r+1}}({y_r}-{y_{r+1}})1_\bi & \text{if } |i_r-i_{r+1}|=1, \\
		1_\bi & \text{otherwise};
	\end{cases}\label{Rquad} \tag{R5} \\
	& \psi_r \psi_s = \psi_s \psi_r \text{ if } |r-s|>1;\label{Rcomm} \tag{R6} \\
	& (\psi_{r+1} \psi_{r} \psi_{r+1} - \psi_{r}\psi_{r+1}\psi_{r}) 1_\bi =
	\begin{cases}
		\eps_{i_r,i_{r+1}}1_\bi & \text{if } |i_r-i_{r+1}|=1 \text{ and } i_r = i_{r+2}, \\
		0 & \text{otherwise}.
	\end{cases}\label{Rbraid} \tag{R7}
\end{align}
The right-hand sides of relations (\ref{Rdots}) and (\ref{Rbraid}), when they are nonzero, will be referred to as \emph{error terms}.
The algebra $R_\theta$ is graded with $\deg 1_{\bi}=0$; $\deg (y_s)= 2$; $\deg (\psi_{r}1_{\bi})=-\ttc_{i_r,i_{r+1}}$.

We will use the Khovanov-Lauda~\cite{KL} diagrammatic notation for elements of $R_\theta$.
In particular, for $\bi=i_1\cdots i_d\in I^\theta$, $1\leq r<d$ and $1\leq s\leq d$, we denote
$$1_\bi=
\begin{braid}\tikzset{baseline=1em}
	\draw (0,0)--(0,2) node[above]{$i_1$};
	\draw (1,0)--(1,2) node[above]{$i_2$};
	\draw (2,1) node{$\cdots$};
	\draw (2,2) node[above]{$\vphantom{i_1}\cdots\vphantom{i_1}$};
	\draw (3,0)--(3,2) node[above]{$i_d$};
\end{braid}, \quad
1_\bi\psi_r =
\begin{braid}\tikzset{baseline=1em}
	\draw (0,0)--(0,2) node[above]{$\vphantom{i_{r+1}}i_1$};
	\draw (1.25,1) node{$\cdots$};
	\draw (1,2) node[above]{$\vphantom{i_{r+1}}\cdots$};
	\draw (2.5,0)--(2.5,2) node[above]{$i_{r-1}$};
	\draw (4.25,0)--(6,2) node[above]{$i_{r+1}$};
	\draw (6,0)--(4.25,2) node[above]{$\vphantom{i_{r+1}}i_{r}$};
	\draw (7.75,0)--(7.75,2) node[above]{$i_{r+2}$};
	\draw (9,1) node{$\cdots$};
	\draw (9.25,2) node[above]{$\vphantom{i_{r+1}}\cdots$};
	\draw (10.25,0)--(10.25,2) node[above]{$\vphantom{i_{r+1}}i_{d}$};
\end{braid}, \quad
1_\bi y_s=
\begin{braid}\tikzset{baseline=1em}
	\draw (0,0)--(0,2) node[above]{$\vphantom{i_{r+1}}i_1$};
	\draw (1.25,1) node{$\vphantom{i_{r+1}}\cdots$};
	\draw (1,2) node[above]{$\vphantom{i_{r+1}}\cdots$};
	\draw (2.5,0)--(2.5,2) node[above]{$i_{s-1}$};
	\draw (4,0)--(4,2) node[above]{$\vphantom{i_{r+1}}i_{s}$};
	\draw (5.5,0)--(5.5,2) node[above]{$i_{s+1}$};
	\draw (6.75,1) node{$\vphantom{i_{r+1}}\cdots$};
	\draw (7,2) node[above]{$\vphantom{i_{r+1}}\cdots$};
	\draw (8,0)--(8,2) node[above]{$\vphantom{i_{r+1}}i_{d}$};
	\blackdot (4,1);
\end{braid}.$$

For each element $w \in \Sgp_n$, fix a reduced expression $w=s_{r_1}\cdots s_{r_l}$ which determines an element $\psi_w = \psi_{r_1} \cdots \psi_{r_l}$.
This element depends on the reduced expression of $w$.

\begin{Theorem}~\cite[Theorem 2.5]{KL},~\cite[Theorem 3.7]{R}\label{TBasis}
	Let $\theta\in Q_+$ and $d=\height(\theta)$.
	Then the following sets are $\k$-bases of $R_\theta$:
	\begin{align*}
		\{\psi_w y_1^{k_1}\dots y_d^{k_d}1_\bi\mid w\in \Sgp_d,\ k_1,\dots,k_d\in\Z_{\geq 0}, \ \bi\in I^\theta\}, \\
 		\{y_1^{k_1}\dots y_d^{k_d}\psi_w 1_\bi\mid w\in \Sgp_d,\ k_1,\dots,k_d\in\Z_{\geq 0}, \ \bi\in I^\theta\}.
	\end{align*}
\end{Theorem}

We identify the polynomial algebra
\begin{equation}\label{EY}
	\CY_d := \k[y_1, \ldots, y_d]
\end{equation}
with the subalgebra of $R_\theta$ generated by $\{ y_1, \ldots, y_d\}$ according to Theorem~\ref{TBasis}.

The following lemma often simplifies calculations in $R_\theta$.

\begin{Lemma}\label{LPsiProduct}
	Let $i_1,\dots,i_l\in I$ be distinct, $\theta := \alpha_{i_1}+\dots+\alpha_{i_l}$, $\bi = i_1\cdots i_l\in I^{\theta}$, and $w, w' \in \Sgp_l$.
	In $R_\theta$, we have $\psi_{w'} \psi_{w} 1_{\bi} = \psi_{w'w} 1_{\bi}$ unless there are $r, s\in [1, l]$ such that $|i_r-i_s|=1$, $r < s$, $w(r) > w(s)$, and $w'w(r) < w'w(s)$.
\end{Lemma}
\begin{proof}
	As $i_1,\dots,i_l$ are distinct, the braid relation (\ref{Rbraid}) holds without error term in $R_\theta$.
	Moreover, so long as there is no pair $r,s$ as in the statement, the only quadratic relations we need to use are of the form $\psi_t^2 1_\bj = 1_\bj$.
	Therefore $\psi_w \psi_{w'} 1_\bi$ simplifies directly to $\psi_{ww'}1_\bi$ as claimed.
\end{proof}

\subsection{Parabolic subalgebras and divided power idempotents}\label{SSDiv}

Let $\theta_1,\dots,\theta_t\in Q_+$ and set $\theta:=\theta_1+\dots+\theta_t$.
Let
$$1_{\theta_1,\dots,\theta_t} := \sum_{\bi^1\in I^{\theta_1},\dots,\bi^t\in I^{\theta_t}}1_{\bi^1\cdots\bi^t}\in R_\theta.$$
Then we have an algebra embedding
\begin{equation}\label{EEmb}
	\iota_{\theta_1,\dots,\theta_t} : R_{\theta_1}\otimes\dots\otimes R_{\theta_t} \into 1_{\theta_1,\dots,\theta_t}R_{\theta_1+\dots+\theta_t}1_{\theta_1,\dots,\theta_t}
\end{equation}
obtained by horizontal concatenation of the Khovanov-Lauda diagrams.
For $r_1\in R_{\theta_1},\dots, r_t\in R_{\theta_t}$ we often write
$$r_1\circ\dots\circ r_t := \iota_{\theta_1,\dots,\theta_t}(r_1\otimes\dots\otimes r_t).$$
For example,
\begin{equation}\label{EEmbOne}
	1_{\bi^1}\circ\dots\circ 1_{\bi^t} = 1_{\bi^1\cdots \bi^t}\qquad(\bi^1\in I^{\theta_1},\dots,\bi^t\in I^{\theta_t}).
\end{equation}

We fix for the moment $i\in I$, $d\in\Z_{\geq 0}$ and take $\theta = d\alpha_i$.
Then we have an isomorphism $\varphi: \NH_d \iso R_{d\alpha_i}, \ x_r\mapsto y_r, \ \tau_s \mapsto \psi_s$.
Recalling (\ref{ENHIdem}) and Lemma~\ref{LW0X0}, the following element is an idempotent in $R_{d\alpha_i}$:
$$1_{i^{(d)}} := \varphi(e_d').$$

Now let $\theta\in Q_+$ be arbitrary.
We define $I^{\theta}_{\di}$ to be the set of all expressions of the form $i_1^{(d_1)} \cdots i_r^{(d_r)}$ with $d_1,\ldots,d_r\in \Z_{\ge 0}$, $i_1,\dots,i_r\in I$ and $d_1 \alpha_{i_1} + \cdots + d_r \alpha_{i_r} = \theta$.
We refer to such expressions as \emph{divided power words}.
We identify $I^\theta$ with the subset of $I^\theta_\di$ which consists of all divided power words as above with all $d_k=1$.
We use the same notation for concatenation of divided power words as for concatenation of words.
For $\bi=i_1^{(d_1)}\cdots i_r^{(d_r)}\in I^\theta_{\di}$, we define
\begin{equation}\label{EFactorial}
	\bi_{(\pm)}^! := [d_1]_{(\pm)}^! \cdots [d_r]_{(\pm)}^!, \quad \text{and} \quad \hat\bi := i_1^{d_1} \cdots i_r^{d_r} \in I^{\theta},
\end{equation}
and the corresponding \emph{divided power idempotent} is
$$1_{\bi}=1_{i_1^{(d_1)}\cdots i_r^{(d_r)}}:=1_{i_1^{(d_1)}}\circ\dots\circ 1_{i_r^{(d_r)}} \in R_{\theta}.$$
We have the following generalization of (\ref{EEmbOne}):
$$1_{\bi^1}\circ\dots\circ  1_{\bi^t}=1_{\bi^1\cdots \bi^t}\qquad(\bi^1\in I^{\theta_1}_{\di},\dots,\bi^t\in I^{\theta_t}_{\di}).$$

\begin{Lemma}\label{LDivPower}
	In the algebra $R_{d\alpha_i}$, if $r_1+\dots+r_t=d$ then $1_{i^{(r_1)} \cdots i^{(r_t)}}\psi_{w_{0,d}}=\psi_{w_{0,d}}$ and $1_{i^{(r_1)} \cdots i^{(r_t)}}1_{i^{(d)}}=1_{i^{(d)}}$.
\end{Lemma}
\begin{proof}
	Write $\psi_{w_{0,d}}=(\psi_{w_{0,r_1}}\circ\dots\circ \psi_{w_{0,r_t}})\psi_u$ for some $u\in\Sgp_d$ and use Lemma~\ref{LW0X0}.
\end{proof}

To be used as part of the Khovanov-Lauda diagrammatics, we denote
$$\psi_{w_{0,d}} =:
\begin{braid}\tikzset{baseline=0.25em}
	\braidbox{0}{2}{0}{1}{$w_0$}
\end{braid}, \quad y_{0,d} =:
\begin{braid}\tikzset{baseline=0.25em}
	\braidbox{0}{2}{0}{1}{$y_0$}
\end{braid}, \quad 1_{i^{(d)}} =
\begin{braid}\tikzset{baseline=1.5em}
	\draw (0,2)--(0,3) node[above]{$i$};
	\draw (1,2.5) node{$\cdots$};
	\draw (1,3) node[above]{$\cdots$};
	\draw (2,2)--(2,3) node[above]{$i$};
	\braidbox{0}{2}{0}{1}{$y_0$}
	\braidbox{0}{2}{1}{2}{$w_0$}
\end{braid} =:
\begin{braid}\tikzset{baseline=0.25em}
	\braidbox{0}{2}{0}{1.1}{$i^d$}
\end{braid}.$$
For example, if $d=3$, we have
$$1_{i^3}\psi_{w_0}=
\begin{braid}\tikzset{baseline=1em}
	\draw(0,1)--(0,2) node[above]{$i$};
	\draw(1,1)--(1,2) node[above]{$i$};
	\draw(2,1)--(2,2) node[above]{$i$};
	\braidbox{0}{2}{0}{1}{$w_0$}
\end{braid} =
\begin{braid}\tikzset{baseline=1em}
	\draw (0,0)--(2,2) node[above]{$i$};
	\draw (1,0)--(0,1)--(1,2) node[above]{$i$};
	\draw (2,0)--(0,2) node[above]{$i$};
\end{braid}, \quad 1_{i^3}y_0 =
\begin{braid}\tikzset{baseline=1em}
	\draw (0,1)--(0,2) node[above]{$i$};
	\draw (1,1)--(1,2) node[above]{$i$};
	\draw (2,1)--(2,2) node[above]{$i$};
	\braidbox{0}{2}{0}{1}{$y_0$}
\end{braid} =
\begin{braid}\tikzset{baseline=1em}
	\draw (0,0)--(0,2) node[above]{$i$};
	\draw (1,0)--(1,2) node[above]{$i$};
	\draw (2,0)--(2,2) node[above]{$i$};
	\blackdot (1,1);
	\blackdot (2,0.7);
	\blackdot (2,1.3);
\end{braid}, \quad 1_{i^{(3)}} =
\begin{braid}\tikzset{baseline=1em}
	\draw (0,0)--(2,2) node[above]{$i$};
	\draw (1,0)--(0,1)--(1,2) node[above]{$i$};
	\draw (2,0)--(0,2) node[above]{$i$};
	\blackdot (0.75,0.25);
	\blackdot (1.3,0.7);
	\blackdot (1.7,0.3);
\end{braid}.$$
More generally, we denote
$$1_{i_1^{(d_1)}\cdots\, i_r^{(d_r)}} =:
\begin{braid}\tikzset{baseline=0.5em}
	\braidbox{0}{3}{0}{1.5}{$i_1^{d_1}$}
	\draw (4.5,0.75) node{\normalsize $\cdots$};
	\braidbox{6}{9}{0}{1.5}{$i_r^{d_r}$}
\end{braid}.$$

\subsection{Modules over $R_\theta$}\label{SSModules}

Let $\theta\in Q_+$.
We denote by $\LMod{R_\theta}$ the category of graded left $R_\theta$-modules.
The morphisms in this category are all homogeneous degree zero  $R_\theta$-homomorphisms, which we denote $\hom_{R_\theta}(-,-)$.
For $V\in\LMod{R_\theta}$, let $q^d V$ denote its grading shift by $d$, so if $V_m$ is the degree $m$ component of $V$, then $(q^dV)_m= V_{m-d}$.
More generally, for a Laurent series $a = a(q) = \sum_{d}a_dq^d\in\Z((q))$ with non-negative coefficients, we set $aV := \bigoplus_d(q^d V)^{\oplus a_d}$.
For $U,V\in \LMod{R_\theta}$, we set $\Hom_{R_\theta}(U, V) := \bigoplus_{d \in \Z} \Hom_{R_\theta}(U, V)_d$, where
$$\Hom_{R_\theta}(U, V)_d := \hom_{R_\theta}(q^d U, V) = \hom_{R_\theta}(U, q^{-d}V).$$
We define $\Ext^m_{R_\theta}(U,V)$ and $\End_{R_\theta}(U)$ similarly from $\operatorname{ext}^m_{R_\theta}(U,V)$ and $\operatorname{end}_{R_\theta}(U)$.

For a free $\k$-module $V$ of finite rank we denote the rank of $V$ by $\dim V$.
A graded $\k$-module $V = \bigoplus_{m\in \Z}V_m$ is called \emph{Laurentian} if the graded components $V_m$ are free of finite rank for all $m\in \Z$ and $V_m = 0$ for $m \ll 0$.
For example $R_\theta$ itself is Laurentian by Theorem~\ref{TBasis}.
The \emph{graded rank} of a Laurentian $\k$-module $V$ is
$$\qdim V:=\sum_{m\in \Z}(\dim V_m) q^m\in\Z((q)).$$

Recall that the ground ring $\k$ is assumed to be a PID.
The following standard result often allows us to reduce to the case where $\k$ is a field.

\begin{Lemma}\label{LGDim}
	If $\varphi: V\to W$ is a degree $0$ homomorphism of Laurentian $\k$-modules such that the induced map $\overline{\varphi}:V/JV \to W/JW$ is an isomorphism of $\k/J$-vector spaces for every maximal ideal $J$, then $\varphi$ is an isomorphism.
\end{Lemma}

We say that an $R_\theta$-module $V$ is \emph{Laurentian} if it is so as a $\k$-module.
Recalling (\ref{EFactorial}), for a Laurentian $R_\theta$-module $V$ and $\bi \in I_\di^\theta$, by~\cite[\S2.5]{KL}, we have
\begin{equation}\label{EUndivided}
	\qdim (1_{\bi} V) = \frac{1}{\bi_+^!} \qdim(1_{\hat\bi} V),
\end{equation}
which explains the usage of the term ``divided power word'' for $\bi \in I_\di^\theta$.
If $V$ is a Laurentian $R_\theta$-module then each $1_\bi V$ is a Laurentian $\k$-module, and so we can define the \emph{formal character} of $V$ as follows:
$$\qch V := \sum_{\bi\in I^\theta}(\qdim 1_\bi V)\cdot\bi \in \Z((q))\cdot I^\theta.$$
Note that $\qch(q^dV)=q^d\qch(V)$.

For $\theta_1,\ldots,\theta_t\in Q_+$ and $\theta := \theta_1+\cdots+\theta_t$, recalling (\ref{EEmb}), we have a functor
$$\Ind_{\theta_1,\ldots,\theta_t} = R_\theta 1_{\theta_1,\ldots,\theta_t}\otimes_{R_{\theta_1}\otimes\cdots\otimes R_{\theta_t}}-:\LMod{(R_{\theta_1}\otimes\cdots\otimes R_{\theta_t})}\to\LMod{R_{\theta}}.$$
For $V_1\in\LMod{R_{\theta_1}}, \ldots, V_t\in\LMod{R_{\theta_t}}$, we denote by $V_1\boxtimes \cdots\boxtimes V_t$ the $\k$-module $V_1\otimes \cdots\otimes V_t$, considered naturally as an $(R_{\theta_1}\otimes\cdots\otimes R_{\theta_t})$-module, and set
$$V_1\circ\cdots\circ V_t := \Ind_{\theta_1,\ldots,\theta_t} V_1\boxtimes \cdots\boxtimes V_t.$$
By Theorem~\ref{TBasis}, setting $d_k:=\height(\theta_k)$, we have
\begin{equation}\label{EShuffleMod}
	V_1 \circ \cdots \circ V_t = \bigoplus_{w\in\D^{(d_1,\dots,d_t)}}\psi_w\otimes V_1\otimes \cdots\otimes V_t.
\end{equation}
If $V_1, \ldots, V_t$ are Laurentian, then by~\cite[Lemma 2.20]{KL}, recalling (\ref{EQSh}), we have
\begin{equation}\label{EShuffle}
	\qch(V_1\circ \cdots \circ V_t) = \qch(V_1)\circ \cdots \circ\qch(V_t).
\end{equation}
For $v_1\in V_1,\dots, v_t\in V_t$, we denote
$$v_1\circ \dots\circ v_t := 1_{\theta_1,\dots,\theta_t}\otimes v_1\otimes\dots\otimes v_t\in V_1\circ \dots\circ V_t.$$

If $\bi^1\in I_\di^{\theta_1}, \ldots, \bi^t\in I_\di^{\theta_t}$, it is easy to check that
\begin{equation}\label{EIndIdemp}
	R_{\theta_1}1_{\bi^1}\circ \cdots \circ R_{\theta_t} 1_{\bi^t}\cong R_{\theta}1_{\bi^1\cdots\bi^t}.
\end{equation}
Since $R_\theta 1_{\theta_1,\dots,\theta_t}$ is a free right $R_{\theta_1}\otimes\dots\otimes R_{\theta_t}$-module of finite rank by Theorem~\ref{TBasis}, we get the following well-known properties:

\begin{Proposition}\label{PCircExact}
	The functor $\Ind_{\theta_1,\dots,\theta_t}$ is exact and sends finitely generated projectives to finitely generated projectives.
\end{Proposition}

Let again $\theta_1, \ldots, \theta_t \in Q_+$ and $\theta = \theta_1 + \cdots + \theta_t$.
Suppose $(C_\bullet^k, d_k)$ is a chain complexes of $R_{\theta_k}$-modules for $1 \leq k \leq t$.
Let $(C_\bullet^1 \circ \cdots \circ C_\bullet^t)_n := \bigoplus_{p_1 + \cdots + p_t = n} C_{p_1}^1 \circ \cdots \circ C_{p_t}^t$, and for $x_1 \in C_{p_1}^1, \ldots, x_t \in C_{p_t}^t$, define
$$d(x_1\circ \cdots \circ x_t) := \sum_{k=1}^t (-1)^{p_{k+1}+ \cdots + p_t} x_1 \circ \cdots \circ x_{k-1} \circ d_k(x_k) \circ x_{k+1} \circ \cdots \circ x_t.$$
Then $(C_\bullet^1 \circ \cdots \circ C_\bullet^t, d)$ is a chain complex of $R_\theta$-modules.
Proposition~\ref{PCircExact} and~\cite[Lemma 2.7.3]{We} immediately imply the following.

\begin{Lemma}\label{LCircRes}
	If $C_\bullet^k$ is a projective resolution of $M_k \in \LMod{R_{\theta_k}}$ for $1 \leq k \leq t$, then $C_\bullet^1 \circ \cdots \circ C_\bullet^t$ is a projective resolution of $M_1\circ \cdots \circ M_t \in \LMod{R_\theta}$.
\end{Lemma}

\subsection{Standard modules}\label{SSStdMod}

The algebra $R_\theta$ is affine quasihereditary in the sense of~\cite{Kdonkin}.
In particular, it comes with an important class of \emph{standard modules}, which we now describe explicitly following~\cite{BKM}.

Fix $\beta = \alpha_i + \cdots + \alpha_j \in \Phi_+$ of height $l := j-i+1$, and set $\bi_\beta := i\,(i+1)\cdots j\in I^\beta$.
We define the $R_\beta$-module $\Delta(\beta)$ to be a cyclic $R_\beta$-module generated by a vector $v_\beta$ of degree $0$ with defining relations
\begin{itemize}
	\item $1_\bi v_\beta = \delta_{\bi,\bi_\beta}v_\beta$ for all  $\bi\in I^\beta$;
	\item $\psi_r v_\beta = 0$ for all $1\leq r < l$;
	\item $y_r v_\beta = y_s v_\beta$ for all $1\leq r, s \leq l$.
\end{itemize}
The module $\Delta(\beta)$ can be considered as an $(R_\beta, \k[x])$-bimodule with the right action given by $v_\beta x := y_1 v_\beta.$
Diagrammatically, we represent
$$v_\beta =
\begin{braid}\tikzset{baseline=0.5em}
	\draw (1.5,0) node[below]{$\bi_\beta$}--(1.5,1)--(0,2);
	\draw (1.5,0)--(1.5,1)--(1,2);
	\draw (1.5,0)--(1.5,1)--(3,2);
	\draw (2,1.9) node{$\cdots$};
\end{braid},$$
$$v_\beta x =
\begin{braid}\tikzset{baseline=0.5em}
	\draw (1.5,0) node[below]{$\bi_\beta$}--(1.5,1)--(0,2);
	\draw (1.5,0)--(1.5,1)--(1,2);
	\draw (1.5,0)--(1.5,1)--(3,2);
	\draw (2,1.9) node{$\cdots$};
	\blackdot (1.5,0.5);
\end{braid} = y_1 v_\beta =
\begin{braid}\tikzset{baseline=0.5em}
	\draw (1.5,0) node[below]{$\bi_\beta$}--(1.5,1)--(0,2);
	\draw (1.5,0)--(1.5,1)--(1,2);
	\draw (1.5,0)--(1.5,1)--(3,2);
	\draw (2,1.9) node{$\cdots$};
	\blackdot (0.75,1.5);
\end{braid} = \dots = y_l v_\beta =
\begin{braid}\tikzset{baseline=0.5em}
	\draw (1.5,0) node[below]{$\bi_\beta$}--(1.5,1)--(0,2);
	\draw (1.5,0)--(1.5,1)--(1,2);
	\draw (1.5,0)--(1.5,1)--(3,2);
	\draw (2,1.9) node{$\cdots$};
	\blackdot (2.25,1.5);
\end{braid},$$
$$\psi_1 v_\beta =
\begin{braid}\tikzset{baseline=1em}
	\draw (1.5,0) node[below]{$\bi_\beta$}--(1.5,1)--(0,2)--(1,3);
	\draw (1.5,0)--(1.5,1)--(1,2)--(0,3);
	\draw (1.5,0)--(1.5,1)--(3,2)--(3,3);
	\draw (2,2) node{$\cdots$};
\end{braid} = \cdots = \psi_{l-1} v_\beta =
\begin{braid}\tikzset{baseline=1em}
	\draw (1.5,0) node[below]{$\bi_\beta$}--(1.5,1)--(0,2)--(0,3);
	\draw (1.5,0)--(1.5,1)--(2,2)--(3,3);
	\draw (1.5,0)--(1.5,1)--(3,2)--(2,3);
	\draw (1,2) node{$\cdots$};
\end{braid}=0.$$

The following lemma is easy to check.

\begin{Lemma}\label{LDeltaRootBasis}
	Let $\beta\in \Phi_+$.
	Then there is an isomorphism of right $\k[x]$-modules $\k[x] \to \Delta(\beta), \ 1 \mapsto v_\beta$.
\end{Lemma}

For $m\in \Z_{\geq 0}$ and $\beta \in \Phi_+$, the $R_{m\beta}$-module $\Delta(\beta)^{\circ m}$ is cyclicly generated by $v_\beta^{\circ m}$.
As explained in~\cite[\S3.2]{BKM}, $\NH_m$ acts on $\Delta(\beta)^{\circ m}$ on the right so that
$$v_\beta^{\circ m}x_r = v_\beta^{\circ (r-1)}\circ (v_\beta x)\circ v_\beta^{\circ (m-r)}, \quad v_\beta^{\circ m}\tau_s = v_\beta^{\circ (s-1)}\circ (\psi_{w_{l,l}} (v_\beta \circ v_\beta))\circ v_\beta^{\circ (m-s-1)},$$
where $w_{l,l}$ is the longest element of $\D^{(l,l)}$.
Diagrammatically, we represent
$$v_\beta^{\circ m} =
\begin{braid}\tikzset{baseline=0.5em}
	\draw (1.5,0) node[below]{$\bi_\beta$}--(1.5,1)--(0,2);
	\draw (1.5,0)--(1.5,1)--(1,2);
	\draw (1.5,0)--(1.5,1)--(3,2);
	\draw (2,1.9) node{$\cdots$};

	\draw (5.5,0) node[below]{$\bi_\beta$}--(5.5,1)--(4,2);
	\draw (5.5,0)--(5.5,1)--(5,2);
	\draw (5.5,0)--(5.5,1)--(7,2);
	\draw (6,1.9) node{$\cdots$};

	\draw (8,1) node{\Large$\cdots$};

	\draw (10.5,0) node[below]{$\bi_\beta$}--(10.5,1)--(9,2);
	\draw (10.5,0)--(10.5,1)--(10,2);
	\draw (10.5,0)--(10.5,1)--(12,2);
	\draw (11,1.9) node{$\cdots$};
\end{braid},$$
$$(v_\beta\circ v_\beta)\tau_1 =
\begin{braid}\tikzset{baseline=1.5em}
	\draw (1.5,0) node[below]{$\bi_\beta$}--(5.5,2)--(5.5,3)--(4,4);
	\draw (1.5,0)--(5.5,2)--(5.5,3)--(5,4);
	\draw (1.5,0)--(5.5,2)--(5.5,3)--(7,4);
	\draw (6,3.9) node{$\cdots$};

	\draw (5.5,0) node[below]{$\bi_\beta$}--(1.5,2)--(1.5,3)--(0,4);
	\draw (5.5,0)--(1.5,2)--(1.5,3)--(1,4);
	\draw (5.5,0)--(1.5,2)--(1.5,3)--(3,4);
	\draw (2,3.9) node{$\cdots$};
\end{braid} = \psi_{w_{l,l}}(v_\beta \circ v_\beta) =
\begin{braid}\tikzset{baseline=1.5em}
	\draw (1.5,0) node[below]{$\bi_\beta$}--(1.5,1)--(0,2)--(4,4);
	\draw (1.5,0)--(1.5,1)--(1,2)--(5,4);
	\draw (1.5,0)--(1.5,1)--(3,2)--(7,4);
	\draw (2,1.9) node{$\cdots$};
	\draw (2.2,3.9) node{$\cdots$};

	\draw (5.5,0) node[below]{$\bi_\beta$}--(5.5,1)--(4,2)--(0,4);
	\draw (5.5,0)--(5.5,1)--(5,2)--(1,4);
	\draw (5.5,0)--(5.5,1)--(7,2)--(3,4);
	\draw (6,1.9) node{$\cdots$};
	\draw (5.8,3.9) node{$\cdots$};
\end{braid}.$$

Let $\leq$ be the lexicographic total order on $\Phi_+$, i.e. for $\beta = \alpha_i + \cdots + \alpha_j\in \Phi_+$ and $\beta' = \alpha_{i'} + \cdots + \alpha_{j'}\in \Phi_+$, we have $\beta < \beta'$ if and only if either $i < i'$ or $i = i'$ and $j < j'$.
Given $\theta\in Q_+$, a \emph{Kostant partition} of $\theta$ is a sequence $\pi = (\beta_1^{m_1}, \ldots, \beta_t^{m_t})$ such that $m_1, \ldots, m_t\in \Z_{>0}$, $\beta_1 > \cdots > \beta_t$ are positive roots, and $m_1\beta_1 + \cdots + m_t\beta_t = \theta$.
We denote by $\KP(\theta)$ the set of all Kostant partitions of $\theta$.
For $\pi = (\beta_1^{m_1}, \ldots, \beta_t^{m_t})\in \KP(\theta)$,
\begin{equation}\label{EDeltaHat}
	\hat{\Delta}(\pi) := q^{\binom{m_1}{2}+\cdots+\binom{m_t}{2}} \Delta(\beta_1)^{\circ m_1}\circ \cdots \circ \Delta(\beta_t)^{\circ m_t}
\end{equation}
can now be considered as an $(R_\theta, \NH_{m_1}\otimes \cdots \otimes \NH_{m_t})$-bimodule.
Recalling (\ref{ENHIdem}), we define the corresponding \emph{standard module} as
$$\Delta(\pi) := \hat{\Delta}(\pi)(e_{m_1}\otimes \cdots \otimes e_{m_t}).$$
Setting
\begin{equation}\label{EPartSym}
	\Lambda_\pi := \Lambda_{m_1, \ldots, m_t},
\end{equation}
 by Theorem~\ref{TNHCenter}, $\Delta(\pi)$ is naturally an $(R_\theta, \Lambda_\pi)$-bimodule.
In fact, by~\cite[Theorem 2.17]{KS}, the bimodule structure yields the isomorphism
\begin{equation}\label{EDeltaEnd}
	\End_{R_\theta}(\Delta(\pi))^{\op}\cong \Lambda_\pi.
\end{equation}

The module $\Delta(\pi)$ is cyclic as a left $R_\theta$-module with \emph{standard generator}
\begin{equation}\label{EDeltaGen}
	v_\pi := (v_\beta^{\circ m_1}\circ \cdots \circ v_\beta^{\circ m_t})(e_{m_1}\otimes \cdots \otimes e_{m_t}).
\end{equation}
Noting that $\Delta(\pi) = \Delta(\beta_1^{m_1})\circ \cdots \circ \Delta(\beta_1^{m_t})$, by~\cite[Lemma 3.10]{BKM}, we have an isomorphism of $R_\theta$-modules
\begin{equation}\label{EDeltaIso}
	\hat{\Delta}(\pi) \cong [m_1]_+^! \cdots [m_t]_+^! \Delta(\pi).
\end{equation}

If $\k$ is a field, the modules $\{\Delta(\pi)\mid \pi\in \KP(\theta)\}$ are the standard modules for an affine quasihereditary structure on the algebra $R_\theta$, see~\cite{BKM,Kdonkin}.
If $\k=\Z$ or $\Z_p$, they can be thought of as integral forms of the standard modules, see~\cite[\S4]{KS}.

For $\rho\in \KP(\theta)$, suppose we have a projective resolution of $\Delta(\rho)$ of the form $P_\bullet^\rho = \cdots \overset{d_1}{\to} P_1^\rho \overset{d_0}{\to} P_0^\rho \overset{\eps_\rho}{\to} \Delta(\rho)$ with $P_n^\rho = \bigoplus_{x\in X_n} q^{s_x}R_\theta 1_x$ for some index set $X_n$, integers $s_x$, and idempotents $1_x$.
The map $d_n: P_{n+1}^\rho \to P_n^\rho$ can be described as right multiplication by an $X_{n+1}\times X_n$ matrix $D_n = (d_n^{y,x})$ for some $d_n^{y,x} \in 1_y R_\theta 1_x$.
Using the isomorphism $\Hom_{R_\theta}(q^n R_\theta 1_x, \Delta(\sigma)) \iso q^{-n}1_x \Delta(\sigma)$ and recalling (\ref{EDeltaEnd}), we obtain:

\begin{Lemma}\label{LCochain}
	There is an isomorphism of complexes of (right) $\Lambda_\sigma$-modules
	$$\Hom^\bullet_{R_\theta}(P_\bullet^\rho, \Delta(\sigma))\cong T_\bullet^{P_\bullet^\rho}(\Delta(\sigma)),$$
	where $T_\bullet^{P_\bullet^\rho}(\sigma) = \cdots \overset{d^1}{\from} T_1^{P_\bullet^\rho}(\sigma) \overset{d^0}{\from} T_0^{P_\bullet^\rho}(\sigma)$ with $T_n^{P_\bullet^\rho}(\sigma)=\bigoplus_{x\in X_n}q^{-s_x}1_x \Delta(\sigma)$ and $d^n$ given by left multiplication with the $X_{n+1}\times X_n$ matrix $D_n$.
\end{Lemma}

This yields an isomorphism $\CE_\theta(\rho, \sigma) \cong H(T_\bullet^{P_\bullet^\rho}(\sigma))$ of $\Lambda_\sigma$-modules.
One can also use the resolutions $P_\bullet^\rho$, $P_\bullet^\sigma$, and $P_\bullet^\tau$ to describe the composition map $\CE_\theta(\sigma, \tau) \otimes \CE_\theta(\rho, \sigma) \to \CE_\theta(\rho, \tau)$.
Indeed, let $\Hom_{R_\theta}^m(P_\bullet^\rho, P_\bullet^\sigma)$ denote the homological degree $m$ homomorphisms.
Then $\Hom_{R_\theta}^\bullet(P_\bullet^\rho, P_\bullet^\sigma)$ is a complex with respect to the differential $\delta$ given by
\begin{equation}\label{EDGDiff}
	\delta(\varphi) := d\varphi - (-1)^m\varphi d
\end{equation}
for $\varphi\in \Hom_{R_\theta}^m(P_\bullet^\rho, P_\bullet^\sigma)$.
We have an isomorphism
\begin{equation}\label{EExtIso}
	H(\Hom_{R_\theta}^\bullet(P_\bullet^\rho, P_\bullet^\sigma)) \cong \CE_\theta(\rho, \sigma)
\end{equation}
induced by the maps
\begin{equation}\label{EExtInduced}
	\Hom_{R_\theta}^m(P_\bullet^\rho, P_\bullet^\sigma) \to \Hom_{R_\theta}(P_m^\rho, \Delta(\sigma)), \ \varphi \mapsto (-1)^{\frac{m(m+1)}{2}}\eps_\sigma(\varphi|_{P_m^\rho}).
\end{equation}
Now the composition map $\CE_\theta^n(\sigma, \tau) \otimes \CE_\theta^m(\rho, \sigma) \to \CE_\theta^{n+m}(\rho, \tau)$ is induced from the composition of homomorphisms $\Hom_{R_\theta}^n(P_\bullet^\sigma, P_\bullet^\tau) \otimes \Hom_{R_\theta}^m(P_\bullet^\rho, P_\bullet^\sigma) \to \Hom_{R_\theta}^{n+m}(P_\bullet^\rho, P_\bullet^\tau)$.

\section{The case where $\theta$ is a root}\label{SRoot}

Throughout the section, $\theta=\alpha_a+\dots+\alpha_{b+1}$ (with $a\leq b+1$) is a positive root of height $l=b+2-a$.
There is a bijection from $\KP(\theta)$ to the set of all subsets of $[a,b]$
$$\pi = (\pi_1, \ldots, \pi_u) \mapsto C_\pi:=\{ \max(\supp \pi_2), \max(\supp \pi_3), \ldots, \max(\supp \pi_u)\},$$
where, for a root $\alpha = \alpha_i + \dots+ \alpha_j$, we let $\supp \alpha := [i,j]$.
For $\pi, \tau \in \KP(\theta)$, if $C_\tau \supseteq C_\pi$, we say that $\tau$ is a \emph{refinement} of $\pi$ and write $\tau \supseteq \pi$.
If, in addition, $|C_\tau \setminus C_\pi| = n$, we write $\tau \supseteq_n \pi$ and say that $\tau$ is an \emph{$n$-refinement} of $\pi$.
If $C_\tau \setminus C_\pi = \{ i\}$ for some $i \in [a,b]$, we write $\refn^i(\pi) := \tau$.
For example, we have $\refn^i((\theta)) = (\alpha_{i+1} + \cdots + \alpha_{b+1}, \alpha_a + \cdots + \alpha_i)$.

If $\tau = (\tau_1, \ldots, \tau_t) \in \KP(\theta)$, the elements of
\begin{equation}\label{ERootShuffles}
	\D^\tau := \D^{(\height(\tau_1), \ldots, \height(\tau_t))}
\end{equation}
are called \emph{$\tau$-shuffles}.
Set
$$d_v := \height(\tau_1) + \cdots + \height(\tau_v)\qquad(0 \leq v \leq t).$$
We say that integers $r,s \in [1,l]$ are \emph{$\tau$-equivalent} if there is some $v \in [1,t]$ with $d_{v-1} < r,s \leq d_v$.
Recalling (\ref{EPartSym}) and (\ref{EDeltaEnd}), we have
$$\Lambda_\tau=\k[x_1,\ldots,x_t] \cong \End_{R_\theta}(\Delta(\tau))^\op.$$
We have a surjection
\begin{equation}\label{EYToLambda}
	\ttp_\tau: \CY_l \onto \Lambda_\tau, \ y_r \mapsto x_v \text{ if } d_{v-1} < r \leq d_v,
\end{equation}
so that $\ttp_\tau(y_r)=\ttp_\tau(y_s)$ if and only if $r$ and $s$ are $\tau$-equivalent.
If, in addition, $\tau \supseteq \pi = (\pi_1, \ldots, \pi_u)$, then $\ttp_\pi$ factors as $\ttp_\pi= \ttp^\tau_\pi  \ttp_\tau$, where the surjection $\ttp^\tau_\pi$ is defined as follows: for each $m\in [0,u]$, we have $\height(\pi_1) + \cdots + \height(\pi_m) = d_{v_m}$ for some $v_m \in [0,t]$; now
\begin{equation}\label{ELambdaToLambda}
	\ttp_\pi^\tau: \Lambda_\tau \onto \Lambda_\pi, \ x_r \mapsto x_m \text{ if } v_{m-1} < r \leq v_m.
\end{equation}

\subsection{The resolution $S_\bullet^\rho$}\label{SSRootRes}

In this subsection, we fix $\rho = (\rho_1, \ldots, \rho_t)\in \KP(\theta)$.
For $\alpha = \alpha_i + \alpha_{i+1} + \cdots + \alpha_j\in \Phi_+$, let
$$\bj_\alpha := i(i+1) \cdots j\in I^\alpha,\quad  e_\alpha := 1_{\bj_\alpha}\in R_\alpha.$$
Then set
\begin{equation}\label{ERootIdem}
	\bj_\rho := \bj_{\rho_1} \cdots \bj_{\rho_t}\in I^\theta, \quad e_\rho := 1_{\bj_\rho}\in R_\theta.
\end{equation}
For $\pi,\tau\in \KP(\theta)$, let $w(\tau,\pi)\in \Sgp_l$ be the unique permutation with
\begin{equation}\label{ERootPerm}
	w(\tau,\pi)\cdot \bj_\pi=\bj_\tau,
\end{equation}
so that $e_\tau\psi_{w(\tau,\pi)} = \psi_{w(\tau,\pi)}e_\pi$.
If $\tau=\refn^i(\pi)$, we set
$$s(\tau, \pi) := (-1)^{|C_\pi \cap [a,i)|}.$$

For $n\in \Z_{\geq 0}$, we set
$$S_n^{\rho} := \bigoplus_{\pi \supseteq_n \rho}q^nR_\theta e_{\pi}.$$
The boundary map $S_{n+1}^\rho \to S_{n}^\rho$ is defined to be right multiplication with the matrix
\begin{equation}\label{ERootDiff}
	d_n := \left(d_n^{\tau, \pi}\right)_{\subalign{\tau &\supseteq_{n+1} \rho \\ \pi &\supseteq_n \rho}}, \text{ where } d_n^{\tau, \pi} :=
	\begin{cases}
		s(\tau, \pi)e_\tau\psi_{w(\tau, \pi)}e_\pi & \text{if } \tau \supseteq \pi, \\
		0 & \text{otherwise}.
	\end{cases}
\end{equation}
We define the augmentation map by
$$\eps_\rho:S_0^\rho=R_\theta e_\rho \to \Delta(\rho), \ he_{\rho}\mapsto hv_\rho,$$
where $v_\rho$ is the standard generator for $\Delta(\rho)$, see (\ref{EDeltaGen}).

\begin{Lemma}\label{LSIsRes}
	The following is a projective resolution of $\Delta(\rho)$:
	$$0 \longrightarrow S_{l-t}^\rho \longrightarrow \cdots \longrightarrow S_{n+1}^\rho \stackrel{d_n}{\longrightarrow} S_n^\rho \longrightarrow \cdots \longrightarrow S_0^\rho \stackrel{\eps_\rho}{\longrightarrow} \Delta(\rho) \longrightarrow 0.$$
\end{Lemma}
\begin{proof}
	Given a complex $C_\bullet$, we denote by $\overline{C}_\bullet$ the same complex but with all the boundary maps negated.
	If $v\in[1,t]$, we let
	$$\widehat{S}_\bullet^{(\rho_v)} =
	\begin{cases}
		\overline{S}_\bullet^{(\rho_v)} & \text{if $t-v$ is odd}, \\
		S_\bullet^{(\rho_v)} & \text{if $t-v$ is even}.
	\end{cases}$$
	Using (\ref{EIndIdemp}) and the fact that the resolutions $\overline{S}_\bullet^{(\rho_v)}$ and $S_\bullet^{(\rho_v)}$ are isomorphic, it is easy to note that
	$$S_\bullet^{\rho}=\widehat{S}_\bullet^{(\rho_1)}\circ \dots \circ \widehat{S}_\bullet^{(\rho_t)}\cong S_\bullet^{(\rho_1)}\circ \dots \circ S_\bullet^{(\rho_t)},$$
	so by Lemma~\ref{LCircRes}, we have reduced to the case $t = 1$, i.e. $\rho = (\rho_1) = (\theta)$.

	To complete the proof, we show that $S_\bullet^{(\theta)}\cong P_\bullet$, where $P_\bullet$ is a resolution of $\Delta(\theta)$ constructed in~\cite[\S4.5]{BKM} (see also~\cite{BKS}) and which we now recall.
	For $\pi\in \KP(\theta)$, put
	$$\bi^\pi := a^{\delta_{a\not\in C_\pi}}(a+1)^{\delta_{a+1\not\in C_\pi}}\cdots b^{\delta_{b\not\in C_\pi}}(b+1)b^{\delta_{b\in C_\pi}}\cdots (a+1)^{\delta_{a+1\in C_\pi}}a^{\delta_{a\in C_\pi}} \in I^\theta.$$
	For $n\in\Z_{\geq 0}$, set $P_n := \bigoplus_{\pi\supseteq_n (\theta)} q^n R_\theta 1_{\bi^\pi}$.
	If $\pi, \tau \in \KP(\theta)$ with $\tau = \refn^i(\pi)$, let $u(\tau,\pi)\in\Sgp_l$ be determined from $u(\tau,\pi)\cdot\bi^\pi=\bi^\tau$, and define the matrix
	$$\partial_n := (\partial_n^{\tau,\pi})_{\subalign{\tau &\supseteq_{n+1}(\theta) \\ \pi &\supseteq_n (\theta)}}, \text{ where } \partial_n^{\tau,\pi} :=
	\begin{cases}
		s(\tau, \pi)1_{\bi^\tau}\psi_{u(\tau,\pi)}1_{\bi^\pi} & \text{if } \tau \supseteq \pi, \\
		0 & \text{otherwise}.
	\end{cases}$$
	Right multiplication with $\partial_n$ defines a map $P_{n+1}\to P_n$.
	By~\cite[Theorem 4.12]{BKM} (see also~\cite[Theorem A]{BKS}), noting that $P_0 = S_0^{(\theta)}$, we have that
	$$0 \longrightarrow P_{b+1-a} \longrightarrow \cdots \longrightarrow P_{n+1} \stackrel{\partial_n}{\longrightarrow} P_n \longrightarrow \cdots \longrightarrow P_0 \stackrel{\eps_{(\theta)}}{\longrightarrow} \Delta(\theta) \longrightarrow 0$$
	is a projective resolution of $\Delta(\theta)$.

	For $\pi\in \KP(\theta)$, let $w(\pi)\in \Sgp_l$ be the unique permutation with $w(\pi)\cdot\bi^\pi=\bj_\pi$ so that $e_\pi\psi_{w(\pi)} = \psi_{w(\pi)}1_{\bi^\pi}$.
	We have the map
	\begin{align*}
		S^{(\theta)}_n = \bigoplus_{\pi \supseteq_n (\theta)}q^n R_\theta e_\pi &\to \bigoplus_{\pi \supseteq_n (\theta)}q^n R_\theta 1_{\bi^\pi}=P_n, \\
		(h_\pi e_\pi)_{\pi \supseteq_n (\theta)} &\mapsto (h_\pi e_\pi \psi_{w(\pi)}1_{\bi^\pi})_{\pi \supseteq_n (\theta)}.
	\end{align*}
	To check that this yields an isomorphism of complexes, let $\tau := \refn^i(\pi)$ for some $i\in [a,b]\setminus C_\pi$, and check the following using Lemma~\ref{LPsiProduct}:
	\begin{itemize}
		\item $e_\tau \psi_{w(\tau)}\psi_{u(\tau, \pi)} = e_\tau\psi_{w(\tau, \pi)}\psi_{w(\pi)}$,
		\item $e_\pi \psi_{w(\pi)}\psi_{w(\pi)^{-1}} = e_\pi$,
		\item $1_{\bi^\pi}\psi_{w(\pi)^{-1}}\psi_{w(\pi)} = 1_{\bi^\pi}$. \qedhere
	\end{itemize}
\end{proof}

\subsection{The $\k$-module $\CE_\theta(\rho, \sigma)$}\label{SSRootExtMod}

Throughout this subsection, we fix $\rho, \sigma \in \KP(\theta)$.
Recall the word $\bj_\sigma\in I^\theta$ and the idempotent $e_\sigma$ from (\ref{ERootIdem}).
Writing $\bj_\sigma = j_1 j_2 \cdots j_l$, for $i\in [a, b+1]$, there exists a unique $r\in [1, l]$ such that $j_r = i$, and we denote $\idx_{\sigma}(i) := r$.

We have the standard module $\Delta(\sigma)$ with generator $v_\sigma\in e_{\sigma}\Delta(\sigma)$, see (\ref{EDeltaGen}).
As in \S\ref{SSStdMod}, we consider $\Delta(\sigma)$ as an $(R_\theta,\Lambda_\sigma)$-bimodule.
Recalling Lemma~\ref{LCochain}, we write $T_\bullet^\rho(\sigma) := T_\bullet^{S_\bullet^\rho}(\Delta(\sigma))$, so that $T_n^{\rho}(\sigma)=\bigoplus_{\pi\supseteq_n\rho}q^{-n}e_\pi\Delta(\sigma)$.
Recall (\ref{EYToLambda}) and (\ref{ERootShuffles}).

\begin{Lemma}\label{LDeltaBasis}
	In $\Delta(\sigma)$, we have:
	\begin{enumerate}
		\item\label{LDeltaBasisY} $y_r v_\sigma = v_\sigma\ttp_\sigma(y_r)$; in particular, $y_r v_\sigma = y_s v_\sigma$ if $r$ and $s$ are $\sigma$-equivalent;
		\item\label{LDeltaBasisPsi} $\psi_w v_\sigma = 0$ whenever $w \in \Sgp_l$ is not a $\sigma$-shuffle.
	\end{enumerate}
	Moreover, $\Delta(\sigma)$ is free as a right $\Lambda_\sigma$-module with basis $\{\psi_{w} v_\sigma \mid w \in \D^{\sigma}\}$.
\end{Lemma}
\begin{proof}
	Use Lemma~\ref{LDeltaRootBasis} and (\ref{EShuffleMod}).
\end{proof}

Recalling (\ref{ERootPerm}), we now get:

\begin{Lemma}\label{LRootWeight}
	Let $\pi,\sigma \in \KP(\theta)$.
	If $\sigma \nsupseteq \pi$, then $e_\pi \Delta(\sigma) = 0$.
	If $\sigma \supseteq_{n} \pi$, then there is an isomorphism of right $\Lambda_\sigma$-modules $$q^{n}\Lambda_\sigma \iso e_\pi\Delta(\sigma), \ f \mapsto \psi_{w(\pi,\sigma)} v_\sigma f.$$
\end{Lemma}
\begin{proof}
	If $\sigma$ is a refinement of $\pi$, then every pair of $\sigma$-equivalent integers is also $\pi$-equivalent, so $w(\pi,\sigma)$ is a $\sigma$-shuffle and the result follows from Lemma~\ref{LDeltaBasis} since in this case we have $\deg(\psi_{w(\pi,\sigma)}e_\sigma) = n$.

	On the other hand, if $\sigma \nsupseteq \pi$, then there is some $i\in C_\pi$ with $i\notin C_\sigma$.
	It follows that $\idx_\sigma(i)$ and $\idx_\sigma(i)+1$ are $\sigma$-equivalent but $w(\pi, \sigma)(\idx_\sigma(i)) > w(\pi, \sigma)(\idx_\sigma(i)+1)$, so $w(\pi, \sigma)$ is not a $\sigma$-shuffle, and $e_\pi\Delta(\sigma) = 0$ by Lemma~\ref{LDeltaBasis}(\ref{LDeltaBasisPsi}).
\end{proof}

\begin{Corollary}\label{CZeroComplex}
	We have:
	\begin{enumerate}
		\item\label{CZeroComplexi} If $\sigma \nsupseteq \rho$, then $T_\bullet^\rho(\sigma) = 0$.
		\item\label{CZeroComplexii} If $\sigma \supseteq_m \rho$, then as right $\Lambda_\sigma$-modules,
		$$T_n^\rho(\sigma)=\bigoplus_{\sigma\supseteq \pi \supseteq_n \rho}q^{-n} \psi_{w(\pi, \sigma)} v_\sigma\cdot\Lambda_\sigma \cong \bigoplus_{\sigma\supseteq \pi \supseteq_n \rho}q^{m-2n} \Lambda_\sigma.$$
		In particular, $T_n^\rho(\sigma) = 0$ for $n>m$ and
		$$T_m^\rho(\sigma) = q^{-m} e_\sigma \Delta(\sigma) = v_\sigma\cdot\Lambda_\sigma\cong q^{-m}\Lambda_\sigma.$$
	\end{enumerate}
\end{Corollary}
\begin{proof}
	If $\sigma$ is not a refinement of $\rho$, then it cannot be a refinement of any $\pi\supseteq \rho$, which implies (\ref{CZeroComplexi}).
	Let $\sigma \supseteq_m \rho$.
	Recall that $T^\rho_n(\sigma)=\bigoplus_{\pi\supseteq_n\rho}q^{-n}e_\pi\Delta(\sigma)$.
	Let $\pi\supseteq_n\rho$.
	If $\sigma \nsupseteq \pi$, then $e_\pi\Delta(\sigma)=0$.
	Otherwise $e_\pi\Delta(\sigma)=\psi_{w(\pi, \sigma)} v_\sigma\cdot\Lambda_\sigma$.
\end{proof}

\begin{Lemma}\label{LPsiSquared}
	Let $\pi \in \KP(\theta)$.
	Suppose that $\sigma \supseteq \refn^i(\pi)$ for some $i\in [a, b]$.
	Then
	$$\psi_{w(\refn^i(\pi), \pi)} \psi_{w(\pi, \sigma)} e_\sigma = \psi_{w(\refn^i(\pi), \sigma)}(y_{\idx_\sigma(i)}-y_{\idx_\sigma(i+1)})e_\sigma.$$
\end{Lemma}
\begin{proof}
	Let $\tau = \refn^i(\pi)$.
	We have that $i$ and $i+1$ are the only adjacent elements of $[a, b+1]$ with $\idx_\tau(i) > \idx_\tau(i+1)$ and $\idx_\pi(i) < \idx_\pi(i+1)$.
	Since $\sigma \supseteq \tau = \refn^i(\pi)$, we have $\idx_\sigma(i) > \idx_\sigma(i+1)$.
	Now we compute:
	\begin{align*}
		\psi_{w(\tau, \pi)} \psi_{w(\pi, \sigma)} e_\sigma
		&= \psi_{w(\tau, \pi)s_{\idx_\pi(i)}} \psi_{\idx_\pi(i)}^2 \psi_{s_{\idx_\pi(i)}w(\pi, \sigma)} e_\sigma \\
		&= \psi_{w(\tau, \pi)s_{\idx_\pi(i)}} (y_{\idx_\pi(i+1)} - y_{\idx_\pi(i)}) \psi_{s_{\idx_\pi(i)}w(\pi, \sigma)} e_\sigma \\
		&= \psi_{w(\tau, \pi)s_{\idx_\pi(i)}} \psi_{s_{\idx_\pi(i)}w(\pi, \sigma)} (y_{\idx_\sigma(i)} - y_{\idx_\sigma(i+1)}) e_\sigma \\
		&= \psi_{w(\tau, \sigma)}(y_{\idx_\sigma(i)} - y_{\idx_\sigma(i+1)}) e_\sigma,
	\end{align*}
	where the first equality is obtained using the fact that the braid relations (\ref{Rbraid}) hold without error term in $R_\theta$, the second comes by applying a non-trivial quadratic relation (\ref{Rquad}) on strands colored $i$ and $i+1$, the third is obtained using the relation (\ref{Rdots}), and the last comes from Lemma~\ref{LPsiProduct}.
\end{proof}

The proof of Theorem~\ref{TRootExt} amounts to showing that $T_\bullet^\rho(\sigma)$ is isomorphic to a certain Koszul complex (see~\cite[\S4.5]{We}) which we now define.
Suppose $\sigma\supseteq_m \rho$, and write $C_\sigma\setminus C_\rho = \{ i_1, \ldots, i_m\}$ with $i_1 < \cdots < i_m$.
Let $N$ be the free right $\Lambda_\sigma$-module of graded rank $mq^{-2}$, that is, $N := q^{-2}\Lambda_\sigma^{\oplus m}$.
For $k\in [1, m]$, we denote $\epsilon_k := (0,\ldots, 0, 1, 0, \ldots, 0)\in N$ (with ``$1$'' in the $k\mathrm{th}$ entry).
Recalling (\ref{EYToLambda}), define
\begin{align*}
	z_k &:= s(\refn^{i_k}(\rho),\rho)\ttp_\sigma(y_{\idx_\sigma(i_k)}-y_{\idx_\sigma(i_k+1)})\in \Lambda_\sigma \qquad (k=1,\dots,m), \\
	Z &:= (z_1, \ldots, z_m) = \epsilon_1 z_1 + \cdots + \epsilon_m z_m\in N.
\end{align*}
Note that $Z$ is a homogeneous degree $0$ element of $N$.
We consider the Koszul complex $q^m\Wed^\bullet N$ associated to the regular sequence $Z$ for the ring $\Lambda_\sigma$:
\begin{align*}
	0 \longleftarrow q^m\Wed^m N \longleftarrow \cdots \longleftarrow q^m\Wed^{n+1} N &\longleftarrow q^m\Wed^n N \longleftarrow \cdots \longleftarrow q^m\Wed^0 N \longleftarrow 0 \\
	Z\wedge a &\longmapsfrom a
\end{align*}
where $\Wed^n N$ is the $n\mathrm{th}$ exterior power of the free $\Lambda_\sigma$-module $N$.
Note that $\Wed^n N$ has a $\Lambda_\sigma$-basis $\{ \epsilon_{k_1}\wedge \cdots \wedge \epsilon_{k_n}\mid 1\leq k_1 < \cdots < k_n \leq m\}$.

Let $\sigma\supseteq_m \rho$.
By Corollary~\ref{CZeroComplex}(\ref{CZeroComplexii}), the $m\mathrm{th}$ component of $\Hom_{R_\theta}^\bullet(S_\bullet^\rho, \Delta(\sigma))$ is the last nonzero component and it can be identified with $\Hom_{R_\theta}(q^m R_\theta e_\sigma, \Delta(\sigma))$.
Thus, every element of the $m\mathrm{th}$ component is a cocycle, so there is a surjective map
$$[-]: \Hom_{R_\theta}(q^m R_\theta e_\sigma, \Delta(\sigma))\, \onto \,\CE_\theta^m(\rho, \sigma)=H^m(\Hom_{R_\theta}^\bullet(S_\bullet^\rho, \Delta(\sigma))), \ \varphi \mapsto [\varphi],$$
where $[\varphi]$ is the cohomology class of $\varphi$.
Moreover, by Lemma~\ref{LRootWeight}, we have an isomorphism
$$\xi:q^{-m}\Lambda_\sigma \iso \Hom_{R_\theta}(q^m R_\theta e_\sigma, \Delta(\sigma)), \ f \mapsto (e_\sigma\mapsto v_\sigma f).$$
We consider $\Lambda_\rho$ to be a $\Lambda_\sigma$-module via the homomorphism $\ttp_\rho^\sigma: \Lambda_\sigma \onto \Lambda_\rho$, see (\ref{ELambdaToLambda}).

\begin{Theorem}\label{TRootExt}
	Let $\rho, \sigma \in \KP(\theta)$.
	If $\sigma \nsupseteq \rho$, then $\CE_\theta(\rho, \sigma) = 0$.
	If $\sigma \supseteq_m \rho$, then $\CE_\theta(\rho, \sigma) = \CE_\theta^m(\rho, \sigma)$ and there is an isomorphism of $\Lambda_\sigma$-modules $\CE_\theta^m(\rho, \sigma) \iso q^{-m}\Lambda_\rho$ which makes the following diagram of $\Lambda_\sigma$-modules commute:
	$$\begin{tikzcd}
		q^{-m}\Lambda_\sigma \arrow[r,"\sim","\xi" below] \arrow[d,"\ttp_\rho^\sigma"] & \Hom_{R_\theta}(q^m R_\theta e_\sigma, \Delta(\sigma)) \arrow[d,"{[-]}"] \\
		q^{-m}\Lambda_\rho & \CE_\theta^m(\rho, \sigma) \arrow[l,dashed,"\sim"].
	\end{tikzcd}$$
\end{Theorem}
\begin{proof}
	If $\sigma \nsupseteq \rho$, then $\CE_\theta(\rho, \sigma) = 0$ by Lemma~\ref{LCochain} and Corollary~\ref{CZeroComplex}(\ref{CZeroComplexi}).

	Assume now that $\sigma \supseteq_m \rho$ and write $C_\sigma\setminus C_\rho = \{ i_1 < \cdots < i_m\}$.
	If $\sigma \supseteq \pi \supseteq_n \rho$, we set $B(\pi) := \{ k\in [1, m] \mid i_k\in C_\pi\}$.
	Note that $|B(\pi)| = n$.
	We define a map $\Theta_n: T_n^\rho(\sigma)\to q^m\Wed^n N$ of $\Lambda_\sigma$-modules by defining it on the $\Lambda_\sigma$-basis $\{ \psi_{w(\pi,\sigma)}v_\sigma \mid \sigma\supseteq \pi\supseteq_n \rho\}$ of $T_n^\rho(\sigma)$, see Corollary~\ref{CZeroComplex}.
	If $B(\pi) = \{ k_1 < \cdots < k_n\}$, then define
	$$\Theta_n(\psi_{w(\pi,\sigma)}v_\sigma) := \epsilon_{k_1}\wedge \cdots \wedge \epsilon_{k_n}.$$
	It is easy to see that $\Theta_n$ is an isomorphism of (graded) $\Lambda_\sigma$-modules.
	To show that $\Theta_n$ defines an isomorphism of complexes $T_\bullet^\rho(\sigma) \to q^m\Wed^\bullet N$, we must verify that the following square commutes:
	$$\begin{tikzcd}
		T_{n+1}^\rho(\sigma) \arrow[d, "\Theta_{n+1}"] & T_n^\rho(\sigma) \arrow[l, "d_n\cdot -" above] \arrow[d, "\Theta_n"] \\
		q^m\Wed^{n+1}N & q^m\Wed^n N, \arrow[l, "Z\wedge -" above]
	\end{tikzcd}$$
	where $d_n$ is the matrix defined in (\ref{ERootDiff}).
	We check this using an arbitrary basis element $\psi_{w(\pi,\sigma)}v_\sigma \in T_n^\rho(\sigma)$.
	We have
	\begin{align*}
		\Theta_{n+1}(d_n \psi_{w(\pi,\sigma)}v_\sigma) &= \sum_{\sigma \supseteq \tau\supseteq_1 \pi}s(\tau,\pi)\Theta_{n+1}(\psi_{w(\tau,\pi)}\psi_{w(\pi,\sigma)}v_\sigma) \\
		&= \sum_{k\in [1,m]\setminus B(\pi)}s(\refn^{i_k}(\pi),\pi)\Theta_{n+1}(\psi_{w(\refn^{i_k}(\pi),\pi)}\psi_{w(\pi,\sigma)}v_\sigma) \\
		&= \sum_{k\in [1,m]\setminus B(\pi)}(-1)^{B(\pi)\cap [1,k)}\Theta_{n+1}(\psi_{w(\refn^{i_k}(\pi),\sigma)}v_\sigma)z_k \\
		&= \sum_{k\in [1,m]\setminus B(\pi)}(-1)^{B(\pi)\cap [1,k)} (\epsilon_{k_1}\wedge \cdots \wedge \epsilon_k \wedge \cdots \wedge \epsilon_{k_n})z_k \\
		&= Z\wedge(\epsilon_{k_1}\wedge \cdots \wedge \epsilon_{k_n}) \\
		&= Z\wedge\Theta_n(\psi_{w(\pi,\sigma)}v_\sigma)
	\end{align*}
	where the second equality follows by noting that $k\mapsto \refn^{i_k}(\pi)$ defines a bijection from $[1,m]\setminus B(\pi)$ to the set of $1$-refinements of $\pi$ which are refined by $\sigma$, the third by Lemma~\ref{LPsiSquared} and the observation that $s(\refn^{i_k}(\pi),\pi) = (-1)^{B(\pi)\cap [1,k)}s(\refn^{i_k}(\rho),\rho)$, and the remaining equalities follow from the definitions.

	Since $q^m\Wed^\bullet N$ is a Koszul complex corresponding to a regular sequence, we now have that $\CE_\theta^n(\rho, \sigma) \cong H^n(q^m\Wed^\bullet N) = 0$ unless $n = m$.
	The proof is complete in view of Lemma~\ref{LCochain} upon noting that the kernel of $\ttp_\rho^\sigma: \Lambda_\sigma \onto \Lambda_\rho$ is the ideal generated by $(z_1, \ldots, z_m)$, so $\ttp_\rho^\sigma$ induces an isomorphism
	\[H^m(q^m\Wed^\bullet N) = \frac{q^{-m}\Lambda_\sigma}{(z_1, \ldots, z_m)} \iso q^{-m}\Lambda_\rho. \qedhere\]
\end{proof}

\subsection{The category $\CE_\theta$}\label{SSRootCat}

Throughout this subsection, we use Theorem~\ref{TRootExt} to identify $\CE_\theta(\rho, \sigma) = \CE_\theta^m(\rho, \sigma)$ with $q^{-m}\Lambda_\rho$ whenever $\sigma\supseteq_m \rho \in \KP(\theta)$.

Let $\sigma \supseteq_m\rho \in \KP(\theta)$.
For any $\hat{f}\in q^{-m}\CY_l$ and $\pi \supseteq_k \sigma$, set
$$\hat{f}_{\rho,\sigma}^\pi:= (-1)^{\frac{m(m+1)}{2}+mk} w(\pi, \sigma)\cdot \hat{f}.$$
We define an element of $\Hom_{R_\theta}^m(S_\bullet^\rho, S_\bullet^\sigma)$ by
\begin{align}\label{EExtLift}
	\varphi_{\rho,\sigma}^{\hat{f}}: S_{m+k}^\rho = \bigoplus_{\pi \supseteq_{m+k} \rho} q^{m+k}R_\theta e_\pi &\to \bigoplus_{\tau \supseteq_k \sigma}q^k R_\theta e_\tau = S_k^\sigma \nonumber \\
	(h_\pi e_\pi)_{\pi \supseteq_{m+k} \rho} &\mapsto (h_\tau \hat{f}_{\rho,\sigma}^\tau e_\tau)_{\tau \supseteq_k \sigma}.
\end{align}
Recalling the differential (\ref{EDGDiff}) on $\Hom_{R_\theta}^\bullet(S_\bullet^\rho, S_\bullet^\sigma)$, we have:

\begin{Lemma}\label{LRootExtLift}
	Let $\sigma \supseteq_m \rho \in \KP(\theta)$.
	If $f\in q^{-m}\Lambda_\rho = \CE_\theta^m(\rho, \sigma)$ and $\hat{f}\in q^{-m}\CY_l$ are such that $\ttp_\rho(\hat{f}) = f$, then
	\begin{enumerate}
		\item\label{LRootExtLifti} $\delta(\varphi_{\rho,\sigma}^{\hat{f}}) = 0$, and
		\item\label{LRootExtLiftii} the isomorphism (\ref{EExtIso}) sends the cohomology class of $\varphi_{\rho,\sigma}^{\hat{f}}$ to $f$.
	\end{enumerate}
\end{Lemma}
\begin{proof}
	We prove (\ref{LRootExtLifti}) by checking that the following diagram either commutes (if $m$ is even) or anticommutes (if $m$ is odd) whenever $\tau \supseteq_1 \pi\supseteq_k \sigma$:
	$$\begin{tikzcd}[column sep = 4em]
		q^{k+m+1}R_\theta e_\tau \arrow[r,"{-\cdot d_{k+m}^{\tau,\pi}}"] \arrow[d,"{-\cdot\hat{f}^\tau_{\rho,\sigma}}" left] & q^{k+m}R_\theta e_\pi \arrow[d,"{-\cdot\hat{f}^\pi_{\rho,\sigma}}"] \\
		q^{k+1}R_\theta e_\tau \arrow[r, "{-\cdot d_k^{\tau,\pi}}"] & q^k R_\theta e_\pi.
	\end{tikzcd}$$
	Modulo the signs, this is checked by the computation:
	\begin{align*}
		\psi_{w(\tau,\pi)}(w(\pi,\sigma)\cdot \hat{f}) e_\pi =
		(w(\tau, \pi)w(\pi,\sigma)\cdot\hat{f})\psi_{w(\tau,\pi)}  e_\pi
		= (w(\tau,\sigma)\cdot \hat{f})\psi_{w(\tau,\pi)} e_\pi,
	\end{align*}
	and the signs are taken care of by
	$$(-1)^{\frac{m(m+1)}{2}+mk}s(\tau,\pi) = (-1)^m(-1)^{\frac{m(m+1)}{2}+m(k+1)}s(\tau, \pi).$$

	To prove (\ref{LRootExtLiftii}), first note that the restriction of $\varphi_{\rho,\sigma}^{\hat{f}}$ to $S_m^\rho$ has image in $S_0^\sigma$ and can therefore be realized as
	\begin{align*}
		\varphi_{\rho,\sigma}^{\hat{f}}|_{S_m^\rho}: S_m^\rho = \bigoplus_{\pi \supseteq_m \rho}q^m R_\theta e_\pi &\to R_\theta e_\sigma = S_0^\sigma \\
		(h_\pi e_\pi)_{\pi \supseteq_m \rho} &\mapsto h_\sigma \hat{f}_{\rho,\sigma}^\sigma e_\sigma = (-1)^{\frac{m(m+1)}{2}} h_\sigma \hat f e_\sigma.
	\end{align*}
	By Corollary~\ref{CZeroComplex}(\ref{CZeroComplexii}), we identify $\Hom_{R_\theta}(S_m^\rho, \Delta(\sigma))$ with $\Hom_{R_\theta}(q^m R_\theta e_\sigma, \Delta(\sigma))$.
	Then the image of $\varphi_{\rho,\sigma}^{\hat{f}}$ under the map (\ref{EExtInduced}) is
	$$(-1)^{\frac{m(m+1)}{2}}\eps_\sigma(\varphi_{\rho,\sigma}^{\hat{f}}|_{S_m^\rho}) = (e_\sigma \mapsto v_\sigma f) \in \Hom_{R_\theta}(q^m R_\theta e_\sigma, \Delta(\sigma)).$$
	An application of Theorem~\ref{TRootExt} completes the proof.
\end{proof}

\begin{Lemma}\label{LAltLift}
	Let $\tau \supseteq \sigma\supseteq_m \rho \in \KP(\theta)$.
	If $f\in q^{-m}\Lambda_\rho$ and $\hat{f}\in q^{-m}\CY_l$ are such that $\ttp_\rho(\hat{f}) = f$, then we also have $\ttp_\rho(w(\tau,\sigma)\cdot\hat{f}) = f$.
\end{Lemma}
\begin{proof}
	Since $\tau\supseteq \sigma \supseteq \rho$, we have that $r$ and $w(\tau, \sigma)(r)$ are $\rho$-equivalent for any $r\in [1, l]$, so $\ttp_\rho(w(\tau,\sigma)\cdot y_r) = \ttp_\rho(y_r)$.
	This implies the result.
\end{proof}

\begin{Lemma}\label{LRootLiftProd}
	Let $\tau \supseteq_n \sigma\supseteq_m \rho \in \KP(\theta)$.
	If $f\in q^{-m}\Lambda_\rho$, and $g\in q^{-n}\Lambda_\sigma$, then there exist $\hat{f}\in q^{-m}\CY_l$, $\hat{g}\in q^{-n}\CY_l$, and $\widehat{\ttp_\rho^\sigma(g)f}\in q^{-(m+n)}\CY_l$ with $\ttp_\rho(\hat{f}) = f$, $\ttp_\sigma(\hat{g}) = g$, and $\ttp_\rho(\widehat{\ttp_\rho^\sigma(g)f}) = \ttp_\rho^\sigma(g)f$, such that
	$$\varphi_{\sigma,\tau}^{\hat{g}}\varphi_{\rho,\sigma}^{\hat{f}} = \varphi_{\rho,\tau}^{\widehat{\ttp_\rho^\sigma(g)f}}.$$
\end{Lemma}
\begin{proof}
	Choose any two lifts $\hat{f}\in q^{-m}\CY_l$ and $\hat{g}\in q^{-n}\CY_l$ of $f$ and $g$, respectively.
	By Lemma~\ref{LAltLift}, since $\tau\supseteq \sigma$, we have that $w(\tau, \sigma)\cdot \hat{f}$ is also a lift of $f$, and so $\widehat{\ttp_\rho^\sigma(g)f} := \hat{g}(w(\tau, \sigma)\cdot\hat{f})\in q^{-(m+n)}\CY_l$ is a lift of $\ttp_\rho^\sigma(g)f$.
	By (\ref{EExtLift}), it suffices to show that for any $\pi\supseteq \tau$ we have $(\widehat{\ttp_\rho^\sigma(g)f})_{\rho,\tau}^\pi = \hat{g}_{\sigma,\tau}^\pi \hat{f}_{\rho,\sigma}^\pi$.
	Modulo the signs, this holds by the following computation:
	\begin{align*}
		w(\pi, \tau) \cdot (\hat{g}(w(\tau,\sigma)\cdot \hat{f})) &=
		(w(\pi,\tau)\cdot\hat{g})(w(\pi,\tau)w(\tau, \sigma)\cdot \hat{f}) \\
		&= (w(\pi, \tau)\cdot\hat{g})(w(\pi, \sigma)\cdot \hat{f}),
	\end{align*}
	and if $\pi\supseteq_k \tau$, the signs are taken care of by
	\[(-1)^{\frac{(m+n)(m+n+1)}{2}+(m+n)k} = (-1)^{\frac{m(m+1)}{2}+m(k+n)}(-1)^{\frac{n(n+1)}{2}+nk}.\qedhere\]
\end{proof}

We combine Lemmas~\ref{LRootLiftProd} and~\ref{LRootExtLift} to obtain the following theorem.

\begin{Theorem}\label{TRootAlg}
	Let $\tau\supseteq_n \sigma\supseteq_m \rho \in \KP(\theta)$.
	The composition in the category $\CE_\theta$ is given by
	$$\begin{tikzcd}[row sep = -0.5em]
		\CE_\theta^n(\sigma, \tau) \otimes \CE_\theta^m(\rho, \sigma) \arrow[r] & \CE_\theta^{m+n}(\rho, \tau) \\
		\rotatebox{90}{$\cong$} & \rotatebox{90}{$\cong$} \\
		q^{-n}\Lambda_\sigma \otimes q^{-m}\Lambda_\rho \arrow[r] & q^{-(m+n)}\Lambda_\rho \\
		g \otimes f \arrow[r, mapsto] & \ttp_\rho^\sigma(g)f.
	\end{tikzcd}$$
\end{Theorem}

\section{The $\ttA[2]$ case}\label{SA2}

Throughout this section, we use a special notation
$$\alpha:=\alpha_1,\ \beta:=\alpha_2,\ \gamma:=\alpha_1+\alpha_2,$$
so that $\alpha,\beta,\gamma$ are now the positive roots of the root system of type $\ttA[2]$.
We fix
$$\theta := a\alpha + b\beta$$ with $a, b\in\Z_{\geq 0}$.
There is a bijection
$$\sigma: [0, \min\{ a, b\} ]\iso \KP(\theta),\ s \mapsto \sigma(s) := (\beta^{b-s}, \gamma^s, \alpha^{a-s}).$$
The standard $R_\theta$-modules are
$$\Delta(s) := \Delta(\sigma(s))\qquad(0\leq s\leq \min\{ a, b\}).$$
We denote the standard generator of $\Delta(s)$ by $v_s:=v_{\sigma(s)}$, see (\ref{EDeltaGen}).
Recall that
$$\End_{R_\theta}(\Delta(s))^\op \cong \Lambda_{\sigma(s)} = \Lambda_{b-s,s,a-s}=\Lambda_{b-s}\otimes\Lambda_s\otimes\Lambda_{a-s},$$
see (\ref{EDeltaEnd}).

\subsection{The resolution $P_\bullet^r$}\label{SSA2Res}

Let $0\leq r\leq \min\{ a, b\}$.
We recall the resolution of $\Delta(r)$ defined in~\cite{BKS}.
For $n\in [0, r]$, recalling (\ref{EThickTransposition}), we define:
\begin{align*}
	\bj_{r,n} &:= 2^{b-r}1^{r-n}2^{r}1^{n+a-r} \in I^\theta, \\
	\bi_{r,n} &:= 2^{(b-r)}1^{(r-n)}2^{(r)}1^{(n)}1^{(a-r)} \in I^\theta_\di, \\
	e_{r,n} &:= 1_{\bi_{r,n}} \in R_\theta, \\
	\tts_{r,n} &:= n(r-n+1) - \binom{a-r}{2} - \binom{b-r}{2} - r(r-1), \\
	x_{r,n} &:= U_{b-n-1;1,r+n}=(b+r,b+r-1,\dots,b-n) \in \Sgp_{a+b}, \\
	d_{r,n} &:= e_{r,n+1}\psi_{x_{r,n}} e_{r,n} \in R_\theta, \\
	P_n^r &:= q^{\tts_{r,n}}R_\theta e_{r,n}.
\end{align*}
Note that the right multiplication with $d_{r,n}$ yields the degree zero $R_\theta$-homomorphism $-\cdot d_{r,n}:P^r_{n+1}\to P^r_n$.
Define $u_r\in \Sgp_{2r}$ by
\begin{equation}\label{ELongGatherer}
	u_r(i) = \begin{cases} r-\frac{i-1}{2} & \text{if }i\text{ is odd}, \\ 2r-\frac{i-2}{2} & \text{if }i\text{ is even}. \end{cases}
\end{equation}
Let
\begin{equation}\label{EA2Aug}
	\eps_r: P_0^r \onto \Delta(r), \ x e_{r,0} \mapsto x(\psi_{w_{0,b-r}}\circ \psi_{u_r} \circ \psi_{w_{0,a-r}}) v_r.
\end{equation}
By~\cite[Theorem A]{BKS} and Lemma~\ref{LCircRes}, we have:

\begin{Lemma}\label{LA2Res}
	The following sequence is a projective resolution of $\Delta(r)$:
	$$0 \longrightarrow P_r^r \longrightarrow \cdots \longrightarrow P_{n+1}^r \stackrel{-\cdot d_{r,n}}{\longrightarrow} P_n^r \longrightarrow \cdots \longrightarrow P_0^r \stackrel{\eps_r}{\longrightarrow} \Delta(r) \longrightarrow 0.$$
\end{Lemma}

\subsection{Weight spaces of standard modules}

The following lemmas are useful for finding bases for certain weight spaces of the standard modules $\Delta(s)$.
The first of them concerns the nil-Hecke algebra and is well-known and easy to check.
Recall the notation from \S\ref{SSSym} and the $(R_{c\alpha_i}, \NH_c)$-bimodule structure on $\hat{\Delta}(\alpha_i^c)$ from \S\ref{SSStdMod}.

\begin{Lemma}\label{LSimpleBasis}
	Let $i\in I$ and $c\in\Z_{\geq 0}$.
	The map
	$$q^{-\binom{c}{2}}\CX_c \to \hat{\Delta}(\alpha_i^c), \ f\mapsto v_{\alpha_i}^{\circ c} f \tau_{w_0}$$
	is an injective map of right $\Lambda_c$-modules with image $\Delta(\alpha_i^c)$.
\end{Lemma}

Recalling (\ref{ELongGatherer}), we have:

\begin{Lemma}\label{L12Basis}
	The map
	$$q^{-2\binom{s}{2}}\Lambda_s \to \hat{\Delta}(\gamma^s), \ f\mapsto \psi_{u_s}v_{\gamma}^{\circ s} f$$
	is an injective map of right $\Lambda_s$-modules with image $1_{1^{(s)}2^{(s)}}\Delta(\gamma^s)$.
\end{Lemma}
Diagrammatically, the map in the lemma is given by
$$f\mapsto
\begin{braid}\tikzset{baseline=2em}
	\draw (0.5,0) -- (0.5,0.5);
	\draw (4.5,0) -- (4.5,0.5);
	\draw (0.5,1.5) -- (0.5,2) -- (0,3) -- (0,4);
	\draw (4.5,1.5) -- (4.5,2) -- (4,3) -- (2,4);
	\draw (0.5,1.5) -- (0.5,2) -- (1,3) -- (3,4);
	\draw (4.5,1.5) -- (4.5,2) -- (5,3) -- (5,4);
	\draw (0,5) -- (0,5.5);
	\draw (2,5) -- (2,5.5);
	\draw (3,5) -- (3,5.5);
	\draw (5,5) -- (5,5.5);
	\braidbox{0.5}{4.5}{0.5}{1.5}{\scriptsize $f$}
	\braidbox{0}{2}{4}{5}{$w_0$}
	\braidbox{3}{5}{4}{5}{$w_0$}
	\draw (2.5,2) node{\tiny $\cdots$};
	\draw (1,6) node{\tiny $1^s$};
	\draw (4,6) node{\tiny $2^s$};
\end{braid}.$$
\begin{proof}
	In view of Lemma~\ref{LGDim}, we may assume that $\k$ is a field.
	Let $M$ be the free graded $\k$-module with basis $\D^{(2^s)}$, where the degree of the basis element $w\in \D^{(2^s)}$ is set to equal to the degree of $\psi_w 1_{(12)^s}$ in $R_{s\alpha+s\beta}$.
	Lemma~\ref{LDeltaRootBasis} and (\ref{EShuffleMod}) show that the map
	$\xi: q^{\binom{s}{2}}M\otimes_\k \CX_s \to \hat{\Delta}(\gamma^s), \ w\otimes f\mapsto \psi_w v_{\gamma}^{\circ s} f$
	is an isomorphism of right $\CX_s$-modules.

	Let $\varphi$ be the map in the statement.
	Then $\varphi$ is the composition of the inclusion $q^{-2\binom{s}{2}}\Lambda_s \into q^{\binom{s}{2}}M\otimes_\k \CX_s, \ f \mapsto u_s\otimes f$ with $\xi$, so $\varphi$ is injective.
	It is now enough to show that $\varphi(f) \in 1_{1^{(s)}2^{(s)}}\Delta(\gamma^s)$ and that $\qdim(1_{1^{(s)}2^{(s)}}\Delta(\gamma^s)) = q^{-2\binom{s}{2}}\qdim(\Lambda_s)$.
	For the first claim, note that $\psi_{u_s} v_{\gamma}^{\circ s}$ is a nonzero element of smallest possible degree in $1_{1^s2^s}\hat{\Delta}(\gamma^s)$, so recalling (\ref{ENHIdem}), we have
	\begin{align*}
		1_{1^{(s)}2^{(s)}}\psi_{u_s}v_{\gamma}^{\circ s} f e_s &=
		\psi_{u_s} v_{\gamma}^{\circ s} f x_0\tau_{w_0} \\
		&= \psi_{u_s} v_{\gamma}^{\circ s} f x_0\tau_1\tau_{s_1^{-1}w_0} \\
		&= \psi_{u_s} v_{\gamma}^{\circ s} \tau_1 f x_0 \tau_{s_1^{-1}w_0} + \psi_{u_s} v_{\gamma}^{\circ s} \partial_1(f x_0) \tau_{s_1^{-1}w_0} \\
		&= \psi_{u_s} v_{\gamma}^{\circ s} \partial_1(f x_0) \tau_{s_1^{-1}w_0} \\
		&\phantom{=}\vdots \\
		&= \psi_{u_s} v_{\gamma}^{\circ s} \partial_{w_0}(fx_0) \\
		&= \psi_{u_s} v_{\gamma}^{\circ s} f,
	\end{align*}
	where the first equality comes from Lemma~\ref{LDivPower} and (\ref{ENHIdem}), the third comes from the relations in $\NH_s$, the fourth follows because $\deg(\psi_{u_s} v_{\gamma}^{\circ s}\tau_1) =\deg(\psi_{u_s} v_{\gamma}^{\circ s})-2< \deg(\psi_{u_s} v_{\gamma}^{\circ s})$, and the last holds by Lemma~\ref{LLambda}.

	As for graded dimension, we have
	\begin{align*}
		\qdim(1_{1^{(s)}2^{(s)}}\Delta(\gamma^s))
		&= \frac{1}{\left([s]_+^!\right)^3} \qdim(1_{1^{s}2^{s}}\hat{\Delta}(\gamma^s)) \\
		&= \frac{q^{2\binom{s}{2}}\left([s]_-^!\right)^2}{\left([s]_+^!\right)^3} (\qdim \CX_1)^s \\
		&= q^{-2\binom{s}{2}}\qdim\Lambda_s.
	\end{align*}
	where the first equality follows from (\ref{EDeltaIso}) and (\ref{EUndivided}), the second from (\ref{EDeltaHat}), (\ref{EShuffle}), and Lemma~\ref{LDeltaRootBasis}, and the last from an elementary computation.
\end{proof}

\subsection{The $\k$-module $\CE_\theta(r, s)$}\label{SSA2ExtMod}

In this subsection we fix $r,s \in [0, \min\{ a, b\} ]$.
Recalling Lemma~\ref{LCochain}, we write $T_\bullet^r(s) := T_\bullet^{P_\bullet^r}(\Delta(s))$.
The terms of this complex are of the form
$$T_n^r(s) = q^{-\tts_{r,n}}e_{r,n}\Delta(s)\qquad(n =0,\ldots,r).$$
If $r \geq s$ and $0\leq n\leq r-s$, we define
\begin{align*}
	\omega_n(r, s) &:= -(r-s)(1+(a-r)+(b-r))+(r-s-n+1)(r-s-n), \\
	K_n &:= q^{\omega_n(r,s)}\Lambda_{b-r, r-s, s, r-s-n, n, a-r}, \\
	w_n &:= U_{b-r;s,r-s-n}U_{b-n;r-s,s}U_{b-r+s;r-s,r-s-n}U_{b-r;r-s,s}U_{b;s,r-s-n}.
\end{align*}
The diagram for $1_{\bj_{r,n}}\psi_{w_n} 1_{\bj_{s,0}}$ is the top part of the diagram below.
Observe that $\Lambda_{\sigma(s)}=\Lambda_{b-s,s,a-s}\subseteq K_n$ in a natural way, so we may consider $K_n$ as a $\Lambda_{\sigma(s)}$-module.
Recalling (\ref{ELongGatherer}), we have

\begin{Lemma}\label{LWeightBasis}
	Suppose $0 \leq n \leq r$.
	\begin{enumerate}
		\item\label{LWeightBasisi} If $n > r-s$, then $T_n^r(s) = 0$.
		\item\label{LWeightBasisii} If $n\leq r-s$, then the map
		\begin{align*}
			\Xi_n: K_n &\to q^{-\tts_{r,n}}\hat{\Delta}(s), \\
			f &\mapsto \psi_{w_n}(v_{\beta}^{\circ b-s}\circ \psi_{u_s}v_{\gamma}^{\circ s} \circ v_{\alpha}^{\circ a-s})\, f\, (\tau_{w_{0,b-s}}\otimes 1 \otimes \tau_{w_{0,a-s}})
		\end{align*}
		is an injective degree zero map of $\Lambda_{\sigma(s)}$-modules with image $T_n^r(s)$.
	\end{enumerate}
\end{Lemma}

Diagrammatically, $\Xi_n$ is given by
$$f\mapsto
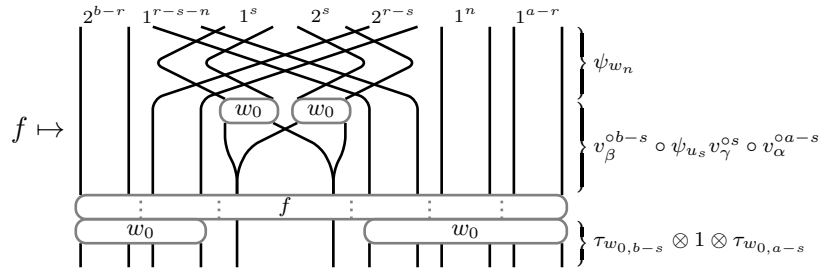
\begin{braid}\tikzset{baseline=4.5em}
	\draw (0,0) -- (0,1);
	\draw (2,0) -- (2,1);
	\draw (3,0) -- (3,1);
	\draw (5,0) -- (5,1);
	\draw (6.5,0) -- (6.5,2);
	\draw (10.5,0) -- (10.5,2);
	\draw (12,0) -- (12,1);
	\draw (14,0) -- (14,1);
	\draw (15,0) -- (15,1);
	\draw (17,0) -- (17,1);
	\draw (18,0) -- (18,1);
	\draw (20,0) -- (20,1);
	\draw (0,3) -- (0,10);
	\draw (2,3) -- (2,10);
	\draw (3,3) -- (3,7) -- (12,10);
	\draw (5,3) -- (5,7) -- (14,10);
	\draw (6.5,3) -- (6.5,4) -- (6,5) -- (6,6);
	\draw (6.5,3) -- (6.5,4) -- (7,5) -- (9,6);
	\draw (10.5,3) -- (10.5,4) -- (10,5) -- (8,6);
	\draw (10.5,3) -- (10.5,4) -- (11,5) -- (11,6);
	\draw (12,3) -- (12,7) -- (3,10);
	\draw (14,3) -- (14,7) -- (5,10);
	\draw (15,3) -- (15,10);
	\draw (17,3) -- (17,10);
	\draw (18,3) -- (18,10);
	\draw (20,3) -- (20,10);
	\draw (6,7) -- (3,8.5) -- (6,10);
	\draw (8,7) -- (5,8.5) -- (8,10);
	\draw (9,7) -- (12,8.5) -- (9,10);
	\draw (11,7) -- (14,8.5) -- (11,10);
	\braidbox{0}{5}{1}{2}{$w_0$};
	\braidbox{12}{20}{1}{2}{$w_0$};
	\braidbox{0}{20}{2}{3}{};
	\draw (8.5,2.5) node{\scriptsize $f$};
	\braidbox{6}{8}{6}{7}{$w_0$};
	\braidbox{9}{11}{6}{7}{$w_0$};
	\draw[dotted, color=gray] (2.5,2) -- (2.5,3);
	\draw[dotted, color=gray] (5.75,2) -- (5.75,3);
	\draw[dotted, color=gray] (11.25,2) -- (11.25,3);
	\draw[dotted, color=gray] (14.5,2) -- (14.5,3);
	\draw[dotted, color=gray] (17.5,2) -- (17.5,3);
	\draw (1,10.5) node{\tiny $2^{b-r}$};
	\draw (4,10.5) node{\tiny $1^{r-s-n}$};
	\draw (7,10.5) node{\tiny $1^s$};
	\draw (10,10.5) node{\tiny $2^s$};
	\draw (13,10.5) node{\tiny $2^{r-s}$};
	\draw (16,10.5) node{\tiny $1^n$};
	\draw (19,10.5) node{\tiny $1^{a-r}$};
	\draw (20,8.5) node[right]{\rotatebox{-90}{$\overbrace{\hspace{4em}}^{\rotatebox{90}{\scriptsize $\psi_{w_n}$}}$}};
	\draw (20,5) node[right]{\rotatebox{-90}{$\overbrace{\hspace{5em}}^{\rotatebox{90}{\scriptsize $v_{\beta}^{\circ b-s}\circ \psi_{u_s}v_{\gamma}^{\circ s} \circ v_{\alpha}^{\circ a-s}$}}$}};
	\draw (20,1) node[right]{\rotatebox{-90}{$\overbrace{\hspace{2.5em}}^{\rotatebox{90}{\scriptsize $\tau_{w_{0,b-s}}\otimes 1 \otimes \tau_{w_{0,a-s}}$}}$}};
\end{braid}$$

\begin{proof}
	(\ref{LWeightBasisi}) The condition $n > r-s$ is equivalent to the condition $a-r+n > a-s$, which easily implies, by (\ref{EShuffle}), that $1_{\bj_{r,n}}\Delta(s) = 0$ hence $e_{r,n}\Delta(s)=0$.

	(\ref{LWeightBasisii}) It is straightforward to check that $\Xi_n$ is homogeneous and $\Lambda_{\sigma(s)}$-equivariant.
	Define
	\begin{align*}
		\xi_\alpha &: q^{-\binom{a-s}{2}}\Lambda_{r-s-n,n,a-r}\longinto q^{-\binom{a-s}{2}} \CX_{a-s} \iso \Delta(\alpha^{a-s}) \\
		\xi_\beta &: q^{-\binom{b-s}{2}}\Lambda_{b-r,r-s}\longinto q^{-\binom{b-s}{2}} \CX_{b-s} \iso \Delta(\beta^{b-s}) \\
		\xi_\gamma &: q^{-2\binom{s}{2}}\Lambda_s \iso 1_{1^{(s)}2^{(s)}}\Delta(\gamma^s)
	\end{align*}
	where the isomorphisms are from Lemmas~\ref{LSimpleBasis} and~\ref{L12Basis}.
	Let $d = (r-s)(r-s-n)-s(r-s-n)-s(r-s)$ and define
	$$\xi : \Delta(\beta^{b-s}) \otimes 1_{1^{(s)}2^{(s)}}\Delta(\gamma^s) \otimes \Delta(\alpha^{a-s}) \rightarrow q^{-d} 1_{\bj_{r,n}}\Delta(s), \ x\otimes y \otimes z \mapsto \psi_{w_n}(x\circ y \circ z).$$
	Since $w_n \in \D^{(b-s, 2s, a-s)}$, (\ref{EShuffleMod}) shows that $\xi$ is injective.
	Observing that $\Xi_n = \xi(\xi_\beta\otimes \xi_\gamma\otimes \xi_\alpha)$, we see that $\Xi_n$ is injective with image in $1_{\bj_{r,n}}\Delta(s)$.

	Now we prove that $\im\Xi_n\subseteq e_{r,n}\Delta(s)$ by showing that $e_{r,n}\Xi_n(f) = \Xi_n(f)$.
	We have, using the relations in $R_\theta$,
	\begin{align*}
		\Xi_n(f) &=
		\begin{braid}\tikzset{baseline=4.5em}
			\draw (6.5,0) -- (6.5,3);
			\draw (10.5,0) -- (10.5,3);
			\draw (0,4) -- (0,10);
			\draw (2,4) -- (2,10);
			\draw (3,4) -- (3,6) -- (9,9) -- (12,10);
			\draw (5,4) -- (5,6) -- (11,9) -- (14,10);
			\draw (6.5,4) -- (6.5,5) -- (6,6) -- (3,8);
			\draw (6.5,4) -- (6.5,5) -- (7,6) -- (12,8);
			\draw (10.5,4) -- (10.5,5) -- (10,6) -- (5,8);
			\draw (10.5,4) -- (10.5,5) -- (11,6) -- (14,8);
			\draw (12,4) -- (12,6) -- (6,9) -- (3,10);
			\draw (14,4) -- (14,6) -- (8,9) -- (5,10);
			\draw (15,4) -- (15,10);
			\draw (17,4) -- (17,10);
			\draw (18,4) -- (18,10);
			\draw (20,4) -- (20,10);
			\draw (3,9) -- (6,10);
			\draw (5,9) -- (8,10);
			\draw (12,9) -- (9,10);
			\draw (14,9) -- (11,10);
			\braidbox{0}{5}{0}{1}{$w_0$};
			\braidbox{12}{20}{0}{1}{$w_0$};
			\braidbox{0}{2}{1}{2}{$x_0$};
			\braidbox{3}{5}{1}{2}{$x_0$};
			\braidbox{12}{14}{1}{2}{$x_0$};
			\braidbox{15}{17}{1}{2}{$x_0$};
			\braidbox{18}{20}{1}{2}{$x_0$};
			\braidbox{0}{2}{2}{3}{$w_0$};
			\braidbox{3}{5}{2}{3}{$w_0$};
			\braidbox{12}{14}{2}{3}{$w_0$};
			\braidbox{15}{17}{2}{3}{$w_0$};
			\braidbox{18}{20}{2}{3}{$w_0$};
			\braidbox{0}{20}{3}{4}{};
			\draw (8.5,3.5) node{\scriptsize $f$};
			\braidbox{3}{5}{8}{9}{$w_0$};
			\braidbox{12}{14}{8}{9}{$w_0$};
			\draw[dotted, color=gray] (2.5,3) -- (2.5,4);
			\draw[dotted, color=gray] (5.75,3) -- (5.75,4);
			\draw[dotted, color=gray] (11.25,3) -- (11.25,4);
			\draw[dotted, color=gray] (14.5,3) -- (14.5,4);
			\draw[dotted, color=gray] (17.5,3) -- (17.5,4);
			\draw (1,10.5) node{\tiny $2^{b-r}$};
			\draw (4,10.5) node{\tiny $1^{r-s-n}$};
			\draw (7,10.5) node{\tiny $1^s$};
			\draw (10,10.5) node{\tiny $2^s$};
			\draw (13,10.5) node{\tiny $2^{r-s}$};
			\draw (16,10.5) node{\tiny $1^n$};
			\draw (19,10.5) node{\tiny $1^{a-r}$};
		\end{braid} \\
		&=
		\begin{braid}\tikzset{baseline=4em}
			\draw (6.5,0) -- (6.5,2);
			\draw (10.5,0) -- (10.5,2);
			\draw (0,3) -- (0,8);
			\draw (2,3) -- (2,8);
			\draw (3,3) -- (3,5) -- (9,7);
			\draw (5,3) -- (5,5) -- (11,7);
			\draw (6.5,3) -- (6.5,4) -- (6,5) -- (3,7);
			\draw (6.5,3) -- (6.5,4) -- (7,5) -- (12,7);
			\draw (10.5,3) -- (10.5,4) -- (10,5) -- (5,7);
			\draw (10.5,3) -- (10.5,4) -- (11,5) -- (14,7);
			\draw (12,3) -- (12,5) -- (6,7);
			\draw (14,3) -- (14,5) -- (8,7);
			\draw (15,3) -- (15,8);
			\draw (17,3) -- (17,8);
			\draw (18,3) -- (18,8);
			\draw (20,3) -- (20,8);
			\draw (3,8) -- (6,9);
			\draw (5,8) -- (8,9);
			\draw (6,8) -- (3,9);
			\draw (8,8) -- (5,9);
			\draw (9,8) -- (12,9);
			\draw (11,8) -- (14,9);
			\draw (12,8) -- (9,9);
			\draw (14,8) -- (11,9);
			\braidbox{0}{5}{0}{1}{$w_0$};
			\braidbox{12}{20}{0}{1}{$w_0$};
			\braidbox{0}{2}{1}{2}{$x_0$};
			\braidbox{3}{5}{1}{2}{$x_0$};
			\braidbox{12}{14}{1}{2}{$x_0$};
			\braidbox{15}{17}{1}{2}{$x_0$};
			\braidbox{18}{20}{1}{2}{$x_0$};
			\braidbox{0}{20}{2}{3}{};
			\draw (8.5,2.5) node{\scriptsize $f$};
			\braidbox{3}{5}{7}{8}{$w_0$};
			\braidbox{6}{8}{7}{8}{$w_0$};
			\braidbox{9}{11}{7}{8}{$w_0$};
			\braidbox{12}{14}{7}{8}{$w_0$};
			\braidbox{0}{2}{8}{9}{$w_0$};
			\braidbox{15}{17}{8}{9}{$w_0$};
			\braidbox{18}{20}{8}{9}{$w_0$};
			\draw[dotted, color=gray] (2.5,2) -- (2.5,3);
			\draw[dotted, color=gray] (5.75,2) -- (5.75,3);
			\draw[dotted, color=gray] (11.25,2) -- (11.25,3);
			\draw[dotted, color=gray] (14.5,2) -- (14.5,3);
			\draw[dotted, color=gray] (17.5,2) -- (17.5,3);
			\draw (1,9.5) node{\tiny $2^{b-r}$};
			\draw (4,9.5) node{\tiny $1^{r-s-n}$};
			\draw (7,9.5) node{\tiny $1^s$};
			\draw (10,9.5) node{\tiny $2^s$};
			\draw (13,9.5) node{\tiny $2^{r-s}$};
			\draw (16,9.5) node{\tiny $1^n$};
			\draw (19,9.5) node{\tiny $1^{a-r}$};
		\end{braid} \\
		&=
		\begin{braid}\tikzset{baseline=3.5em}
			\draw (6.5,0) -- (6.5,2);
			\draw (10.5,0) -- (10.5,2);
			\draw (0,3) -- (0,7);
			\draw (2,3) -- (2,7);
			\draw (3,3) -- (3,5) -- (9,7);
			\draw (5,3) -- (5,5) -- (11,7);
			\draw (6.5,3) -- (6.5,4) -- (6,5) -- (3,7);
			\draw (6.5,3) -- (6.5,4) -- (7,5) -- (12,7);
			\draw (10.5,3) -- (10.5,4) -- (10,5) -- (5,7);
			\draw (10.5,3) -- (10.5,4) -- (11,5) -- (14,7);
			\draw (12,3) -- (12,5) -- (6,7);
			\draw (14,3) -- (14,5) -- (8,7);
			\draw (15,3) -- (15,7);
			\draw (17,3) -- (17,7);
			\draw (18,3) -- (18,7);
			\draw (20,3) -- (20,7);
			\braidbox{0}{5}{0}{1}{$w_0$};
			\braidbox{12}{20}{0}{1}{$w_0$};
			\braidbox{0}{2}{1}{2}{$x_0$};
			\braidbox{3}{5}{1}{2}{$x_0$};
			\braidbox{12}{14}{1}{2}{$x_0$};
			\braidbox{15}{17}{1}{2}{$x_0$};
			\braidbox{18}{20}{1}{2}{$x_0$};
			\braidbox{0}{20}{2}{3}{};
			\draw (8.5,2.5) node{\scriptsize $f$};
			\braidbox{3}{8}{7}{8}{$w_0$};
			\braidbox{9}{14}{7}{8}{$w_0$};
			\braidbox{0}{2}{7}{8}{$w_0$};
			\braidbox{15}{17}{7}{8}{$w_0$};
			\braidbox{18}{20}{7}{8}{$w_0$};
			\draw[dotted, color=gray] (2.5,2) -- (2.5,3);
			\draw[dotted, color=gray] (5.75,2) -- (5.75,3);
			\draw[dotted, color=gray] (11.25,2) -- (11.25,3);
			\draw[dotted, color=gray] (14.5,2) -- (14.5,3);
			\draw[dotted, color=gray] (17.5,2) -- (17.5,3);
			\draw (1,8.5) node{\tiny $2^{b-r}$};
			\draw (5.5,8.5) node{\tiny $1^{r-n}$};
			\draw (11.5,8.5) node{\tiny $2^r$};
			\draw (16,8.5) node{\tiny $1^n$};
			\draw (19,8.5) node{\tiny $1^{a-r}$};
		\end{braid}
	\end{align*}
	where the first equality follows from Lemma~\ref{LDivPower}, the second follows because $f$ is symmetric in the variables as indicated by the vertical dotted lines, hence, by Theorem~\ref{TNHCenter}, commutes with the parabolic subalgebra $\NH_{b-r}\otimes \NH_{r-s}\otimes\NH_{s}\otimes \NH_{r-s-n}\otimes \NH_n \otimes \NH_{a-r}$, and the last equality is straightforward.
	Now we have $e_{r,n}\Xi_n(f) = \Xi_n(f)$ by Lemma~\ref{LDivPower}.

	To complete the proof, in view of Lemma~\ref{LGDim}, we assume that $\k$ is a field and check that $\qdim K_n = \qdim q^{-\tts_{r,n}}e_{r,n}\Delta(s)$.
	For brevity, we denote $[M] := \qdim M$.
	Let $z = (r-s)s + (r-s-n)s + (r-s)(r-s-n)$.
	We have
	\begin{align*}
		\left[1_{\bj_{r,n}} \Delta(s)\right] &= q^{z} \qbinom{r-n}{s}_- \qbinom{r}{s}_- \left[\Delta(\beta^{b-s})\right] \left[1_{1^s 2^s}\Delta(\gamma^s)\right] \left[\Delta(\alpha^{a-s})\right] \\
		&= q^z \qbinom{r-n}{s}_- \qbinom{r}{s}_- ([s]_+^!)^2 \left[\Delta(\alpha^{a-s})\right] \left[1_{1^{(s)} 2^{(s)}}\Delta(\gamma^s)\right]\left[\Delta(\alpha^{a-s})\right] \\
		&= q^{\omega_n(r,s)+\tts_{r,n}} \frac{[r-n]_+^![r]_+^!}{[r-n-s]_+^![r-s]_+^!} [\CX_{b-s}] [\Lambda_s][\CX_{a-s}],
	\end{align*}
	where the first equality follows from (\ref{EShuffle}), then second from (\ref{EUndivided}), and the last from Lemmas~\ref{LSimpleBasis} and~\ref{L12Basis}.
	Thus,
	\begin{align*}
		q^{-\tts_{r,n}}\left[e_{r,n}\Delta(s)\right] &= \frac{q^{-\tts_{r,n}}}{[b-r]_+^![r-n]_+^![r]_+^![n]_+^![a-r]_+^!}\left[1_{\bj_{r,n}}\Delta(s)\right] \\
		&= q^{\omega_n(r,s)}\frac{[\CX_{b-s}]}{[b-r]_+^![r-s]_+^!}[\Lambda_s]\frac{[\CX_{a-s}]}{[r-s-n]_+^![n]_+^![a-r]_+^!} \\
		&= q^{\omega_n(r,s)}[\Lambda_{b-r}][\Lambda_{r-s}][\Lambda_{s}][\Lambda_{r-s-n}][\Lambda_{n}][\Lambda_{a-r}],
	\end{align*}
	where the first equality is (\ref{EUndivided}), the second is by the above computation, and the last is by an elementary computation.
\end{proof}

Because of Lemma~\ref{LWeightBasis}(\ref{LWeightBasisi}), we assume for the rest of the subsection that $r \geq s$.
We use Lemma~\ref{LWeightBasis}(\ref{LWeightBasisii}) to understand the complex $T_\bullet^r(s)$ and compute its cohomology.
First, we re-express the coboundary map of $T_\bullet^r(s)$.
For $0 \leq n < r-s$, set
$$g_n := \prod_{k=1}^{r-s}(x_{b+r-s-n} - x_{b-r+k}) \in \CX_{a+b-s}.$$
If $f \in K_n$, observe that $fg_n \in \Lambda_{b-r, r-s, s, r-s-n-1, 1, n, a-r}$ so by Proposition~\ref{PThickPartial}, we have a map
$$\delta_n: K_n \to K_{n+1}, \ f \mapsto \partial_{U_{b+r-s-n-1;1,n}}(fg_n).$$
Recalling the isomorphisms $\Xi_n:K_n\iso T_n^r(s)$ from  Lemma~\ref{LWeightBasis}(\ref{LWeightBasisii}), we have:

\begin{Lemma}\label{LSymDiff}
	If $0 \leq n < r-s$, then the following diagram commutes:
	$$\begin{tikzcd}
		K_{n+1} \arrow[d,"\Xi_{n+1}"] & K_n \arrow[l,"{\delta_n}" above] \arrow[d, "\Xi_n"] \\
		T_{n+1}^r(s) & T_n^r(s) \arrow[l, "d_{r,n}\cdot -" above].
	\end{tikzcd}$$
	In particular, the maps $\delta_n$ make $K_\bullet$ into a complex isomorphic to $T_\bullet^r(s)$.
\end{Lemma}
\begin{proof}
	There exist polynomials $h_j \in \k[y_{b-r+1}, \ldots, y_{b-s}]$ and $k_j \in \k[y_{b+r-n}]$ such that $\prod_{i=1}^{r-s} (y_{b+r-n} - y_{b-r+i}) = \sum_j h_j k_j$.
	Let $f\in K_n$ and note that
	$$d_{r,n}\Xi_n(f) = e_{r,n+1}\psi_{x_{r,n}} e_{r,n}\Xi_n(f) = e_{r,n+1}\psi_{x_{r,n}}\Xi_n(f)$$
	since $\im \Xi_n\subseteq T_n^r(s)=q^{-\tts_{r,n}}e_{r,n}\Delta(s)$ by Lemma~\ref{LWeightBasis}.
	We compute $e_{r,n+1}\psi_{x_{r,n}}\Xi_n(f)$:
	\begin{align*}
		&
		\begin{braid}\tikzset{baseline=5.5em}
			\draw (0,2) -- (0,13);
			\draw (2,2) -- (2,13);
			\draw (3,2) -- (3,6) -- (15,10) -- (14,13);
			\draw (5,2) -- (5,6) -- (17,10) -- (16,13);
			\draw (6.5,2) -- (6.5,3) -- (6,4) -- (6,5);
			\draw (12.5,2) -- (12.5,3) -- (12,4) -- (9,5);
			\draw (6.5,2) -- (6.5,3) -- (7,4) -- (10,5);
			\draw (12.5,2) -- (12.5,3) -- (13,4) -- (13,5);
			\draw (6,6) -- (3,8) -- (7,10) -- (7,13);
			\draw (8,6) -- (5,8) -- (9,10) -- (9,13);
			\draw (9,6) -- (6,8) -- (10,10) -- (20,13);
			\draw (10,6) -- (14,8) -- (11,10) -- (10,13);
			\draw (13,6) -- (17,8) -- (14,10) -- (13,13);
			\draw (14,2) -- (14,6) -- (3,10) -- (3,13);
			\draw (17,2) -- (17,6) -- (6,10) -- (6,13);
			\draw (18,2) -- (18,10) -- (17,13);
			\draw (20,2) -- (20,10) -- (19,13);
			\draw (21,2) -- (21,13);
			\draw (23,2) -- (23,13);
			\braidbox{0}{5}{0}{1}{$w_0$}
			\braidbox{14}{23}{0}{1}{$w_0$}
			\braidbox{0}{23}{1}{2}{\scriptsize $f$}
			\braidbox{6}{9}{5}{6}{$w_0$}
			\braidbox{10}{13}{5}{6}{$w_0$}
			\braidbox[0.5]{0}{2}{13}{14}{\tiny $2^{b-r}$}
			\braidbox{3}{9}{13}{14}{\tiny $1^{r-n-1}$}
			\braidbox{10}{16}{13}{14}{\tiny $2^{r}$}
			\braidbox{17}{20}{13}{14}{\tiny $1^{n+1}$}
			\braidbox[0.5]{21}{23}{13}{14}{\tiny $1^{a-r}$}
			\draw (0,7) node[left]{\normalsize $\phantom{=}$};
		\end{braid}
		& \\
		&
		\begin{braid}\tikzset{baseline=5.5em}
			\draw (0,2) -- (0,14);
			\draw (2,2) -- (2,14);
			\draw (3,2) -- (3,6) -- (15,10) -- (14,14);
			\draw (5,2) -- (5,6) -- (17,10) -- (16,14);
			\draw (6.5,2) -- (6.5,3) -- (6,4) -- (6,5) -- (9,6) -- (6,8) -- (10,10) -- (20,14);
			\draw (8.5,2) -- (8.5,3) -- (8,4) -- (7,5) -- (6,6) -- (3,8) -- (7,10) -- (7,13);
			\draw (12.5,2) -- (12.5,3) -- (12,4) -- (9,5) -- (8,6) -- (5,8) -- (9,10) -- (9,13);
			\draw (6.5,2) -- (6.5,3) -- (7,4) -- (10,5) -- (14,8) -- (11,10) -- (10,13);
			\draw (8.5,2) -- (8.5,3) -- (9,4) -- (11,5) -- (15,8) -- (12,10) -- (11,13);
			\draw (12.5,2) -- (12.5,3) -- (13,4) -- (13,5) -- (17,8) -- (14,10) -- (13,13);
			\draw (14,2) -- (14,6) -- (3,10) -- (3,14);
			\draw (17,2) -- (17,6) -- (6,10) -- (6,14);
			\draw (18,2) -- (18,10) -- (17,14);
			\draw (20,2) -- (20,10) -- (19,14);
			\draw (21,2) -- (21,14);
			\draw (23,2) -- (23,14);
			\braidbox{0}{5}{0}{1}{$w_0$}
			\braidbox{14}{23}{0}{1}{$w_0$}
			\braidbox{0}{23}{1}{2}{\scriptsize $f$}
			\braidbox{7}{9}{13}{14}{$w_0$}
			\braidbox{10}{13}{13}{14}{$w_0$}
			\draw (1,14.5) node{\tiny $2^{b-r}$};
			\draw (4.5,14.5) node{\tiny $1^{r-s-n}$};
			\draw (8,14.5) node{\tiny $1^{s-1}$};
			\draw (11.5, 14.5) node{\tiny $2^{s}$};
			\draw (15,14.5) node{\tiny $2^{r-s}$};
			\draw (18,14.5) node{\tiny $1^{n}$};
			\draw (20,14.5) node{\tiny $1$};
			\draw (22,14.5) node{\tiny $1^{a-r}$};
			\draw (0,7) node[left]{\normalsize $=$};
		\end{braid}
		& \begin{tabular}{p{11em}} by the relation (\ref{Rbraid}) and Lemma~\ref{LDivPower} \end{tabular} \\
		&
		\begin{braid}\tikzset{baseline=5.5em}
			\draw (0,2) -- (0,13);
			\draw (2,2) -- (2,13);
			\draw (3,2) -- (3,5) -- (11,8) -- (11,10) -- (14,12) -- (14,13);
			\draw (5,2) -- (5,5) -- (13,8) -- (13,10) -- (16,12) -- (16,13);
			\draw (6.5,2) -- (6.5,3) -- (6,4) -- (6,5) -- (12,7) -- (10,8) -- (20,12) -- (20,13);
			\draw (8.5,2) -- (8.5,3) -- (8,4) -- (7,5) -- (3,7) -- (3,10) -- (7,12);
			\draw (12.5,2) -- (12.5,3) -- (12,4) -- (9,5) -- (5,7) -- (5,10) -- (9,12);
			\draw (6.5,2) -- (6.5,3) -- (7,4) -- (10,5) -- (14,8) -- (14,10) -- (10,12);
			\draw (8.5,2) -- (8.5,3) -- (9,4) -- (11,5) -- (15,8) -- (15,10) -- (11,12);
			\draw (12.5,2) -- (12.5,3) -- (13,4) -- (13,5) -- (17,8) -- (17,10) -- (13,12);
			\draw (14,2) -- (14,5) -- (6,7) -- (6,10) -- (3,12) -- (3,13);
			\draw (17,2) -- (17,5) -- (9,7) -- (9,10) -- (6,12) -- (6,13);
			\draw (18,2) -- (18,11) -- (17,12) -- (17,13);
			\draw (20,2) -- (20,11) -- (19,12) -- (19,13);
			\draw (21,2) -- (21,13);
			\draw (23,2) -- (23,13);
			\braidbox{0}{5}{0}{1}{$w_0$}
			\braidbox{14}{23}{0}{1}{$w_0$}
			\braidbox{0}{23}{1}{2}{\scriptsize $f$}
			\braidbox{7}{9}{12}{13}{$w_0$}
			\braidbox{10}{13}{12}{13}{$w_0$}
			\draw (1,13.5) node{\tiny $2^{b-r}$};
			\draw (4.5,13.5) node{\tiny $1^{r-s-n}$};
			\draw (8,13.5) node{\tiny $1^{s-1}$};
			\draw (11.5, 13.5) node{\tiny $2^{s}$};
			\draw (15,13.5) node{\tiny $2^{r-s}$};
			\draw (18,13.5) node{\tiny $1^{n}$};
			\draw (20,13.5) node{\tiny $1$};
			\draw (22,13.5) node{\tiny $1^{a-r}$};
			\draw (0,6.5) node[left]{\normalsize $=$};
		\end{braid}
		& \begin{tabular}{p{11em}} by the relation (\ref{Rbraid}) \end{tabular} \\
		&
		\begin{braid}\tikzset{baseline=4.5em}
			\draw (0,2) -- (0,11);
			\draw (2,2) -- (2,11);
			\draw (3,3) -- (3,5) -- (10,8) -- (14,10) -- (14,11);
			\draw (5,3) -- (5,5) -- (12,8) -- (16,10) -- (16,11);
			\draw (6.5,2) -- (6.5,3) -- (6,4) -- (6,5) -- (17,9);
			\draw (8.5,2) -- (8.5,3) -- (8,4) -- (7,5) -- (3,8) -- (7,10);
			\draw (12.5,2) -- (12.5,3) -- (12,4) -- (9,5) -- (5,8) -- (9,10);
			\draw (6.5,2) -- (6.5,3) coordinate (x);
			\draw [dashed] (x) -- (7,4) -- (10,5) -- (14,7) -- (13,8) coordinate (x);
			\draw (x) -- (10,10);
			\draw (8.5,2) -- (8.5,3) -- (9,4) -- (11,5) -- (15,7) -- (14,8) -- (11,10);
			\draw (12.5,2) -- (12.5,3) -- (13,4) -- (13,5) -- (17,7) -- (16,8) -- (13,10);
			\draw (14,2) -- (14,5) -- (6,8) -- (3,10) -- (3,11);
			\draw (17,2) -- (17,5) -- (9,8) -- (6,10) -- (6,11);
			\draw (18,2) -- (18,10) -- (17,11);
			\draw (20,2) -- (20,10) -- (19,11);
			\draw (21,2) -- (21,11);
			\draw (23,2) -- (23,11);
			\draw (17,10) -- (20,11);
			\braidbox{0}{5}{0}{1}{$w_0$}
			\braidbox{14}{23}{0}{1}{$w_0$}
			\braidbox{0}{23}{1}{2}{\scriptsize $f$}
			\braidbox{3}{5}{2}{3}{\scriptsize $h_i$}
			\braidbox{7}{9}{10}{11}{$w_0$}
			\braidbox{10}{13}{10}{11}{$w_0$}
			\braidbox{16.6}{17.4}{9}{10}{\scriptsize $k_i$}
			\draw (1,11.5) node{\tiny $2^{b-r}$};
			\draw (4.5,11.5) node{\tiny $1^{r-s-n}$};
			\draw (8,11.5) node{\tiny $1^{s-1}$};
			\draw (11.5, 11.5) node{\tiny $2^{s}$};
			\draw (15,11.5) node{\tiny $2^{r-s}$};
			\draw (18,11.5) node{\tiny $1^{n}$};
			\draw (20,11.5) node{\tiny $1$};
			\draw (22,11.5) node{\tiny $1^{a-r}$};
			\draw (0,5.5) node[left]{\normalsize $= \sum\limits_j$};
		\end{braid}
		& \begin{tabular}{p{11em}} by several applications of the relation (\ref{Rquad}) \end{tabular} \\
		&
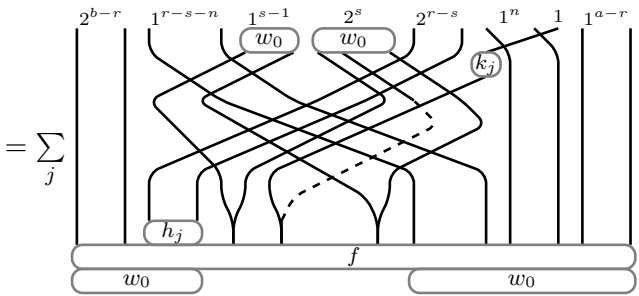
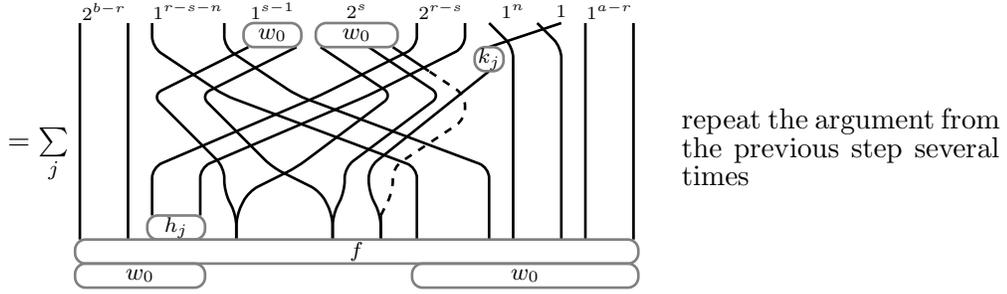
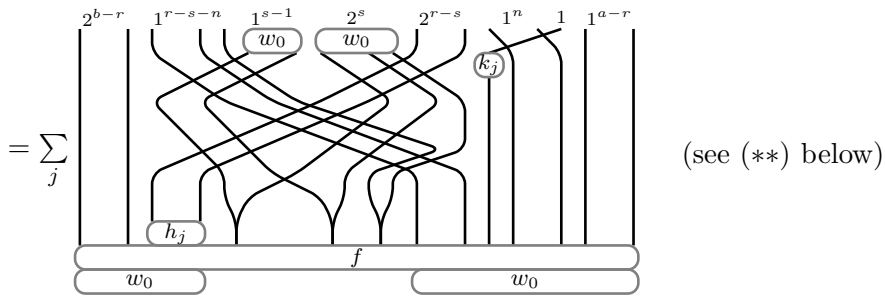
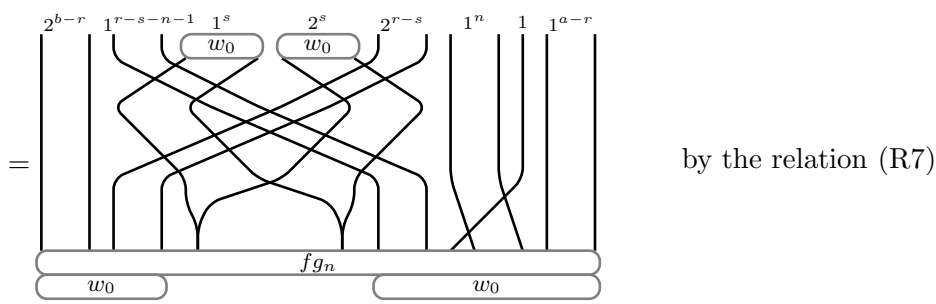
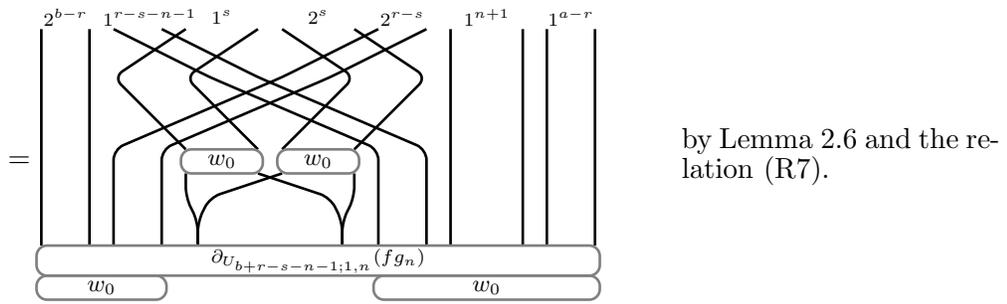
		\begin{braid}\tikzset{baseline=4.5em}
			\draw (0,2) -- (0,11);
			\draw (2,2) -- (2,11);
			\draw (3,3) -- (3,5) -- (10,8) -- (14,10) -- (14,11);
			\draw (5,3) -- (5,5) -- (12,8) -- (16,10) -- (16,11);
			\draw (6.5,2) -- (6.5,3) -- (6,4) -- (6,5) -- (3,8) -- (7,10);
			\draw (8.5,2) -- (8.5,3) -- (8,4) -- (8,5) -- (17,9);
			\draw (12.5,2) -- (12.5,3) -- (12,4) -- (5,8) -- (9,10);
			\draw (6.5,2) -- (6.5,3) -- (7,4) -- (7,5) -- (13,8) -- (10,10);
			\draw (8.5,2) -- (8.5,3) coordinate (x);
			\draw [dashed] (x) -- (9,4) -- (15,7) -- (14,8) coordinate (x);
			\draw (x) -- (11,10);
			\draw (12.5,2) -- (12.5,3) -- (13,4) -- (13,5) -- (17,7) -- (16,8) -- (13,10);
			\draw (14,2) -- (14,5) -- (6,8) -- (3,10) -- (3,11);
			\draw (17,2) -- (17,5) -- (9,8) -- (6,10) -- (6,11);
			\draw (18,2) -- (18,10) -- (17,11);
			\draw (20,2) -- (20,10) -- (19,11);
			\draw (21,2) -- (21,11);
			\draw (23,2) -- (23,11);
			\draw (17,10) -- (20,11);
			\braidbox{0}{5}{0}{1}{$w_0$}
			\braidbox{14}{23}{0}{1}{$w_0$}
			\braidbox{0}{23}{1}{2}{\scriptsize $f$}
			\braidbox{3}{5}{2}{3}{\scriptsize $h_j$}
			\braidbox{7}{9}{10}{11}{$w_0$}
			\braidbox{10}{13}{10}{11}{$w_0$}
			\braidbox{16.6}{17.4}{9}{10}{\scriptsize $k_j$}
			\draw (1,11.5) node{\tiny $2^{b-r}$};
			\draw (4.5,11.5) node{\tiny $1^{r-s-n}$};
			\draw (8,11.5) node{\tiny $1^{s-1}$};
			\draw (11.5, 11.5) node{\tiny $2^{s}$};
			\draw (15,11.5) node{\tiny $2^{r-s}$};
			\draw (18,11.5) node{\tiny $1^{n}$};
			\draw (20,11.5) node{\tiny $1$};
			\draw (22,11.5) node{\tiny $1^{a-r}$};
			\draw (0,5.5) node[left]{\normalsize $= \sum\limits_j$};
		\end{braid}
		& \begin{tabular}{p{11em}} (see $(*)$ below) \end{tabular} \\
		&
		\begin{braid}\tikzset{baseline=4.5em}
			\draw (0,2) -- (0,11);
			\draw (2,2) -- (2,11);
			\draw (3,3) -- (3,5) -- (10,8) -- (14,10) -- (14,11);
			\draw (5,3) -- (5,5) -- (12,8) -- (16,10) -- (16,11);
			\draw (6.5,2) -- (6.5,3) -- (6,4) -- (6,5) -- (3,8) -- (7,10);
			\draw (10.5,2) -- (10.5,3) -- (10,4) -- (8,5) -- (5,8) -- (9,10);
			\draw (12.5,2) -- (12.5,3) -- (12,4) -- (12,5) -- (17,9);
			\draw (6.5,2) -- (6.5,3) -- (7,4) -- (9,5) -- (13,8) -- (10,10);
			\draw (10.5,2) -- (10.5,3) -- (11,4) -- (11,5) -- (15,8) -- (12,10);
			\draw (12.5,2) -- (12.5,3) coordinate (x);
			\draw [dashed] (x) -- (13,4) -- (13,5) -- (16,7) -- (16,8) -- (14.5,9) coordinate (x);
			\draw (x) -- (13,10);
			\draw (14,2) -- (14,5) -- (6,8) -- (3,10) -- (3,11);
			\draw (17,2) -- (17,5) -- (9,8) -- (6,10) -- (6,11);
			\draw (18,2) -- (18,10) -- (17,11);
			\draw (20,2) -- (20,10) -- (19,11);
			\draw (21,2) -- (21,11);
			\draw (23,2) -- (23,11);
			\draw (17,10) -- (20,11);
			\braidbox{0}{5}{0}{1}{$w_0$}
			\braidbox{14}{23}{0}{1}{$w_0$}
			\braidbox{0}{23}{1}{2}{\scriptsize $f$}
			\braidbox{3}{5}{2}{3}{\scriptsize $h_j$}
			\braidbox{7}{9}{10}{11}{$w_0$}
			\braidbox{10}{13}{10}{11}{$w_0$}
			\braidbox{16.6}{17.4}{9}{10}{\scriptsize $k_j$}
			\draw (1,11.5) node{\tiny $2^{b-r}$};
			\draw (4.5,11.5) node{\tiny $1^{r-s-n}$};
			\draw (8,11.5) node{\tiny $1^{s-1}$};
			\draw (11.5, 11.5) node{\tiny $2^{s}$};
			\draw (15,11.5) node{\tiny $2^{r-s}$};
			\draw (18,11.5) node{\tiny $1^{n}$};
			\draw (20,11.5) node{\tiny $1$};
			\draw (22,11.5) node{\tiny $1^{a-r}$};
			\draw (0,5.5) node[left]{\normalsize $= \sum\limits_j$};
		\end{braid}
		& \begin{tabular}{p{11em}} repeat the argument from the previous step several times \end{tabular} \\
		&
		\begin{braid}\tikzset{baseline=4.5em}
			\draw (0,2) -- (0,11);
			\draw (2,2) -- (2,11);
			\draw (3,3) -- (3,5) -- (10,8) -- (14,10) -- (14,11);
			\draw (5,3) -- (5,5) -- (12,8) -- (16,10) -- (16,11);
			\draw (6.5,2) -- (6.5,3) -- (6,4) -- (6,5) -- (3,8) -- (7,10);
			\draw (10.5,2) -- (10.5,3) -- (10,4) -- (8,5) -- (5,8) -- (9,10);
			\draw (12.5,2) -- (12.5,3) -- (12,4) -- (12,5) -- (15,6) -- (9,8) -- (6,10) -- (6,11);
			\draw (6.5,2) -- (6.5,3) -- (7,4) -- (9,5) -- (13,8) -- (10,10);
			\draw (10.5,2) -- (10.5,3) -- (11,4) -- (11,5) -- (15,8) -- (12,10);
			\draw (12.5,2) -- (12.5,3) -- (13,4) -- (13,5) -- (16,6) -- (16,8) -- (13,10);
			\draw (14,2) -- (14,5) -- (6,8) -- (3,10) -- (3,11);
			\draw (16,2) -- (16,5) -- (8,8) -- (5,10) -- (5,11);
			\draw (17,2) -- (17,9);
			\draw (18,2) -- (18,10) -- (17,11);
			\draw (20,2) -- (20,10) -- (19,11);
			\draw (21,2) -- (21,11);
			\draw (23,2) -- (23,11);
			\draw (17,10) -- (20,11);
			\braidbox{0}{5}{0}{1}{$w_0$}
			\braidbox{14}{23}{0}{1}{$w_0$}
			\braidbox{0}{23}{1}{2}{\scriptsize $f$}
			\braidbox{3}{5}{2}{3}{\scriptsize $h_j$}
			\braidbox{7}{9}{10}{11}{$w_0$}
			\braidbox{10}{13}{10}{11}{$w_0$}
			\braidbox{16.6}{17.4}{9}{10}{\scriptsize $k_j$}
			\draw (1,11.5) node{\tiny $2^{b-r}$};
			\draw (4.5,11.5) node{\tiny $1^{r-s-n}$};
			\draw (8,11.5) node{\tiny $1^{s-1}$};
			\draw (11.5, 11.5) node{\tiny $2^{s}$};
			\draw (15,11.5) node{\tiny $2^{r-s}$};
			\draw (18,11.5) node{\tiny $1^{n}$};
			\draw (20,11.5) node{\tiny $1$};
			\draw (22,11.5) node{\tiny $1^{a-r}$};
			\draw (0,5.5) node[left]{\normalsize $= \sum\limits_j$};
		\end{braid}
		& \begin{tabular}{p{11em}} (see $(**)$ below) \end{tabular} \\
		&
		\begin{braid}\tikzset{baseline=4.5em}
			\draw (0,2) -- (0,11);
			\draw (2,2) -- (2,11);
			\draw (3,2) -- (3,5) -- (10,8) -- (14,10) -- (14,11);
			\draw (5,2) -- (5,5) -- (12,8) -- (16,10) -- (16,11);
			\draw (6.5,2) -- (6.5,3) -- (6,4) -- (6,5) -- (3,8) -- (6,10);
			\draw (12.5,2) -- (12.5,3) -- (12,4) -- (9,5) -- (6,8) -- (9,10);
			\draw (6.5,2) -- (6.5,3) -- (7,4) -- (10,5) -- (13,8) -- (10,10);
			\draw (12.5,2) -- (12.5,3) -- (13,4) -- (13,5) -- (16,8) -- (13,10);
			\draw (14,2) -- (14,5) -- (7,8) -- (3,10) -- (3,11);
			\draw (16,2) -- (16,5) -- (9,8) -- (5,10) -- (5,11);
			\draw (17,2) -- (20,5) -- (20,11);
			\draw (18,2) -- (17,5) -- (17,11);
			\draw (20,2) -- (19,5) -- (19,11);
			\draw (21,2) -- (21,11);
			\draw (23,2) -- (23,11);
			\braidbox{0}{5}{0}{1}{$w_0$}
			\braidbox{14}{23}{0}{1}{$w_0$}
			\braidbox{0}{23}{1}{2}{\scriptsize $fg_n$}
			\braidbox{6}{9}{10}{11}{$w_0$}
			\braidbox{10}{13}{10}{11}{$w_0$}
			\draw (1,11.5) node{\tiny $2^{b-r}$};
			\draw (4.5,11.5) node{\tiny $1^{r-s-n-1}$};
			\draw (7.5,11.5) node{\tiny $1^{s}$};
			\draw (11.5, 11.5) node{\tiny $2^{s}$};
			\draw (15,11.5) node{\tiny $2^{r-s}$};
			\draw (18,11.5) node{\tiny $1^{n}$};
			\draw (20,11.5) node{\tiny $1$};
			\draw (22,11.5) node{\tiny $1^{a-r}$};
			\draw (0,5.5) node[left]{\normalsize $=$};
		\end{braid}
		& \begin{tabular}{p{11em}} by the relation (\ref{Rbraid}) \end{tabular} \\
		&
		\begin{braid}\tikzset{baseline=4.5em}
			\draw (0,2) -- (0,11);
			\draw (2,2) -- (2,11);
			\draw (3,2) -- (3,6) -- (10,9) -- (14,11);
			\draw (5,2) -- (5,6) -- (12,9) -- (16,11);
			\draw (6.5,2) -- (6.5,3) -- (6,4) -- (6,5);
			\draw (12.5,2) -- (12.5,3) -- (12,4) -- (9,5);
			\draw (6.5,2) -- (6.5,3) -- (7,4) -- (10,5);
			\draw (12.5,2) -- (12.5,3) -- (13,4) -- (13,5);
			\draw (14,2) -- (14,6) -- (7,9) -- (3,11);
			\draw (16,2) -- (16,6) -- (9,9) -- (5,11);
			\draw (17,2) -- (17,11);
			\draw (20,2) -- (20,11);
			\draw (21,2) -- (21,11);
			\draw (23,2) -- (23,11);
			\draw (6,6) -- (3,9) -- (6,11);
			\draw (9,6) -- (6,9) -- (9,11);
			\draw (10,6) -- (13,9) -- (10,11);
			\draw (13,6) -- (16,9) -- (13,11);
			\braidbox{0}{5}{-0.25}{0.75}{$w_0$}
			\braidbox{14}{23}{-0.25}{0.75}{$w_0$}
			\braidbox{0}{23}{0.75}{2}{\tiny $\partial_{U_{b+r-s-n-1;1,n}}(fg_n)$}
			\braidbox{6}{9}{5}{6}{$w_0$}
			\braidbox{10}{13}{5}{6}{$w_0$}
			\draw (1,11.5) node{\tiny $2^{b-r}$};
			\draw (4.5,11.5) node{\tiny $1^{r-s-n-1}$};
			\draw (7.5,11.5) node{\tiny $1^{s}$};
			\draw (11.5, 11.5) node{\tiny $2^{s}$};
			\draw (15,11.5) node{\tiny $2^{r-s}$};
			\draw (18.5,11.5) node{\tiny $1^{n+1}$};
			\draw (22,11.5) node{\tiny $1^{a-r}$};
			\draw (0,5.5) node[left]{\normalsize $=$};
		\end{braid}
		& \begin{tabular}{p{11em}} by Lemma~\ref{LPartialPoly} and the relation (\ref{Rbraid}). \end{tabular}
	\end{align*}

	$(*)$ To obtain this equality, we attempt to move the portion of dashed strand in the previous diagram to the left by applying a special case of relation (\ref{Rbraid}):
	\begin{equation}\label{EEasyBraid}
		\psi_{t+1}\psi_t\psi_{t+1}1_{\bi} = \psi_t\psi_{t+1}\psi_t 1_{\bi} + 1_\bi \ (\text{if } i_{t+1} = i_t+1 \text{ and } i_t = i_{t+2})
	\end{equation}
	several times.
	In all except the last application, the error term $1_\bi$ causes the rest of the diagram to become $0$, so we only keep the term $\psi_t\psi_{t+1}\psi_t 1_{\bi}$.
	In the last application, because of the defining relations in $\Delta(s)$, the term $\psi_t\psi_{t+1}\psi_t 1_{\bi}$ causes the rest of the diagram to become $0$, so we only keep the error term $1_\bi$, which yields the desired diagram.

	$(**)$ To obtain this equality, we again attempt to move the portion of dashed strand in the previous diagram to the left by applying the relation (\ref{EEasyBraid}) several times.
	In the first application, the error term $1_\bi$ yields the desired diagram, so we wish to show that the term $\psi_t\psi_{t+1}\psi_t 1_{\bi}$ causes the rest of the diagram to become $0$.
	Since $f$ is symmetric in the variables $x_{b+1}, \ldots, x_{b+r-s-n}$, the error term in any application of the relation (\ref{EEasyBraid}) other than the first causes the rest of the diagram to become $0$.
	However, after the last application of the relation, the term $\psi_t\psi_{t+1}\psi_t 1_{\bi}$ also causes the rest of the diagram to become $0$ because of the defining relations in $\Delta(s)$.

	The expression represented by the last diagram above is $\Xi_{n+1}(\delta_n(f))$, which completes the proof of the lemma.
\end{proof}

Given an interval $(c,d]$ and a polynomial in $d-c$ variables, we denote
$$f(\underline{x}_{(c,d]}):=f(x_{c+1},\dots,x_d).$$
For example, if $0\leq m\leq d-c$ then we have the $m\mathrm{th}$ elementary symmetric function
$$E_m(\underline{x}_{(c,d]})=\sum_{c<i_1<\dots<i_m\leq d}x_{i_1}\cdots x_{i_m}.$$
Now, for $0\leq k < r-s$, we define
$$z_k := (-1)^{r-s-k}(E_{r-s-k}(\underline{x}_{(b-r, b-s]}) - E_{r-s-k}(\underline{x}_{(b, b+r-s]})).$$
These elements are considered as elements of the algebra
\begin{equation}\label{EPartiallyIdentified}
	\Lambda^{r,s} := \Lambda_{b-r, r-s, s, r-s, a-r}.
\end{equation}
Note that since $\Lambda^{r,s}\subseteq \Lambda_{b-r,r-s,s,r-s-n,n,a-r}$, each $K_n$ is naturally a (right) $\Lambda^{r,s}$-module.
We use this to interpret the right-hand side of the lemma below as an element of $K_{n+1}$.

\begin{Lemma}\label{LKDiff}
	For $0\leq n < r-s$ and a symmetric polynomial $f$ in $n$ variables, we have
	$$\delta_n(f(\underline{x}_{(b+r-s-n, b+r-s]})) = \sum_{k=0}^{r-s-1} \left(x_{b+r-s-n}^k \star f(\underline{x}_{(b+r-s-n, b+r-s]})\right)z_k.$$
\end{Lemma}
\begin{proof}
	For brevity, write $f = f(\underline{x}_{(b+r-s-n, b+r-s]})$.
	We observe
	$$g_n = \prod_{k=1}^{r-s}(x_{b+r-s-n} - x_{b-r+k}) = \sum_{k=0}^{r-s} (-1)^{r-s-k}x_{b+r-s-n}^k E_{r-s-k}(\underline{x}_{(b-r, b-s]}),$$
	so that
	\begin{equation}\label{EDiffSchur}
		\delta_n(f) = \sum_{k=0}^{r-s} (-1)^{r-s-k} \partial_{U_{b+r-s-n-1;1,n}}(x_{b+r-s-n}^k f)E_{r-s-k}(\underline{x}_{(b-r, b-s]}).
	\end{equation}
	Next, using the identity
	$$x_{b+r-s-n}^{r-s} = -\sum_{k=0}^{r-s-1}(-1)^{r-s-k}x_{b+r-s-n}^k E_{r-s-k}(\underline{x}_{(b, b+r-s]}),$$
	(\ref{EDiffSchur}) becomes
	$$\delta_n(f) = \sum_{k=0}^{r-s-1}\partial_{U_{b+r-s-n-1;1,n}}(x_{b+r-s-n}^k f)z_k,$$
	and the result follows from Proposition~\ref{PThickPartial}.
\end{proof}

The proof of Theorem~\ref{TA2Ext} below amounts to showing that $T_\bullet^\rho(\sigma) \cong K_\bullet$ is isomorphic to a certain Koszul complex which we now define.
Let $N$ be the free right $\Lambda^{r,s}$-module of graded rank $\sum_{k=0}^{r-s-1}q^{2k-2(r-s)}$.
For $k=0,\dots,r-s-1$, we have the basis element $\epsilon_k := 1\in q^{2k-2(r-s)}\Lambda^{r,s}\subseteq N$.
We set
$$Z := (z_0, \ldots, z_{r-s-1}) =  \epsilon_0 z_0 + \cdots + \epsilon_{r-s-1} z_{r-s-1}\in N.$$
Note that $Z$ is a homogeneous degree $0$ element of $N$.
We consider the Koszul complex $q^{\omega_0(r,s)}\Wed^{\bullet}N$ associated to the regular sequence $Z$ for the algebra $\Lambda^{r,s}$ (see~\cite[\S4.5]{We}), which has the form
\begin{align*}
	\cdots \longleftarrow q^{\omega_0(r,s)}\Wed^{n+1} N &\longleftarrow q^{\omega_0(r,s)}\Wed^n N \longleftarrow \cdots \\
	Z\wedge a &\longmapsfrom a
\end{align*}
where $\Wed^n N$ is the $n\rm{th}$ exterior power of the free $\Lambda^{r,s}$-module $N$.
Note that $\Wed^n N$ has basis $\{ \epsilon_{i_1}\wedge \cdots \wedge \epsilon_{i_n} \mid 0 \leq i_1 < \cdots < i_n < r-s\}$.
Recall (\ref{EStar}).

By Lemmas~\ref{LCochain} and~\ref{LWeightBasis}(\ref{LWeightBasisi}), the complex $\Hom_{R_\theta}^\bullet(P_\bullet^r, \Delta(s))$ is zero in degrees larger than $r-s$, so every element of the $r-s$ component is a cocycle and there is a surjective map
$$[-]: \Hom_{R_\theta}(q^{\tts_{r,r-s}} R_\theta e_{r,r-s}, \Delta(s)) \onto \CE_\theta^{r-s}(r, s) = H^{r-s}(\Hom_{R_\theta}^\bullet(P_\bullet^r, \Delta(s))), \ \varphi \mapsto [\varphi],$$
where $[\varphi]$ is the cohomology class of $\varphi$.
Moreover, by Lemma~\ref{LWeightBasis}(\ref{LWeightBasisii}), we have an isomorphism
$$\xi: K_{r-s} \iso \Hom_{R_\theta}(q^{\tts_{r,r-s}}R_\theta e_{r,r-s},\Delta(s)), \ f \mapsto (e_{r-s,s} \mapsto v_s f).$$
For $r, s \in \Z_{\geq 0}$ with $\min\{ a, b\} \geq r \geq s$, we define
$$\Lambda(r,s) := q^{\omega_{r-s}(r,s)}\Lambda_{b-r, r-s, s, a-r}.$$
Note that there is a surjection
\begin{equation}\label{ELambdaSur}
	\ttp_{r,s}: q^{\omega_{r-s}(r,s)}\Lambda^{r,s} \onto \Lambda(r,s),\ f_1\otimes f_2\otimes f_3\otimes f_4\otimes f_5\mapsto f_1\otimes f_2f_4\otimes f_3\otimes f_5,
\end{equation}
obtained by identifying the two $\Lambda_{r-s}$ components.
Since $\Lambda_{\sigma(s)}=\Lambda_{b-s,s,a-s}\subseteq \Lambda^{r,s}$, we consider $\Lambda^{r,s}$ to be a right $\Lambda_{\sigma(s)}$-module, and we consider $\Lambda(r,s)$ to be a right $\Lambda_{\sigma(s)}$-module via the composition of algebra homomorphisms
$$\Lambda_{\sigma(s)} \longinto \Lambda^{r,s} \overset{\ttp_{r,s}}{\longonto} \Lambda_{b-r, r-s, s, a-r}=q^{-\omega_{r-s}(r,s)}\Lambda(r,s)$$
and then degree shift.
Note that $\ttp_{r,s}$ is a $\Lambda_{\sigma(s)}$-homomorphism.

\begin{Theorem}\label{TA2Ext}
	Let $0 \leq r,s \leq \min\{ a, b\}$.
	If $r < s$, then $\CE_\theta(r, s) = 0$.
	If $r \geq s$, then $\CE_\theta(r,s) = \CE_\theta^{r-s}(r,s)$ and there is an isomorphism of right $\Lambda_{\sigma(s)}$-modules $\CE_\theta^{r-s}(r,s)\iso \Lambda(r,s)$ such that the following diagram of right $\Lambda_{\sigma(s)}$-modules commutes:
	$$\begin{tikzcd}
		q^{\omega_{r,s}(r-s)}\Lambda^{r,s} \arrow[r, "\sim", "\xi" below] \arrow[d, "\ttp_{r,s}"] & \Hom_{R_\theta}(q^{\tts_{r,r-s}}R_\theta e_{r,r-s}, \Delta(s)) \arrow[d, "{[-]}"] \\
		\Lambda(r,s) & \CE_\theta^{r-s}(r,s). \arrow[l, dashed, "\sim"]
	\end{tikzcd}$$
\end{Theorem}
\begin{proof}
	If $r < s$, then $\CE_\theta(r, s) = 0$ by Lemmas~\ref{LCochain} and~\ref{LWeightBasis}(\ref{LWeightBasisi}), so assume $r \geq s$.
	For $0\leq n\leq r-s$ and $\lambda \in \Par(n, r-s-n)$, let  $s_\lambda := s_\lambda(\underline{x}_{(b+r-s-n, b+r-s]}) \in K_n$.
	By Proposition~\ref{PPartSymBasis}, $\{ s_\lambda\mid  \lambda \in \Par(n, r-s-n)\}$ is a basis of $K_n$ as an $\Lambda^{r,s}$-module, so there exists an isomorphism $\Theta_n: K_n \to q^{\omega_0(r,s)}\Wed^n N$ of $\Lambda^{r,s}$-modules such that
	$$\Theta_n(s_\lambda) = \epsilon_{\lambda_n}\wedge \epsilon_{\lambda_{n-1}+1}\wedge \cdots \wedge \epsilon_{\lambda_1+n-1}.$$

	We claim that the maps $\Theta_n$ define an isomorphism of complexes between $K_\bullet$ and $q^{\omega_0(r,s)}\Wed^\bullet N$.
	We must verify that the following square commutes:
	$$\begin{tikzcd}
		K_{n+1} \arrow[d, "\Theta_{n+1}"] & K_n \arrow[l, "\delta_n" above] \arrow[d, "\Theta_n"] \\
		q^{\omega_0(r,s)}\Wed^{n+1}N & q^{\omega_0(r,s)}\Wed^n N. \arrow[l, "Z\wedge -" above]
	\end{tikzcd}$$
	Fix some $\lambda = (\lambda_1, \ldots, \lambda_n)\in \Par(n, r-s-n)$ and set $X := \{\lambda_n, \lambda_{n-1}+1, \ldots, \lambda_1+n-1\}$.
	We then have by Lemma~\ref{LKDiff} and (\ref{EStar}):
	\begin{align*}
		\Theta_{n+1}(\delta_n(s_\lambda)) &= \sum_{k\in [0,r-s)\setminus X}(-1)^{|X \cap [0,k)|} (\epsilon_{\lambda_n}\wedge\cdots \wedge \epsilon_k\wedge \cdots \wedge \epsilon_{\lambda_1+n-1} )z_k\\
		&= \sum_{k\in [0,r-s)\setminus X} (\epsilon_k\wedge \epsilon_{\lambda_n}\wedge \cdots \wedge \epsilon_{\lambda_1+n-1})z_k \\
		&= Z \wedge \Theta_n(s_\lambda).
	\end{align*}

	Since $Z$ is a regular sequence, we have $\CE_\theta^n(r, s) \cong H^n(q^{\omega_0(r,s)}\Wed^\bullet N) = 0$ unless $n = r-s$.
	We complete the proof using Lemmas~\ref{LCochain} and~\ref{LSymDiff} and the observation that by the fundamental theorem of elementary symmetric polynomials, the kernel of $\ttp_{r,s}: q^{\omega_{r-s}(r,s)}\Lambda^{r,s} \onto \Lambda(r,s)$ is the ideal generated by $(z_0, \ldots, z_{r-s-1})$, so $\ttp_{r,s}$ induces an isomorphism of $\Lambda^{r,s}$-modules (and therefore of $\Lambda_{\sigma(s)}$-modules)
	\[H^{r-s}(q^{\omega_0(r,s)}\Wed^\bullet N) = q^{\omega_{r-s}(r,s)}\Lambda^{r,s}/(z_0, \ldots, z_{r-s-1}) \iso \Lambda(r,s).\qedhere\]
\end{proof}

\subsection{The category $\CE_\theta$}\label{SSA2Cat}

Throughout this subsection, we use Theorem~\ref{TA2Ext} to identify $\CE_\theta^{r-s}(r, s)$ with $\Lambda(r,s) = q^{\omega_{r-s}(r,s)}\Lambda_{b-r,r-s,s,a-r}$ whenever $\min\{ a, b\} \geq r \geq s \geq 0$.

In this subsection, we will need to consider not only partially symmetric polynomials in the variables $x$ but also partially symmetric polynomials in the variables $y$.
This will be important since elements of $\CY_{a+b}$ will be considered as elements of $R_\theta$, cf. (\ref{EY}).
For any $d$ we have an isomorphism
$$\iota_{y\to x}:\CY_d\iso \CX_d,\ y_r\mapsto x_r.$$
We will use the notation $\Lambda_m^\CY:=\CY_m^{\Sgp_m}$ for the symmetric polynomials in $y_1,\dots,y_m$.
More generally, given a composition $\mu = (\mu_1, \ldots, \mu_k)$ of $d$, we have the algebra of $\mu$-partially symmetric polynomials $\Lambda_\mu^\CY := \CY_d^{\Sgp_\mu}$.
We often write $\Lambda_{\mu_1,\dots,\mu_k}^\CY$ for $\Lambda_\mu^\CY$ and identify it with $\Lambda_{\mu_1}^\CY \otimes \cdots \otimes \Lambda_{\mu_k}^\CY$.
The isomorphism $\iota_{y\to x}$ restricts to the isomorphism $\iota_{y\to x}:\Lambda_\mu^{\CY}\iso\Lambda_\mu$.

For integers $r,s,t$ with $\min\{ a, b\} \geq r\geq s \geq t \geq 0$, define the following:
\begin{align*}
	\hat{\Lambda}(r,s) &:= q^{\omega_{r-s}(r,s)}\Lambda^\CY_{b-r, r-s, s, s, r-s, a-r}\subseteq \CY_{a+b}, \\
	\hat{\Lambda}(r,s,t) &:= q^{\omega_{r-t}(r,t)-4(r-s)(s-t)}\Lambda^\CY_{b-r, r-s, s-t, t, t, s-t, r-s, a-r}\subseteq \CY_{a+b}.
\end{align*}
Recalling (\ref{EPartiallyIdentified}), there is a surjection
$$\hat{\ttp}_{r,s}: \hat{\Lambda}(r,s) \onto q^{\omega_{r-s}(r,s)}\Lambda^{r,s},\ f_1\otimes f_2\otimes f_3\otimes f_4\otimes f_5\otimes f_6\mapsto \iota_{y\to x}\big(f_1\otimes f_2\otimes f_3 f_4\otimes f_5\otimes f_6\big).$$
Recalling (\ref{ELambdaSur}), let
$$\ttq_{r,s} := \ttp_{r,s}\hat{\ttp}_{r,s}: \hat{\Lambda}(r,s) \to \Lambda(r,s).$$

Let again $\min\{ a, b\} \geq r\geq s \geq t \geq 0$, and $0 \leq n \leq r-s$.
Set
\begin{align*}
	D(r,s,n) &:= \prod_{\substack{i\in (b+s-n, b+s] \\ j\in (b-r, b-s]}}(y_i-y_j) \in \CY_{a+b}, \\
	u(s,n) &:= U_{b-n; n, s} \in \Sgp_{a+b}, \\
	v(r,s,n) &:= U_{b-r; r-s, 2s-n}U_{b+s-n; n, r-s} \in \Sgp_{a+b}, \\
	w(s,t) &:= U_{b-s; t, s-t}U_{b; s-t, t}U_{b-s+t; s-t, s-t} \in \Sgp_{a+b}, \\
	x(r,s,t) &:= U_{b-r; r-s, s-t}U_{b+t; s-t, r-s} \in \Sgp_{a+b}.
\end{align*}

Let $\min\{ a, b\} \geq r\geq s \geq n \geq 0$.
For any $\hat{f} \in \hat{\Lambda}(r,s)$, set
\begin{align*}
	\hat{f}_{r,s}^n &:= (-1)^{\binom{r-s+1}{2}+n(r-s)}e_{r,r-s+n}\psi_{v(r,s,n)}D(r,s,n)(u(s,n)\cdot \hat{f})e_{s,n} \\
	&= \pm
	\begin{braid}\tikzset{baseline=3em}
		\draw (0,2) -- (0,6);
		\draw (2,2) -- (2,6);
		\draw (3,2) -- (9,6);
		\draw (5,2) -- (11,6);
		\draw (6,2) -- (3,6);
		\draw (8,2) -- (5,6);
		\draw (9,2) -- (6,6);
		\draw (11,2) -- (8,6);
		\draw (12,2) -- (15,6);
		\draw (14,2) -- (17,6);
		\draw (15,2) -- (12,6);
		\draw (17,2) -- (14,6);
		\draw (18,2) -- (18,6);
		\draw (20,2) -- (20,6);
		\braidbox{0}{5}{0}{1}{\tiny $2^{b-s}$}
		\braidbox{6}{8}{0}{1}{\tiny $1^{s-n}$}
		\braidbox{9}{11}{0}{1}{\tiny $2^{s}$}
		\braidbox{12}{14}{0}{1}{\tiny $1^{n}$}
		\braidbox{15}{20}{0}{1}{\tiny $1^{a-s}$}
		\braidbox{0}{20}{1}{2}{\tiny $D(r,s,n)(u(s,n)\cdot \hat{f})$}
		\braidbox[0.4]{0}{2}{6}{7}{\tiny $2^{b-r}$}
		\braidbox{3}{5}{6}{7}{\tiny $1^{s-n}$}
		\braidbox{6}{11}{6}{7}{\tiny $2^{r}$}
		\braidbox{12}{17}{6}{7}{\tiny $1^{r-s+n}$}
		\braidbox[0.4]{18}{20}{6}{7}{\tiny $1^{a-r}$}
	\end{braid}.
\end{align*}
Define
\begin{align*}
	\varphi_{r,s}^{\hat{f}}: P^r := \bigoplus_{m=0}^r q^{\tts_{r,m}}R_\theta e_{r,m} &\to \bigoplus_{n=0}^s q^{\tts_{s,n}}R_\theta e_{s,n} =: P^s, \\
	(h_m e_{r,m})_{m=0}^r &\mapsto (h_{r-s+n} \hat{f}_{r,s}^n e_{s,n})_{n=0}^s.
\end{align*}
We think of $\varphi_{r,s}^{\hat{f}}$ as an element of $\Hom_{R_\theta}^{r-s}(P_\bullet^r, P_\bullet^s)$.
Recalling the differential (\ref{EDGDiff}) on $\Hom_{R_\theta}^\bullet(P_\bullet^r, P_\bullet^s)$, we have:

\begin{Lemma}\label{LA2ExtLift}
	Suppose $0 \leq s \leq r \leq \min\{ a, b\}$ and let $f\in \Lambda(r,s) = \CE_\theta(r,s)$.
	If $\hat{f}\in \hat{\Lambda}(r,s)$ is such that $\ttq_{r,s}(\hat{f}) = f$, then
	\begin{enumerate}
		\item\label{LA2ExtLifti} $\delta(\varphi_{r,s}^{\hat{f}}) = 0$, and
		\item\label{LA2ExtLiftii} the isomorphism (\ref{EExtIso}) sends the cohomology class of $\varphi_{r,s}^{\hat{f}}$ to $f$.
	\end{enumerate}
\end{Lemma}
\begin{proof}
	We prove (\ref{LA2ExtLifti}) by checking that the following diagram either commutes (if $r-s$ is even) or anticommutes (if $r-s$ is odd) whenever $0 \leq n < r-s$.
	$$\begin{tikzcd}[column sep=5em]
		q^{\tts_{r,n+1+r-s}}R_\theta e_{r,n+1+r-s} \arrow[r, "-\cdot d_{r,r-s+n}"] \arrow[d, "-\cdot\hat{f}_{n+1}"] & q^{\tts_{r,r-s+n}}R_\theta e_{r,r-s+n} \arrow[d, "-\cdot \hat{f}_n"] \\
		q^{\tts_{s,n+1}}R_\theta e_{s,n+1} \arrow[r, "-\cdot d_{s,n}"] & q^{\tts_{s,n}}R_\theta e_{s,n}.
	\end{tikzcd}$$
	We compute $\pm d_{r,r-s+n} \hat{f}_n$:
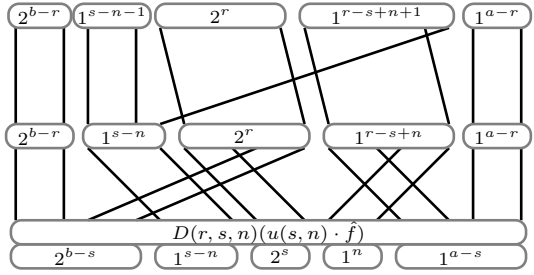
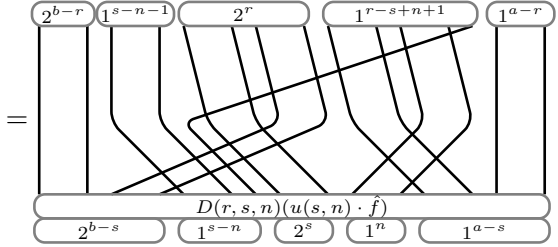
	\begin{align*}
		&
		\begin{braid}\tikzset{baseline=4.5em}
			\draw (0,2) -- (0,5);
			\draw (2,2) -- (2,5);
			\draw (3,2) -- (10,5);
			\draw (5,2) -- (12,5);
			\draw (6,2) -- (3,5);
			\draw (9,2) -- (6,5);
			\draw (10,2) -- (7,5);
			\draw (12,2) -- (9,5);
			\draw (13,2) -- (16,5);
			\draw (15,2) -- (18,5);
			\draw (16,2) -- (13,5);
			\draw (18,2) -- (15,5);
			\draw (19,2) -- (19,5);
			\draw (21,2) -- (21,5);
			\draw (0,6) -- (0,10);
			\draw (2,6) -- (2,10);
			\draw (3,6) -- (3,10);
			\draw (5,6) -- (5,10);
			\draw (6,6) -- (18,10);
			\draw (7,6) -- (6,10);
			\draw (12,6) -- (11,10);
			\draw (13,6) -- (12,10);
			\draw (18,6) -- (17,10);
			\draw (19,6) -- (19,10);
			\draw (21,6) -- (21,10);
			\braidbox{0}{5}{0}{1}{\tiny $2^{b-s}$}
			\braidbox{6}{9}{0}{1}{\tiny $1^{s-n}$}
			\braidbox{10}{12}{0}{1}{\tiny $2^{s}$}
			\braidbox{13}{15}{0}{1}{\tiny $1^{n}$}
			\braidbox{16}{21}{0}{1}{\tiny $1^{a-s}$}
			\braidbox{0}{21}{1}{2}{\tiny $D(r,s,n)(u(s,n)\cdot\hat{f})$}
			\braidbox[0.4]{0}{2}{5}{6}{\tiny $2^{b-r}$}
			\braidbox{3}{6}{5}{6}{\tiny $1^{s-n}$}
			\braidbox{7}{12}{5}{6}{\tiny $2^{r}$}
			\braidbox{13}{18}{5}{6}{\tiny $1^{r-s+n}$}
			\braidbox[0.4]{19}{21}{5}{6}{\tiny $1^{a-r}$}
			\braidbox[0.3]{0}{2}{10}{11}{\tiny $2^{b-r}$}
			\braidbox[0.6]{3}{5}{10}{11}{\tiny $1^{s-n-1}$}
			\braidbox{6}{11}{10}{11}{\tiny $2^{r}$}
			\braidbox{12}{18}{10}{11}{\tiny $1^{r-s+n+1}$}
			\braidbox[0.4]{19}{21}{10}{11}{\tiny $1^{a-r}$}
		\end{braid}
		& \\ 
		&
		\begin{braid}\tikzset{baseline=4em}
			\draw (0,2) -- (0,5) -- (0,9);
			\draw (2,2) -- (2,5) -- (2,9);
			\draw (3,2) -- (10,5) -- (9,9);
			\draw (5,2) -- (12,5) -- (11,9);
			\draw (6,2) -- (3,5) -- (3,9);
			\draw (8,2) -- (5,5) -- (5,9);
			\draw (9,2) -- (6,5) -- (18,9);
			\draw (10,2) -- (7,5) -- (6,9);
			\draw (12,2) -- (9,5) -- (8,9);
			\draw (13,2) -- (16,5) -- (15,9);
			\draw (15,2) -- (18,5) -- (17,9);
			\draw (16,2) -- (13,5) -- (12,9);
			\draw (18,2) -- (15,5) -- (14,9);
			\draw (19,2) -- (19,5) -- (19,9);
			\draw (21,2) -- (21,5) -- (21,9);
			\braidbox{0}{5}{0}{1}{\tiny $2^{b-s}$}
			\braidbox{6}{9}{0}{1}{\tiny $1^{s-n}$}
			\braidbox{10}{12}{0}{1}{\tiny $2^{s}$}
			\braidbox{13}{15}{0}{1}{\tiny $1^{n}$}
			\braidbox{16}{21}{0}{1}{\tiny $1^{a-s}$}
			\braidbox{0}{21}{1}{2}{\tiny $D(r,s,n)(u(s,n)\cdot\hat{f})$}
			\braidbox[0.3]{0}{2}{9}{10}{\tiny $2^{b-r}$}
			\braidbox[0.6]{3}{5}{9}{10}{\tiny $1^{s-n-1}$}
			\braidbox{6}{11}{9}{10}{\tiny $2^{r}$}
			\braidbox{12}{18}{9}{10}{\tiny $1^{r-s+n+1}$}
			\braidbox[0.4]{19}{21}{9}{10}{\tiny $1^{a-r}$}
			\draw (0,5) node[left]{\normalsize $=$};
		\end{braid}
		& \begin{tabular}{p{15em}} by Lemma~\ref{LDivPower} \end{tabular} \\
		&
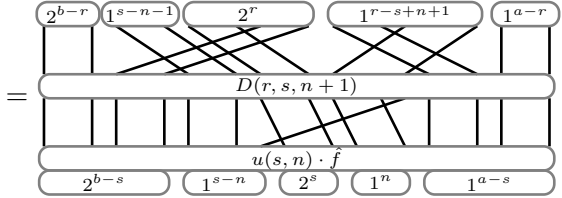
		\begin{braid}\tikzset{baseline=3em}
			\draw (0,2) -- (0,4);
			\draw (2,2) -- (2,4);
			\draw (3,2) -- (3,4);
			\draw (5,2) -- (5,4);
			\draw (6,2) -- (6,4);
			\draw (8,2) -- (8,4);
			\draw (9,2) -- (15,4);
			\draw (10,2) -- (9,4);
			\draw (12,2) -- (11,4);
			\draw (13,2) -- (12,4);
			\draw (15,2) -- (14,4);
			\draw (16,2) -- (16,4);
			\draw (18,2) -- (18,4);
			\draw (19,2) -- (19,4);
			\draw (21,2) -- (21,4);
			\draw (0,5) -- (0,7);
			\draw (2,5) -- (2,7);
			\draw (3,5) -- (9,7);
			\draw (5,5) -- (11,7);
			\draw (6,5) -- (3,7);
			\draw (8,5) -- (5,7);
			\draw (9,5) -- (6,7);
			\draw (11,5) -- (8,7);
			\draw (12,5) -- (15,7);
			\draw (15,5) -- (18,7);
			\draw (16,5) -- (12,7);
			\draw (18,5) -- (14,7);
			\draw (19,5) -- (19,7);
			\draw (21,5) -- (21,7);
			\braidbox{0}{5}{0}{1}{\tiny $2^{b-s}$}
			\braidbox{6}{9}{0}{1}{\tiny $1^{s-n}$}
			\braidbox{10}{12}{0}{1}{\tiny $2^{s}$}
			\braidbox{13}{15}{0}{1}{\tiny $1^{n}$}
			\braidbox{16}{21}{0}{1}{\tiny $1^{a-s}$}
			\braidbox{0}{21}{1}{2}{\tiny $u(s,n)\cdot \hat{f}$}
			\braidbox{0}{21}{4}{5}{\tiny $D(r,s,n+1)$}
			\braidbox[0.3]{0}{2}{7}{8}{\tiny $2^{b-r}$}
			\braidbox[0.6]{3}{5}{7}{8}{\tiny $1^{s-n-1}$}
			\braidbox{6}{11}{7}{8}{\tiny $2^{r}$}
			\braidbox{12}{18}{7}{8}{\tiny $1^{r-s+n+1}$}
			\braidbox[0.4]{19}{21}{7}{8}{\tiny $1^{a-r}$}
			\draw (0,4) node[left]{\normalsize $= \,$};
		\end{braid}
		& \begin{tabular}{p{15em}} by several applications of the relations (\ref{Rquad}) and (\ref{Rdots}), and by Lemma~\ref{TNHCenter} \end{tabular} \\
		&
		\begin{braid}\tikzset{baseline=3em}
			\draw (0,1) -- (0,3);
			\draw (5,1) -- (5,3);
			\draw (6,1) -- (6,3);
			\draw (8,1) -- (8,3);
			\draw (9,1) -- (15,3);
			\draw (10,1) -- (9,3);
			\draw (12,1) -- (11,3);
			\draw (13,1) -- (12,3);
			\draw (15,1) -- (14,3);
			\draw (16,1) -- (16,3);
			\draw (21,1) -- (21,3);
			\draw (0,4) -- (0,6);
			\draw (2,4) -- (2,6);
			\draw (3,4) -- (9,6);
			\draw (5,4) -- (11,6);
			\draw (6,4) -- (3,6);
			\draw (8,4) -- (5,6);
			\draw (9,4) -- (6,6);
			\draw (11,4) -- (8,6);
			\draw (12,4) -- (15,6);
			\draw (15,4) -- (18,6);
			\draw (16,4) -- (12,6);
			\draw (18,4) -- (14,6);
			\draw (19,4) -- (19,6);
			\draw (21,4) -- (21,6);
			\braidbox{0}{5}{0}{1}{\tiny $2^{b-s}$}
			\braidbox{6}{9}{0}{1}{\tiny $1^{s-n}$}
			\braidbox{10}{12}{0}{1}{\tiny $2^{s}$}
			\braidbox{13}{15}{0}{1}{\tiny $1^{n}$}
			\braidbox{16}{21}{0}{1}{\tiny $1^{a-s}$}
			\braidbox{0}{21}{3}{4}{\tiny $D(r,s,n+1)(u(s,n+1)\cdot\hat{f})$}
			\braidbox[0.3]{0}{2}{6}{7}{\tiny $2^{b-r}$}
			\braidbox[0.6]{3}{5}{6}{7}{\tiny $1^{s-n-1}$}
			\braidbox{6}{11}{6}{7}{\tiny $2^{r}$}
			\braidbox{12}{18}{6}{7}{\tiny $1^{r-s+n+1}$}
			\braidbox[0.4]{19}{21}{6}{7}{\tiny $1^{a-r}$}
			\draw (0,3.5) node[left]{\normalsize $= \,$};
		\end{braid}
		& \begin{tabular}{p{15em}} by Lemma~\ref{TNHCenter} and the relation (\ref{Rdots}) \end{tabular} \\
		&
		\begin{braid}\tikzset{baseline=3em}
			\draw (0,1) -- (0,3);
			\draw (5,1) -- (5,3);
			\draw (6,1) -- (6,3);
			\draw (8,1) -- (8,3);
			\draw (9,1) -- (15,3);
			\draw (10,1) -- (9,3);
			\draw (12,1) -- (11,3);
			\draw (13,1) -- (12,3);
			\draw (15,1) -- (14,3);
			\draw (16,1) -- (16,3);
			\draw (21,1) -- (21,3);
			\draw (0,5) -- (0,7);
			\draw (2,5) -- (2,7);
			\draw (3,5) -- (9,7);
			\draw (5,5) -- (11,7);
			\draw (6,5) -- (3,7);
			\draw (8,5) -- (5,7);
			\draw (9,5) -- (6,7);
			\draw (11,5) -- (8,7);
			\draw (12,5) -- (15,7);
			\draw (15,5) -- (18,7);
			\draw (16,5) -- (12,7);
			\draw (18,5) -- (14,7);
			\draw (19,5) -- (19,7);
			\draw (21,5) -- (21,7);
			\braidbox{0}{5}{0}{1}{\tiny $2^{b-s}$}
			\braidbox[0.6]{6}{9}{0}{1}{\tiny $1^{s-n}$}
			\braidbox{10}{12}{0}{1}{\tiny $2^{s}$}
			\braidbox{13}{15}{0}{1}{\tiny $1^{n}$}
			\braidbox{16}{21}{0}{1}{\tiny $1^{a-s}$}
			\braidbox{0}{5}{3}{4}{\tiny $2^{b-s}$}
			\braidbox[0.6]{6}{8}{3}{4}{\tiny $1^{s-n-1}$}
			\braidbox{9}{11}{3}{4}{\tiny $2^{s}$}
			\braidbox{12}{15}{3}{4}{\tiny $1^{n+1}$}
			\braidbox{16}{21}{3}{4}{\tiny $1^{a-s}$}
			\braidbox{0}{21}{4}{5}{\tiny $D(r,s,n+1)(u(s,n+1)\cdot\hat{f})$}
			\braidbox[0.3]{0}{2}{7}{8}{\tiny $2^{b-r}$}
			\braidbox[0.6]{3}{5}{7}{8}{\tiny $1^{s-n-1}$}
			\braidbox{6}{11}{7}{8}{\tiny $2^{r}$}
			\braidbox{12}{18}{7}{8}{\tiny $1^{r-s+n+1}$}
			\braidbox[0.4]{19}{21}{7}{8}{\tiny $1^{a-r}$}
			\draw (0,4) node[left]{\normalsize $= \,$};
		\end{braid}
		& \begin{tabular}{p{15em}} by Lemma~\ref{LDivPower}. \end{tabular}
	\end{align*}
	The last diagram represents $\pm \hat{f}_{n+1}d_{s,n}$.
	The signs are taken care of by
	$$(-1)^{\binom{r-s+1}{2}+(n+1)(r-s)} = (-1)^{r-s}(-1)^{\binom{r-s+1}{2}+n(r-s)}.$$

	To prove (\ref{LA2ExtLiftii}), recalling (\ref{EA2Aug}), we compute $(-1)^{\frac{(r-s)(r-s+1)}{2}}\eps_s(e_{r,r-s}\hat{f}_{r,s}^0)$:
	\begin{align*}
		&
		\begin{braid}\tikzset{baseline=3em}
			\draw (0,0) -- (0,3);
			\draw (5,0) -- (5,3);
			\draw (6.5,0) -- (6.5,1) -- (6,2) -- (6,3);
			\draw (10.5,0) -- (10.5,1) -- (10,2) -- (8,3);
			\draw (6.5,0) -- (6.5,1) -- (7,2) -- (9,3);
			\draw (10.5,0) -- (10.5,1) -- (11,2) -- (11,3);
			\draw (12,0) -- (12,3);
			\draw (17,0) -- (17,3);
			\draw (0,6) -- (0,8);
			\draw (2,6) -- (2,8);
			\draw (3,6) -- (9,8);
			\draw (5,6) -- (11,8);
			\draw (6,6) -- (3,8);
			\draw (8,6) -- (5,8);
			\draw (9,6) -- (6,8);
			\draw (11,6) -- (8,8);
			\draw (12,6) -- (12,8);
			\draw (14,6) -- (14,8);
			\draw (15,6) -- (15,8);
			\draw (17,6) -- (17,8);
			\braidbox{0}{5}{3}{4}{$w_0$}
			\braidbox{6}{8}{3}{4}{$w_0$}
			\braidbox{9}{11}{3}{4}{$w_0$}
			\braidbox{12}{17}{3}{4}{$w_0$}
			\braidbox{0}{5}{4}{5}{\tiny $2^{b-s}$}
			\braidbox{6}{8}{4}{5}{\tiny $1^s$}
			\braidbox{9}{11}{4}{5}{\tiny $2^s$}
			\braidbox{12}{17}{4}{5}{\tiny $1^{a-s}$}
			\braidbox{0}{17}{5}{6}{\tiny $\hat{f}$}
			\braidbox[0.3]{0}{2}{8}{9}{\tiny $2^{b-r}$}
			\braidbox{3}{5}{8}{9}{\tiny $1^s$}
			\braidbox{6}{11}{8}{9}{\tiny $2^r$}
			\braidbox{12}{14}{8}{9}{\tiny $1^{r-s}$}
			\braidbox[0.3]{15}{17}{8}{9}{\tiny $1^{a-r}$}
		\end{braid}
		& \\ 
		&
		\begin{braid}\tikzset{baseline=3em}
			\draw (0,0) -- (0,3);
			\draw (5,0) -- (5,3);
			\draw (6.5,0) -- (6.5,1) -- (6,2) -- (6,3);
			\draw (10.5,0) -- (10.5,1) -- (10,2) -- (8,3);
			\draw (6.5,0) -- (6.5,1) -- (7,2) -- (9,3);
			\draw (10.5,0) -- (10.5,1) -- (11,2) -- (11,3);
			\draw (12,0) -- (12,3);
			\draw (17,0) -- (17,3);
			\draw (0,5) -- (0,7);
			\draw (2,5) -- (2,7);
			\draw (3,5) -- (9,7);
			\draw (5,5) -- (11,7);
			\draw (6,5) -- (3,7);
			\draw (8,5) -- (5,7);
			\draw (9,5) -- (6,7);
			\draw (11,5) -- (8,7);
			\draw (12,5) -- (12,7);
			\draw (14,5) -- (14,7);
			\draw (15,5) -- (15,7);
			\draw (17,5) -- (17,7);
			\braidbox{0}{5}{3}{4}{$w_0$}
			\braidbox{6}{8}{3}{4}{$w_0$}
			\braidbox{9}{11}{3}{4}{$w_0$}
			\braidbox{12}{17}{3}{4}{$w_0$}
			\braidbox{0}{17}{4}{5}{\tiny $\hat{f}$}
			\draw (1,7.5) node{\tiny $2^{b-r}$};
			\draw (4,7.5) node{\tiny $1^s$};
			\draw (7,7.5) node{\tiny $2^s$};
			\draw (10,7.5) node{\tiny $2^{r-s}$};
			\draw (13,7.5) node{\tiny $1^{r-s}$};
			\draw (16,7.5) node{\tiny $1^{a-r}$};
			\draw (0,3.5) node[left]{\normalsize $= \,$};
		\end{braid}
		& \begin{tabular}{p{18em}} by Lemma~\ref{LDivPower} \end{tabular} \\
		&
		\begin{braid}\tikzset{baseline=3.5em}
			\draw (6.5,0) -- (6.5,1);
			\draw (10.5,0) -- (10.5,1);
			\draw (0,2) -- (0,8);
			\draw (2,2) -- (2,8);
			\draw (3,2) -- (3,6) -- (9,8);
			\draw (5,2) -- (5,6) -- (11,8);
			\draw (6.5,2) -- (6.5,3) -- (6,4) -- (6,5);
			\draw (10.5,2) -- (10.5,3) -- (10,4) -- (8,5);
			\draw (6.5,2) -- (6.5,3) -- (7,4) -- (9,5);
			\draw (10.5,2) -- (10.5,3) -- (11,4) -- (11,5);
			\draw (12,2) -- (12,8);
			\draw (14,2) -- (14,8);
			\draw (15,2) -- (15,8);
			\draw (17,2) -- (17,8);
			\draw (6,6) -- (3,8);
			\draw (8,6) -- (5,8);
			\draw (9,6) -- (6,8);
			\draw (11,6) -- (8,8);
			\braidbox{0}{5}{0}{1}{$w_0$}
			\braidbox{12}{17}{0}{1}{$w_0$}
			\braidbox{0}{17}{1}{2}{\tiny $\hat{\ttp}_{r,s}(\hat{f})$}
			\braidbox{6}{8}{5}{6}{$w_0$}
			\braidbox{9}{11}{5}{6}{$w_0$}
			\draw (1,8.5) node{\tiny $2^{b-r}$};
			\draw (4,8.5) node{\tiny $1^s$};
			\draw (7,8.5) node{\tiny $2^s$};
			\draw (10,8.5) node{\tiny $2^{r-s}$};
			\draw (13,8.5) node{\tiny $1^{r-s}$};
			\draw (16,8.5) node{\tiny $1^{a-r}$};
			\draw (0,4) node[left]{\normalsize $=$};
		\end{braid}
		& \begin{tabular}{p{18em}} by Lemma~\ref{TNHCenter}, the relation (\ref{Rdots}), and the defining relations of $\Delta(s)$. \end{tabular}
	\end{align*}
	The last diagram represents $\Xi_{r-s}(\hat{\ttp}_{r,s}(\hat{f}))$, where $\Xi_{r-s}$ is as in Lemma~\ref{LWeightBasis}(\ref{LWeightBasisii}).
	Thus, the image of $\varphi_{r,s}^{\hat{f}}$ under the map (\ref{EExtInduced}) is
	$$(e_{r,r-s}\mapsto \Xi_{r-s}(\hat{\ttp}_{r,s}(\hat{f})))\in \Hom_{R_\theta}(q^{\tts_{r,r-s}}R_\theta e_{r,r-s}, \Delta(s)).$$
	The proof is complete upon an application of Theorem~\ref{TA2Ext}.
\end{proof}

For $\hat{f}\in \hat{\Lambda}(r,s)$ and $\hat{g}\in \hat{\Lambda}(s,t)$, define
\begin{equation}\label{EDiamondHat}
	\blank \diamondhat \blank: \hat{\Lambda}(s,t)\otimes \hat{\Lambda}(r,s) \to \hat{\Lambda}(r,t), \ \hat{g} \otimes \hat{f} \mapsto \partial_{x(r,s,t)}\big(D(r,s,s-t)\hat{g}(w(s,t)\cdot \hat{f})\big).
\end{equation}

\begin{Lemma}\label{LA2LiftProd}
	Let $0 \leq t \leq s \leq r \leq \min\{ a, b\}$.
	If $\hat{f}\in \hat{\Lambda}(r,s)$ and $\hat{g}\in \hat{\Lambda}(s,t)$, then $\varphi_{s,t}^{\hat{g}}\varphi_{r,s}^{\hat{f}} = \varphi_{r,t}^{\hat{g}\diamondhat\hat{f}}$.
\end{Lemma}
\begin{proof}
	Let $y(s,t,n) := U_{b-s;t-n,s-t}\in \Sgp_d$ and note that
	\begin{equation}\label{EProdShuf}
		u(t,n)^{-1}y(s,t,n)u(s,s-t+n) = w(s,t)
	\end{equation}
	and
	\begin{equation}\label{EDProd}
		D(r,s,s-t+n)D(s,t,n) = D(r,t,n)D(r,s,s-t).
	\end{equation}
	We compute $\hat{f}_{r,s}^{s-t+n}\hat{g}_{s,t}^n$:
	\begin{align*}
		&
		\begin{braid}\tikzset{baseline=3em}
			\draw (0,2) -- (0,4);
			\draw (5,2) -- (5,4);
			\draw (6,2) -- (12,4);
			\draw (8,2) -- (14,4);
			\draw (9,2) -- (6,4);
			\draw (11,2) -- (8,4);
			\draw (12,2) -- (9,4);
			\draw (14,2) -- (11,4);
			\draw (15,2) -- (18,4);
			\draw (17,2) -- (20,4);
			\draw (18,2) -- (15,4);
			\draw (20,2) -- (17,4);
			\draw (21,2) -- (21,4);
			\draw (26,2) -- (26,4);
			\draw (0,6) -- (0,9);
			\draw (2,6) -- (2,9);
			\draw (3,6) -- (12,9);
			\draw (5,6) -- (14,9);
			\draw (6,6) -- (3,9);
			\draw (8,6) -- (5,9);
			\draw (9,6) -- (6,9);
			\draw (14,6) -- (11,9);
			\draw (15,6) -- (18,9);
			\draw (20,6) -- (23,9);
			\draw (21,6) -- (15,9);
			\draw (23,6) -- (17,9);
			\draw (24,6) -- (24,9);
			\draw (26,6) -- (26,9);
			\braidbox{0}{8}{0}{1}{\tiny $2^{b-t}$}
			\braidbox{9}{11}{0}{1}{\tiny $1^{t-n}$}
			\braidbox{12}{14}{0}{1}{\tiny $2^t$}
			\braidbox{15}{17}{0}{1}{\tiny $1^n$}
			\braidbox{18}{26}{0}{1}{\tiny $1^{a-t}$}
			\braidbox{0}{26}{1}{2}{\tiny $D(s,t,n)(u(t,n)\cdot \hat{g})$}
			\braidbox{0}{5}{4}{5}{\tiny $2^{b-s}$}
			\braidbox{6}{8}{4}{5}{\tiny $1^{t-n}$}
			\braidbox{9}{14}{4}{5}{\tiny $2^s$}
			\braidbox{15}{20}{4}{5}{\tiny $1^{s-t+n}$}
			\braidbox{21}{26}{4}{5}{\tiny $1^{a-s}$}
			\braidbox{0}{26}{5}{6}{\tiny $D(r,s,s-t+n)(u(s,s-t+n)\cdot \hat{f})$}
			\braidbox{0}{2}{9}{10}{\tiny $2^{b-r}$}
			\braidbox{3}{5}{9}{10}{\tiny $1^{t-n}$}
			\braidbox{6}{14}{9}{10}{\tiny $2^r$}
			\braidbox{15}{23}{9}{10}{\tiny $1^{r-t+n}$}
			\braidbox{24}{26}{9}{10}{\tiny $1^{a-r}$}
		\end{braid}
		& \\ 
		&
		\begin{braid}\tikzset{baseline=3em}
			\draw (0,2) -- (0,7);
			\draw (2,2) -- (2,7);
			\draw (3,2) -- (3,4) -- (12,7);
			\draw (5,2) -- (5,4) -- (14,7);
			\draw (6,2) -- (12,4) -- (9,7);
			\draw (8,2) -- (14,4) -- (11,7);
			\draw (9,2) -- (6,4) -- (3,7);
			\draw (11,2) -- (8,4) -- (5,7);
			\draw (12,2) -- (9,4) -- (6,7);
			\draw (14,2) -- (11,4) -- (8,7);
			\draw (15,2) -- (18,4) -- (21,7);
			\draw (17,2) -- (20,4) -- (23,7);
			\draw (18,2) -- (15,4) -- (18,7);
			\draw (20,2) -- (17,4) -- (20,7);
			\draw (21,2) -- (21,4) -- (15,7);
			\draw (23,2) -- (23,4) -- (17,7);
			\draw (24,2) -- (24,7);
			\draw (26,2) -- (26,7);
			\braidbox{0}{8}{0}{1}{\tiny $2^{b-t}$}
			\braidbox{9}{11}{0}{1}{\tiny $1^{t-n}$}
			\braidbox{12}{14}{0}{1}{\tiny $2^t$}
			\braidbox{15}{17}{0}{1}{\tiny $1^n$}
			\braidbox{18}{26}{0}{1}{\tiny $1^{a-t}$}
			\braidbox{0}{26}{1}{2}{\tiny $D(r,s,s-t+n)D(s,t,n)(u(t,n)\cdot \hat{g})(y(s,t,n)u(s,s-t+n)\cdot \hat{a})$}
			\braidbox{0}{2}{7}{8}{\tiny $2^{b-r}$}
			\braidbox{3}{5}{7}{8}{\tiny $1^{t-n}$}
			\braidbox{6}{14}{7}{8}{\tiny $2^r$}
			\braidbox{15}{23}{7}{8}{\tiny $1^{r-t+n}$}
			\braidbox{24}{26}{7}{8}{\tiny $1^{a-r}$}
			\draw (0,4) node[left]{\normalsize $=$};
		\end{braid}
		& \begin{tabular}{p{11.4em}} by (\ref{Rdots}), Lemma~\ref{LDivPower}, and Theorem~\ref{TNHCenter} \end{tabular} \\
		&
		\begin{braid}\tikzset{baseline=3em}
			\draw (0,2) -- (0,7);
			\draw (2,2) -- (2,7);
			\draw (3,2) -- (6,4) -- (12,7);
			\draw (5,2) -- (8,4) -- (14,7);
			\draw (6,2) -- (3,4) -- (9,7);
			\draw (8,2) -- (5,4) -- (11,7);
			\draw (9,2) -- (9,4) -- (3,7);
			\draw (11,2) -- (11,4) -- (5,7);
			\draw (12,2) -- (12,4) -- (6,7);
			\draw (14,2) -- (14,4) -- (8,7);
			\draw (15,2) -- (15,4) -- (21,7);
			\draw (17,2) -- (17,4) -- (23,7);
			\draw (18,2) -- (21,4) -- (18,7);
			\draw (20,2) -- (23,4) -- (20,7);
			\draw (21,2) -- (18,4) -- (15,7);
			\draw (23,2) -- (20,4) -- (17,7);
			\draw (24,2) -- (24,7);
			\draw (26,2) -- (26,7);
			\braidbox{0}{8}{0}{1}{\tiny $2^{b-t}$}
			\braidbox{9}{11}{0}{1}{\tiny $1^{t-n}$}
			\braidbox{12}{14}{0}{1}{\tiny $2^t$}
			\braidbox{15}{17}{0}{1}{\tiny $1^n$}
			\braidbox{18}{26}{0}{1}{\tiny $1^{a-t}$}
			\braidbox{0}{26}{1}{2}{\tiny $D(r,t,n)D(r,s,s-t)u(t,n)\cdot(\hat{g}(w(s,t)\cdot \hat{f}))$}
			\braidbox{0}{2}{7}{8}{\tiny $2^{b-r}$}
			\braidbox{3}{5}{7}{8}{\tiny $1^{t-n}$}
			\braidbox{6}{14}{7}{8}{\tiny $2^r$}
			\braidbox{15}{23}{7}{8}{\tiny $1^{r-t+n}$}
			\braidbox{24}{26}{7}{8}{\tiny $1^{a-r}$}
			\draw (0,4) node[left]{\normalsize $=$};
		\end{braid}
		& \begin{tabular}{p{11.4em}} by (\ref{EProdShuf}), (\ref{EDProd}), and the relation~(\ref{Rbraid}) \end{tabular} \\
		&
		\begin{braid}\tikzset{baseline=2.5em}
			\draw (0,2) -- (0,5);
			\draw (2,2) -- (2,5);
			\draw (3,2) -- (9,5);
			\draw (8,2) -- (14,5);
			\draw (9,2) -- (3,5);
			\draw (11,2) -- (5,5);
			\draw (12,2) -- (6,5);
			\draw (14,2) -- (8,5);
			\draw (15,2) -- (21,5);
			\draw (17,2) -- (23,5);
			\draw (18,2) -- (15,5);
			\draw (23,2) -- (20,5);
			\draw (24,2) -- (24,5);
			\draw (26,2) -- (26,5);
			\braidbox{0}{8}{0}{1}{\tiny $2^{b-t}$}
			\braidbox{9}{11}{0}{1}{\tiny $1^{t-n}$}
			\braidbox{12}{14}{0}{1}{\tiny $2^t$}
			\braidbox{15}{17}{0}{1}{\tiny $1^n$}
			\braidbox{18}{26}{0}{1}{\tiny $1^{a-t}$}
			\braidbox{0}{26}{1}{2}{\tiny $D(r,t,n)(u(t,n)\cdot\partial_{x(s,t)}(D(r,s,s-t)\hat{g}(w(s,t)\cdot \hat{f})))$}
			\braidbox{0}{2}{5}{6}{\tiny $2^{b-r}$}
			\braidbox{3}{5}{5}{6}{\tiny $1^{t-n}$}
			\braidbox{6}{14}{5}{6}{\tiny $2^r$}
			\braidbox{15}{23}{5}{6}{\tiny $1^{r-t+n}$}
			\braidbox{24}{26}{5}{6}{\tiny $1^{a-r}$}
			\draw (0,3) node[left]{\normalsize $=$};
		\end{braid}
		& \begin{tabular}{p{11.4em}} by Lemma~\ref{LPartialPoly}. \end{tabular}
	\end{align*}
	The last diagram represents $(\hat{g}\diamondhat \hat{f})^n$.
\end{proof}

\begin{Lemma}\label{LSchurDiamondHat}
	Let $0 \leq t \leq s \leq r \leq \min\{ a, b\}$, $\lambda \in \Par(r-s,s)$, and $\mu \in \Par(s-t,t)$, and set $\hat{f} := 1\otimes s_\lambda \otimes 1 \otimes 1 \otimes 1 \otimes 1 \in \hat{\Lambda}(r,s)$ and $\hat{g} := 1\otimes s_\mu \otimes 1 \otimes 1 \otimes 1 \otimes 1 \in \hat{\Lambda}(s,t)$.
	Then
	$$\hat{g}\diamondhat\hat{f} = 1\otimes (s_\lambda \star s_\mu) \otimes 1 \otimes 1 \otimes 1 \otimes 1 \in \hat{\Lambda}(r,t).$$
\end{Lemma}
\begin{proof}
	We use the notation $\lambda^c$ for the conjugate of a partition $\lambda$.
	For a partition $\kappa = (\kappa_1, \ldots, \kappa_{s-t})\in \Par(s-t,r-s)$, let $\hat{\kappa} := (r-s-\kappa_{s-t}, \ldots, r-s-\kappa_1)^c \in \Par(r-s,s-t)$.
	We have $\hat{f} = w(s,t)\hat{f}$ and
	$$D(r,s,s-t) = \sum_{\kappa\in \Par(s-t, r-s)}(-1)^{|\hat{\kappa}|} 1\otimes s_{\hat{\kappa}}\otimes 1 \otimes 1 \otimes 1 \otimes s_{\kappa} \otimes 1 \otimes 1 \in \Lambda(r,s,t)$$
	by~\cite[I.4 Example 5]{Mac}, so that
	\begin{align*}
		\hat{g}\diamondhat\hat{f} &= \sum_{\kappa\in \Par(s-t, r-s)}(-1)^{|\hat{\kappa}|}\partial_{x(r,s,t)}(1\otimes (s_{\hat{\kappa}}s_\lambda)\otimes s_\mu \otimes 1 \otimes 1 \otimes s_\kappa \otimes 1 \otimes 1) \\
		&= \sum_{\kappa\in \Par(s-t, r-s)}(-1)^{|\hat{\kappa}|}1\otimes ((s_{\hat{\kappa}}s_\lambda)\star s_\mu) \otimes 1 \otimes 1 \otimes (s_\kappa \star 1) \otimes 1 \\
		&= 1\otimes (s_\lambda\star s_\mu) \otimes 1 \otimes 1 \otimes (s_{(r-s)^{(s-t)}} \star 1) \otimes 1 \\
		&= 1\otimes (s_\lambda\star s_\mu) \otimes 1 \otimes 1 \otimes 1 \otimes 1. \qedhere
	\end{align*}
\end{proof}

We consider $\Lambda(r,s)$ as a right $\Lambda(r,r)$ module via the natural algebra embedding
$$\Lambda(r,r)=\Lambda_{b-r,r,a-r}\into \Lambda_{b-r,r-s,s,a-r}=q^{-\omega_{r-s}(r,s)}\Lambda(r,s).$$
If $f\in \Lambda_{r-s}$, we write
$$f^{r,s} := 1_{\Lambda_{b-r}} \otimes f \otimes 1_{\Lambda_r} \otimes 1_{\Lambda_{a-r}}\in \Lambda(r,s).$$
By Proposition~\ref{PPartSymBasis}, $\Lambda(r,s)$ is a free right $\Lambda(r,r)$-module with basis $$\{ s_\lambda^{r,s} \mid \lambda \in \Par(r-s,s)\}.$$
We also make $\Lambda(r, s)$ into a left $\Lambda(s, s)$-module via the composition of algebra homomorphisms:
\begin{equation}\label{ELeftModule}
	\Lambda(s,s)=\Lambda_{b-s,s,a-s} \into \Lambda_{b-r,r-s,s,r-s,a-r}=\Lambda^{r,s} \stackrel{\ttp_{r,s}}{\longrightarrow} q^{-\omega_{r-s}(r,s)} \Lambda(r,s),
\end{equation}
the first map being the natural embedding.
(This is similar to the definition of the \emph{right} $\Lambda_{\sigma(s)}$-module structure on $\Lambda(r,s)$ used in the previous subsection; in fact $\Lambda(s,s)$ can be  identified with $\End_{R_\theta}(\Delta(s))=\Lambda_{\sigma(s)}^\op$, see the proof of Theorem~\ref{TA2Alg}.)

If $0 \leq t \leq s \leq r \leq \min\{ a, b\}$, then the tensor product $\Lambda(s,t)\otimes_{\Lambda(s,s)}\Lambda(r,s)$ is now a free right $\Lambda(r,r)$-module with basis
\begin{equation}\label{ERightModuleBasis}
	\{ s_\mu^{s,t} \otimes s_\lambda^{r,s}\mid \mu\in \Par(s-t,t), \lambda\in \Par(r-s,s)\}
\end{equation}
and we define a map of right $\Lambda(r,r)$-modules
$$\Theta: \Lambda(s,t)\otimes_{\Lambda(s,s)}\Lambda(r,s) \to \Lambda(r,t), \ s_\mu^{s,t}\otimes s_\lambda^{r,s} \mapsto (s_\mu\star s_\lambda)^{r,t}.$$
Let
$$\blank\diamond \blank : \Lambda(s,t)\otimes_\k \Lambda(r,s) \to \Lambda(r,t), \ g\otimes f \mapsto \Theta(g\otimes f).$$

\begin{Theorem}\label{TA2Alg}
	Let $0 \leq t \leq s \leq r \leq \min\{ a, b\}$.
	The composition in the category $\CE_\theta$ is given by
	$$\begin{tikzcd}[row sep = -0.5em]
		\CE_\theta(s,t) \otimes \CE_\theta(r,s) \arrow[r] & \CE_\theta(r,t) \\
		\rotatebox{90}{$\cong$} & \rotatebox{90}{$\cong$} \\
		\Lambda(s,t) \otimes \Lambda(r,s) \arrow[r] & \Lambda(r,t) \\
		g \otimes f \arrow[r, mapsto] & g\diamond f.
	\end{tikzcd}$$
\end{Theorem}
\begin{proof}
	According to Lemmas~\ref{LA2ExtLift} and~\ref{LA2LiftProd}, the composition of $g\in \Lambda(s,t)$ with $f\in \Lambda(r,s)$ is given by $\ttq_{r,t}(\hat{g}\diamondhat \hat{f})$, where $\hat{g} \in \hat{\Lambda}(s,t)$ and $\hat{f} \in \hat{\Lambda}(r,s)$ are such that $\ttq_{s,t}(\hat{g}) = g$ and $\ttq_{r,s}(\hat{f}) = f$.

	First suppose $r = s$.
	Since $D(s,s,s-t) = 1$ and $x(s,s,t) = 1$, (\ref{EDiamondHat}) shows that $\ttq_{s,t}(\hat{g}\diamondhat\hat{f}) = gf \in \Lambda(s,t)$.
	Thus, the composition map $\Lambda(s,t)\otimes \Lambda(s,s) \to \Lambda(s,s)$ coincides with the right $\Lambda(s,s)$-module structure of $\Lambda(s,t)$.

	Now suppose $s = t$.
	Since $D(r,s,0) = 1$ and $w(s,s) = x(r,s,s) = 1$, (\ref{EDiamondHat}) shows that $\hat{g}\diamondhat\hat{f} = \hat{g}\hat{f}$, so that $\ttq_{r,s}(\hat{g}\diamondhat\hat{f}) = \varphi(g)f \in \Lambda(r,s)$ where $\varphi$ is the composition in (\ref{ELeftModule}).
	Thus, the composition map $\Lambda(s,s)\otimes \Lambda(r,s) \to \Lambda(r,s)$ coincides with the left $\Lambda(s,s)$-module structure of $\Lambda(r,s)$.

	Associativity in $\CE_\theta$ implies that the composition map $\Lambda(s,t) \otimes \Lambda(r,s) \to \Lambda(r,t)$ is $\Lambda(s,s)$-balanced and $\Lambda(r,r)$-equivariant, and is therefore completely determined by the image of the basis elements in (\ref{ERightModuleBasis}).
	Lemma~\ref{LSchurDiamondHat} now completes the proof.
\end{proof}

\section{Non-formality of the $A_\infty$ structure}\label{SExample}

We provide an example to show that the $A_\infty$-category structure on $\CE_\theta$ is, in general, non-formal.
First, we recall some basic definitions and an important theorem, see~\cite{Ke1, Ke2}.

A ($\k$-linear) \emph{$A_\infty$-category} $\CA$ consists of a class of objects $\ob(\CA)$, a $\Z$-graded $\k$-module $\CA^\bullet(\rho, \sigma)$ for every pair of objects $\rho, \sigma \in \ob(\CA)$, and for each $n \in \Z_{>0}$ and objects $\pi_0, \ldots, \pi_n \in \ob(\CA)$, a degree $2-n$ $\k$-linear map
$$m_n^{\pi_0, \ldots, \pi_n}: \CA^\bullet(\pi_{n-1},\pi_n) \otimes \CA^\bullet(\pi_{n-2}, \pi_{n-1}) \otimes \cdots \otimes \CA^\bullet(\pi_0, \pi_1) \to \CA^\bullet(\pi_0, \pi_n),$$
such that for each $n \in \Z_{>0}$, we have
\begin{equation}\label{EAn} \tag{$*_n$}
	\sum_{r+s+t=n}(-1)^{r+st}m_{r+1+t}(1^{\otimes r} \otimes m_s \otimes 1^{\otimes t}) = 0,
\end{equation}
where $m_n := \bigoplus_{\pi_0, \ldots, \pi_n \in \ob(\CA)}m_n^{\pi_0, \ldots, \pi_n}$.
An \emph{$A_\infty$-functor} $F$ from $\CB = (\CB, M_1, M_2, \ldots)$ to $\CA = (\CA, m_1, m_2, \ldots)$ consists of a function $F: \ob(\CB) \to \ob(\CA)$, and for each $n \in \Z_{\geq 0}$ and objects $\pi_0, \ldots, \pi_n \in \ob(\CA)$, a degree $1-n$ map
$$F_n^{\pi_0, \ldots, \pi_n}: \CB^\bullet(\pi_{n-1},\pi_n) \otimes \CB^\bullet(\pi_{n-2}, \pi_{n-1}) \otimes \cdots \otimes \CB^\bullet(\pi_0, \pi_1) \to \CA^\bullet(F(\pi_0), F(\pi_n)),$$
such that for each $n \in \Z_{>0}$, we have
\begin{equation}\label{EFn} \tag{$**_n$}
	\sum_{r+s+t = n}(-1)^{r+st}F_{r+1+t}(1^{\otimes r} \otimes M_s \otimes 1^{\otimes t}) = \sum_{i_1+\cdots+i_k = n}(-1)^\tts m_k(F_{i_1}\otimes \cdots \otimes F_{i_k}),
\end{equation}
where $\tts := \sum_{j=1}^{k-1}(k-j)(i_j-1)$ and $F_n := \bigoplus_{\pi_0, \ldots, \pi_n \in \ob(\CA)}F_n^{\pi_0, \ldots, \pi_n}$.
The \emph{composition} $GF$ of two $A_\infty$-functors $F: \CB \to \CA$ and $G: \CC \to \CB$ is given on objects by $(FG)(\pi) = F(G(\pi))$ for $\pi \in \ob(\CC)$ with
$$(FG)_n := \sum_{i_1+\cdots+i_k = n}(-1)^\tts F_k(G_{i_1}\otimes \cdots G_{i_k}),$$
where $\tts$ is as above.
The \emph{identity} functor $1_{\CA}$ on $\CA$ is the $A_\infty$-functor $F:\CA\to \CA$ with $F(\pi) = \pi$ for each $\pi \in \ob(\CA)$, $F_1^{\rho,\sigma} = \id_{\CA^\bullet(\rho, \sigma)}$ for each $\rho,\sigma \in \ob(\CA)$, and $F_n = 0$ for $n > 1$.
An \emph{isomorphism} of $A_\infty$-categories is an $A_\infty$-functor $F: \CB \to \CA$ such that there exists an $A_\infty$-functor $G: \CA \to \CB$ such that $GF = 1_{\CB}$ and $FG = 1_{\CA}$.

When we speak of the \emph{homology} $H\CA$ of an $A_\infty$-category $\CA = (\CA, m_1, m_2, \ldots)$, we mean the homology with respect to $m_1$.
It is a graded category with the same object class as $\CA$.
In fact, we have:

\begin{Theorem}~\cite[Theorem 1]{Kad}\label{TKad}
	Let $(\CA, m_1, m_2, \ldots)$ be an $A_\infty$ category and $H\CA$ its homology.
	If each morphism space in $H\CA$ is a free graded $\k$-module, then $H\CA$ carries the structure of an $A_\infty$-category $H\CA = (H\CA, M_1, M_2, \ldots)$ such that
	\begin{enumerate}
		\item\label{TKadi} $M_1 = 0$ and $M_2 = [m_2]$,
		\item\label{TKadii} there exists an $A_\infty$-functor $F: H\CA \to \CA$ such that $F(\pi) = \pi$ for each $\pi \in \ob(\CA)$ and $[F_1] = \id_{H\CA}$.
	\end{enumerate}
	Moreover, the $A_\infty$-category structure on $H\CA$ satisfying (\ref{TKadi}) and (\ref{TKadii}) is unique up to (non-unique) isomorphism of $A_\infty$-categories.
\end{Theorem}

Such an $A_\infty$-category structure on $H\CA$ is called a \emph{minimal model} of $\CA$.
An $A_\infty$-category $\CA$ is called \emph{formal} if its minimal model can be chosen so that $M_n = 0$ for $n\neq 2$.
A graded category $\CB$ is \emph{intrinsically formal} if every $A_\infty$-category $\CA = (\CA, m_1, m_2, \ldots)$ whose homology is isomorphic to $\CB$ as a graded category is formal.
For example, the graded category $\CE_\theta$ in the situation of either Theorems~\ref{TA} or~\ref{TB} is intrinsically formal because by homological degree consideration, there is no way to impose an $A_\infty$-category structure on $\CE_\theta$ with $M_n \neq 0$ for any $n \neq 2$.

Kadeishvili's original proof is constructive and yields an inductive algorithm for producing a minimal model in the special case where $\CA$ is a differential-graded algebra, i.e., $m_n = 0$ for $n > 2$:

\begin{Algorithm}~\cite[Proof of Theorem 1]{Kad}\label{AKad}
	Let $\CA = (\CA, m_1, m_2, \ldots)$ be an $A_\infty$-category with $m_k = 0$ for $k > 2$ and $H\CA$ its homology.
	The following algorithm produces an $A_\infty$-category structure $(H\CA, M_1, M_2, \ldots)$ and an $A_\infty$-functor $F:H\CA \to \CA$ which satisfies the conditions in Theorem~\ref{TKad}.
\end{Algorithm}
\begin{enumerate}
	\item[\underline{Step 1:}]
	Let $M_1 = 0$ and take $F_1^{\mu,\nu}: \CE_\theta(\mu,\nu) \to \CH_\theta(\mu,\nu)$ to be a cycle-choosing homomorphism of $\k$-modules.
	Set $n := 2$.

	\item[\underline{Step 2:}]
	Since $m_k = 0$ for $k > 2$, we may rewrite (\ref{EFn}) as
	\begin{equation}\label{EFunctorSimplified}
		m_1 F_n = F_1 M_n - U_n,
	\end{equation}
	where
	$$U_n := m_2 \sum_{i=1}^{n-1}(-1)^{i-1}(F_i\otimes F_{n-i}) - \sum_{s=2}^{n-1}\sum_{t=0}^{n-2}(-1)^{n-s-t+st}F_{n-s+1}(1^{\otimes n-s-t}\otimes M_s \otimes 1^{\otimes t}).$$
	We will also use the restriction to $H\CA(\pi_{n-1},\pi_n)\otimes \cdots \otimes H\CA(\pi_0,\pi_1)$:
	$$U_n^{\pi_0, \ldots, \pi_n} : H\CA(\pi_{n-1},\pi_n)\otimes \cdots \otimes H\CA(\pi_0,\pi_1) \to \CA(\pi_0,\pi_n).$$
	One can check that $m_1 U_n = 0$.
	Thus, since $M_k, F_k$ have been defined for $k < n$, we take $M_n$ to be the (well-defined) homology class $[U_n]$ of $U_n$.

	\item[\underline{Step 3:}]
	Note that $[F_1 M_n - U_n] = [F_1 M_n] - [U_n] = M_n - M_n = 0$, so $F_1 M_n - U_n$ is a boundary, and choose $F_n$ such that $m_1 F_n = F_1 M_n - U_n$.
	Increment $n$ and return to Step 2.
\end{enumerate}

Let $\theta \in Q_+$.
For each $\pi \in \KP(\theta)$, fix a projective resolution $P_\bullet^\pi$ of $\Delta(\pi)$.
Consider the differential-graded category $\CH_\theta$ whose objects are the Kostant partitions of $\theta$ with morphism spaces $\CH_\theta^\bullet(\rho, \sigma) := \Hom_{R_\theta}^\bullet(P_\bullet^\rho, P_\bullet^\sigma)$.
We denote by
\begin{align*}
	m_1^{\rho,\sigma}&: \CH_\theta^\bullet(\rho,\sigma) \to \CH_\theta^\bullet(\rho,\sigma) \\
	m_2^{\rho,\sigma,\tau}&: \CH_\theta^\bullet(\sigma,\tau) \otimes \CH_\theta^\bullet(\rho,\sigma) \to \CH_\theta^\bullet(\rho,\tau)
\end{align*}
the differential and composition in $\CH_\theta$, respectively.
Note that $m_1^{\rho,\sigma}$ is precisely $\delta$ from (\ref{EDGDiff}).
Being a differential-graded category, $\CH_\theta$ is also an $A_\infty$-category (see~\cite{Ke2}), so its homology $\CE_\theta$ carries a structure of an $A_\infty$-category according to Theorem~\ref{TKad}.

For the rest of this section, we let $\theta := \alpha_1 + 2\alpha_2 + \alpha_3 \in Q_+$ and set
\begin{align*}
	\pi &:= (\alpha_2, \alpha_1+\alpha_2+\alpha_3) & \rho &:= (\alpha_2+\alpha_3, \alpha_1+\alpha_2) \\
	\sigma &:= (\alpha_3, \alpha_2, \alpha_1+\alpha_2) & \tau &:= (\alpha_3, (\alpha_2)^2, \alpha_1).
\end{align*}
Note that there is one other Kostant partition, $(\alpha_2+\alpha_3, \alpha_2, \alpha_1)$, which will not play a role in our construction.
Recall the standard generators (\ref{EDeltaGen}) of the standard modules and the idempotents (\ref{ENHIdem}) in the nil-Hecke algebra.
We define
\begin{align*}
	\hat{v}_\tau &:= v_{\alpha_3}\circ v_{\alpha_2}\circ v_{\alpha_2}\circ v_{\alpha_1} \in \hat{\Delta}(\tau), \\
	e_\tau &:= e_1\otimes e_2 \otimes e_1 \in \NH_1 \otimes \NH_2 \otimes \NH_1,
\end{align*}
so that $v_\tau = \hat{v}_\tau e_\tau$.
We list the resolutions of the corresponding standard modules from~\cite[Theorem A]{BKM} and Lemma~\ref{LCircRes} below:
$$\begin{tikzcd}[column sep=2em, row sep=1em]
	0 \arrow[r] & [-1em] q^2 R_\theta 1_{2321} \arrow[r,"d_1^\pi"] & q R_\theta 1_{2231} \oplus q R_\theta 1_{2132} \arrow[r,"d_0^\pi"] & R_\theta 1_{2123} \arrow[r,"\epsilon_\pi"] & \Delta(\pi) \arrow[r] & [-1em] 0, \\
	0 \arrow[r] & q^2 R_\theta 1_{3221} \arrow[r,"d_1^\rho"] & q R_\theta 1_{2321} \oplus q R_\theta 1_{3212} \arrow[r,"d_0^\rho"] & R_\theta 1_{2312} \arrow[r,"\epsilon_\rho"] & \Delta(\rho) \arrow[r] & 0, \\
	& 0 \arrow[r] & q R_\theta 1_{3221} \arrow[r," d_0^\sigma"] & R_\theta 1_{3212} \arrow[r,"\epsilon_\sigma"] & \Delta(\sigma) \arrow[r] & 0, \\
	& & 0 \arrow[r] & q^{-1} R_\theta 1_{32^{(2)}1} \arrow[r,"\epsilon_\tau"] & \Delta(\tau) \arrow[r] & 0,
\end{tikzcd}$$
where a matrix label stands for right multiplication with that matrix, and
\begin{align*}
	d_1^\pi &:=
	\left[\begin{matrix}
		-\psi_2 1_{2231} & \psi_3\psi_2 1_{2132}
	\end{matrix}\right], &
	d_0^\pi &:=
	\left[\begin{matrix}
		\psi_3\psi_2 1_{2123} \\
		\psi_3 1_{2123}
	\end{matrix}\right], &
	\epsilon_\pi &:=
	\left[\begin{matrix}
		v_\pi
	\end{matrix}\right], \\
	d_1^\rho &:=
	\left[\begin{matrix}
		-\psi_1 1_{2321} & \psi_3 1_{3212}
	\end{matrix}\right], &
	d_0^\rho &:=
	\left[\begin{matrix}
		\psi_3 1_{2312} \\
		\psi_1 1_{2312}
	\end{matrix}\right], &
	\epsilon_\rho &:=
	\left[\begin{matrix}
		v_\rho
	\end{matrix}\right], \\
	&&
	d_0^\sigma &:=
	\left[\begin{matrix}
		\psi_3 1_{3212}
	\end{matrix}\right], &
	\epsilon_\sigma &:=
	\left[\begin{matrix}
		v_\sigma
	\end{matrix}\right], \\
	&&&&
	\epsilon_\tau &:=
	\left[\begin{matrix}
		\psi_2 v_\tau
	\end{matrix}\right] =
	\left[\begin{matrix}
		\hat{v}_\tau \tau_2
	\end{matrix}\right].
\end{align*}

One can easily check, using (\ref{EShuffleMod}), that the complexes  $T_\bullet^\pi(\rho)$, $T_\bullet^\rho(\sigma)$, $T_\bullet^\sigma(\tau)$, and $T_\bullet^\pi(\tau)$ are, respectively, the top complexes in the four diagrams below.
Moreover, the diagrams define isomorphisms of complexes.
$$\begin{tikzcd}[column sep=2em, ampersand replacement=\&]
	\& [-1em] 0 \& [-0.5em] q^{-1}1_{2132}\Delta(\rho) \arrow[l] \& [3em] 1_{2123}\Delta(\rho) \arrow[l,"{\left[(d_0^\pi)_{2,1}\right]}" above] \& [-1em] 0 \arrow[l] \\
	\& 0 \& q^{-1}\CX_2 \arrow[l] \arrow[u,"f \mapsto \psi_2 v_\rho f" right,"\rotatebox{90}{$\sim$}" left] \& q\CX_2 \arrow[l,"{\left[x_2-x_1\right]}" above] \arrow[u,"f \mapsto \psi_3 \psi_2 v_\rho f" right,"\rotatebox{90}{$\sim$}" left] \& 0 \arrow[l], \\
	\& 0 \& q^{-1}1_{3212}\Delta(\sigma) \arrow[l] \& 1_{2312}\Delta(\sigma) \arrow[l,"{\left[(d_0^\rho)_{2,1}\right]}" above] \& 0 \arrow[l] \\
	\& 0 \& q^{-1}\CX_3 \arrow[l] \arrow[u,"f \mapsto v_\sigma f" right,"\rotatebox{90}{$\sim$}" left] \& q\CX_3 \arrow[l,"{\left[x_2-x_1\right]}" above] \arrow[u,"f \mapsto \psi_1 v_\sigma f" right,"\rotatebox{90}{$\sim$}" left] \& 0 \arrow[l], \\
	\& 0 \& 1_{3221}\Delta(\tau) \arrow[l] \& q 1_{3212}\Delta(\tau) \arrow[l,"{\left[(d_0^\sigma)_{1,1}\right]}" above] \& 0 \arrow[l] \\
	\& 0 \& q^{-2}\CX_4 \arrow[l] \arrow[u,"f \mapsto \hat{v}_\tau f \tau_2" right,"\rotatebox{90}{$\sim$}" left] \& \CX_4 \arrow[l,"{\left[x_4-x_3\right]}" above] \arrow[u,"f \mapsto \psi_3 \hat{v}_\tau f \tau_2" right,"\rotatebox{90}{$\sim$}" left] \& 0 \arrow[l], \\
	0 \& q^{-1}1_{2321}\Delta(\tau) \arrow[l] \& 1_{2231}\Delta(\tau) \oplus 1_{2132}\Delta(\tau) \arrow[l,"d_1^\pi" above] \& q 1_{2123}\Delta(\tau) \arrow[l,"d_0^\pi" above] \& 0 \arrow[l] \\
	0 \& q^{-2}\CX_4 \arrow[l] \arrow[u,"f\mapsto \psi_1\hat{v}_\tau f \tau_2" right,"\rotatebox{90}{$\sim$}" left] \& \CX_4 \oplus \CX_4 \arrow[l,"{\left[\begin{matrix} -(x_3-x_1) & x_4-x_3 \end{matrix}\right]}" below=0.5em] \arrow[u,"{(f,g) \mapsto (\psi_2\psi_1 \hat{v}_\tau f \tau_2, \psi_2 \psi_1 \psi_3 \hat{v}_\tau g \tau_2)}" right,"\rotatebox{90}{$\sim$}" left] \& q^2 \CX_4 \arrow[l,"{\left[\begin{matrix} x_4-x_3 \\ x_3-x_1 \end{matrix}\right]}"] \arrow[u,"f \mapsto \psi_3 \psi_2 \psi_3 \psi_1 \hat{v}_\tau f \tau_2" right,"\rotatebox{90}{$\sim$}" left] \& 0 \arrow[l]
\end{tikzcd}$$
Thus, denoting $\CZ_k := \k[z_1, \ldots, z_k]$, we have
\begin{align*}
	\CE_\theta(\pi,\rho) &= \CE_\theta^1(\pi,\rho) \cong q^{-1}\CX_2 / (x_1 = x_2) \iso q^{-1}\CZ_1, \ \bar{x}_1 \mapsto z_1 \\
	\CE_\theta(\rho,\sigma) &= \CE_\theta^1(\rho,\sigma) \cong q^{-1}\CX_3 / (x_1 = x_2) \iso q^{-1}\CZ_2, \ \bar{x}_1 \mapsto z_1, \bar{x}_3 \mapsto z_2 \\
	\CE_\theta(\sigma,\tau) &= \CE_\theta^1(\sigma,\tau) \cong q^{-2}\CX_4 / (x_3 = x_4) \iso q^{-2}\CZ_3, \ \bar{x}_1 \mapsto z_1, \bar{x}_2 \mapsto z_2, \bar{x}_3 \mapsto z_3 \\
	\CE_\theta(\pi,\tau) &= \CE_\theta^2(\pi,\tau) \cong q^{-2}\CX_4 / (x_1 = x_3 = x_4) \iso q^{-2}\CZ_2, \ \bar{x}_1 \mapsto z_1, \bar{x}_2 \mapsto z_2.
\end{align*}

\begin{Example}
	There is an $A_\infty$-category structure $\CE_\theta = (\CE_\theta, M_1, M_2, \ldots)$ satisfying the conditions in Theorem~\ref{TKad} such that
	$$\begin{tikzcd}[row sep = -0.5em]
		\CE_\theta^1(\sigma,\tau) \otimes \CE_\theta^1(\rho,\sigma) \otimes \CE_\theta^1(\pi,\rho) \arrow[r,"M_3^{\pi,\rho,\sigma,\tau}"] & \CE_\theta^2(\pi,\tau) \\
		\rotatebox{90}{$\cong$} & \rotatebox{90}{$\cong$} \\
		q^{-2}\CZ_3 \otimes q^{-1}\CZ_2 \otimes q^{-1}\CZ_1 \arrow[r] & q^{-2}\CZ_2 \\
		z_1^a z_2^b z_3^c \otimes z_1^m z_2^n \otimes z_1^w \arrow[r,mapsto] & z_1^{c+m+w} z_2^{b+n} \dfrac{z_2^a - z_1^a}{z_2-z_1}.
	\end{tikzcd}$$
	Moreover, there is no $A_\infty$-category structure on $\CE_\theta$ satisfying the conditions in Theorem~\ref{TKad} with $M_3 = 0$.
\end{Example}
\begin{proof}
	We apply Algorithm~\ref{AKad}.
	We will need to examine the complexes $T_\bullet^\pi(\sigma)$ and $T_\bullet^\rho(\tau)$.
	They are, respectively, the top complexes in the two diagrams below, and the diagrams define isomorphisms of complexes.
	$$\begin{tikzcd}[column sep=2em, ampersand replacement=\&]
		\& [-1em] 0 \& q^{-1}1_{2132}\Delta(\sigma) \arrow[l] \& [2em] 1_{2123}\Delta(\sigma) \arrow[l,"{\left[(d_0^\pi)_{2,1}\right]}" above] \& [-1em] 0 \arrow[l] \\
		\& 0 \& \CX_3 \arrow[l] \arrow[u,"f \mapsto \psi_2 \psi_1 v_\sigma f" right,"\rotatebox{90}{$\sim$}" left] \& q^2 \CX_3 \arrow[l,"{\left[x_3-x_1\right]}" above] \arrow[u,"f \mapsto \psi_3 \psi_2 \psi_1 v_\sigma f" right,"\rotatebox{90}{$\sim$}" left] \& 0 \arrow[l] \\
		0 \& q^{-1}1_{3221}\Delta(\tau) \arrow[l] \& 1_{2321}\Delta(\tau) \oplus 1_{3212}\Delta(\tau) \arrow[l,"d_1^\rho" above] \& q 1_{2312}\Delta(\tau) \arrow[l,"d_0^\rho" above] \& 0 \arrow[l] \\
		0 \& q^{-3}\CX_4 \arrow[l] \arrow[u,"f \mapsto \hat{v}_\tau f \tau_2" right,"\rotatebox{90}{$\sim$}" left] \& q^{-1}\CX_4 \oplus q^{-1}\CX_4 \arrow[l,"{\left[\begin{matrix} -(x_2-x_1) & x_4-x_3 \end{matrix}\right]}" below=0.5em] \arrow[u,"{(f,g) \mapsto (\psi_1 \hat{v}_\tau f \tau_2, \psi_3 \hat{v}_\tau g \tau_2)}" right,"\rotatebox{90}{$\sim$}" left] \& q\CX_4 \arrow[l,"{\left[\begin{matrix} x_4-x_3 \\ x_2-x_1 \end{matrix}\right]}"] \arrow[u,"f \mapsto \psi_1 \psi_3 \hat{v}_\tau f \tau_2" right,"\rotatebox{90}{$\sim$}" left] \& 0 \arrow[l]
	\end{tikzcd}$$
	so that
	\begin{align*}
		\CE_\theta(\pi,\sigma) &= \CE_\theta^1(\pi,\sigma) \cong \CX_3 / (x_1=x_3) \iso \CZ_2, \ \bar{x}_1 \mapsto z_1, \bar{x}_2 \mapsto z_2 \\
		\CE_\theta(\rho,\tau) &= \CE_\theta^2(\rho,\tau) \cong q^{-3}\CX_4/(x_1=x_2,x_3=x_4) \iso q^{-3}\CZ_2, \bar{x}_1 \mapsto z_1, \bar{x}_3 \mapsto z_2.
	\end{align*}

	Instead of $F_1$, it will only be relevant to determine the restrictions $F_1^{\pi,\rho}$, $F_1^{\rho,\sigma}$, $F_1^{\sigma,\tau}$, and $F_1^{\rho,\tau}$.
	According to the algorithm, we are free to take for these any cycle-choosing homomorphisms.
	We define
	\begin{alignat*}{2}
		&{}F_1^{\pi,\rho}(z_1^a) := &
		\begin{tikzcd}[column sep=2em, row sep=3em, ampersand replacement=\&]
			0 \arrow[r] \& q^2 R_\theta 1_{2321} \arrow[r,"d_1^\pi"] \arrow[dr,"{\left[\begin{matrix} y_4^a 1_{2321} & 0 \end{matrix}\right]}" right=1.2em] \& q R_\theta 1_{2231} \oplus q R_\theta 1_{2132} \arrow[r,"d_0^\pi"] \arrow[dr,"{\left[\begin{matrix} 0 \\ -\psi_2 y_3^a 1_{2312} \end{matrix}\right]}" right=2em] \& R_\theta 1_{2123} \arrow[r] \& 0, \\
			0 \arrow[r] \& q^2 R_\theta 1_{3221} \arrow[r,"d_1^\rho" below] \& q R_\theta 1_{2321} \oplus q R_\theta 1_{3212} \arrow[r,"d_0^\rho" below] \& R_\theta 1_{2312} \arrow[r] \& 0,
		\end{tikzcd} \\
		&{}F_1^{\rho,\sigma}(z_1^a z_2^b) := &
		\begin{tikzcd}[column sep=2em, row sep=3em, ampersand replacement=\&]
			0 \arrow[r] \& q^2 R_\theta 1_{3221} \arrow[r,"d_1^\rho"] \arrow[dr,"{\left[\begin{matrix} y_2^a y_4^b 1_{3221} \end{matrix}\right]}" right=1.2em] \& q R_\theta 1_{2321} \oplus q R_\theta 1_{3212} \arrow[r,"d_0^\rho"] \arrow[dr,"{\left[\begin{matrix} 0 \\ -y_2^a y_3^b 1_{3212}
			\end{matrix}\right]}" right=2em] \& R_\theta 1_{2312} \arrow[r] \& 0, \\
			\& 0 \arrow[r] \& q R_\theta 1_{3221} \arrow[r,"d_0^\sigma" below] \& R_\theta 1_{3212} \arrow[r] \& 0,
		\end{tikzcd} \\
		&{}F_1^{\sigma,\tau}(z_1^a z_2^b z_3^c) := &
		\begin{tikzcd}[column sep=2em, row sep=3em, ampersand replacement=\&]
			0 \arrow[r] \& q R_\theta 1_{3221} \arrow[r,"d_0^\sigma"] \arrow[dr,"{\left[\begin{matrix} -y_1^a y_2^b y_4^c 1_{32^{(2)}1} \end{matrix}\right]}" right=1em] \& R_\theta 1_{3212} \arrow[r] \& 0, \\
			\& 0 \arrow[r] \& q^{-1} R_\theta 1_{32^{(2)}1} \arrow[r] \& 0,
		\end{tikzcd} \\
		&{}F_1^{\rho,\tau}(z_1^a z_2^b) := &
		\begin{tikzcd}[column sep=2em, row sep=3em, ampersand replacement=\&]
			0 \arrow[r] \& [-1em] q^2 R_\theta 1_{3221} \arrow[r,"d_1^\rho"] \arrow[drr,"{\left[\begin{matrix} -y_2^a y_4^b 1_{32^{(2)}1} \end{matrix}\right]}" right=3em] \& q R_\theta 1_{2321} \oplus q R_\theta 1_{3212} \arrow[r,"d_0^\rho"] \& [-1em] R_\theta 1_{2312} \arrow[r] \& [-1em] 0, \\
			\& \& 0 \arrow[r] \& q^{-1} R_\theta 1_{32^{(2)}1} \arrow[r] \& 0.
		\end{tikzcd}
	\end{alignat*}

	Using the above, we now show that the following choices are in accordance with the algorithm:
	\begin{equation}\label{EZeroThings}
		M_2^{\pi,\rho,\sigma} = 0, \quad F_2^{\pi,\rho,\sigma} = 0, \quad M_2^{\pi,\rho,\tau} = 0, \quad F_2^{\pi,\rho,\tau} = 0
	\end{equation}
	By our choices of $F_1^{\pi,\rho}$ and $F_1^{\rho,\sigma}$, we have $U_2^{\pi,\rho,\sigma} = 0$, so $M_2^{\pi,\rho,\sigma} = 0$, and according to (\ref{EFunctorSimplified}), we may take $F_2^{\pi,\rho,\sigma} = 0$.
	Since $\CE_\theta(\pi,\rho)$ is concentrated in homological degree $1$ and $\CE_\theta(\rho,\tau)$ is concentrated in homological degree $2$, the image of $U_2^{\pi,\rho,\tau} = m_2^{\pi,\rho,\tau}(F_1^{\rho,\tau}\otimes F_1^{\pi,\rho})$ is in $\CH_\theta^3(\pi,\tau)$, which is zero since $P_\bullet^\pi$ has length $2$.
	Thus $M_2^{\pi,\rho,\tau} = 0$ and according to (\ref{EFunctorSimplified}), we may take $F_2^{\pi,\rho,\tau} = 0$.

	We now have
	\begin{align*}
		U_3^{\pi,\rho,\sigma,\tau} &= m_2^{\pi,\sigma,\tau}(F_1^{\sigma,\tau}\otimes F_2^{\pi,\rho,\sigma}) - m_2^{\pi,\rho,\tau}(F_2^{\rho,\sigma,\tau}\otimes F_1^{\pi,\rho}) + F_2^{\pi,\sigma,\tau}(1^{\sigma,\tau}\otimes M_2^{\pi,\rho,\sigma}) \\
		&{} \qquad \qquad - F_2^{\pi,\rho,\tau}(M_2^{\rho,\sigma,\tau}\otimes 1^{\pi,\rho}) \\
		&= -m_2^{\pi,\rho,\tau}(F_2^{\rho,\sigma,\tau}\otimes F_1^{\pi,\rho}),
	\end{align*}
	so we only need to make a choice for $F_2^{\rho,\sigma,\tau}$.
	We have $U_2^{\rho,\sigma,\tau} = m_2^{\rho,\sigma,\tau}(F_1^{\sigma,\tau}\otimes F_1^{\rho,\sigma})$ so that $U_2^{\rho,\sigma,\tau}(z_1^a z_2^b z_3^c \otimes z_1^m z_2^n)$ is given by
	$$\begin{tikzcd}[column sep=2em, row sep=3em, ampersand replacement=\&]
		0 \arrow[r] \& q^2 R_\theta 1_{3221} \arrow[r,"d_1^\rho"] \arrow[drr,"{\left[\begin{matrix} -y_1^a y_2^{b+m} y_4^{c+n} 1_{32^{(2)}1}\end{matrix}\right]}" right=3em] \& q R_\theta 1_{2321} \oplus q R_\theta 1_{3212} \arrow[r,"d_0^\rho"] \& R_\theta 1_{2312} \arrow[r] \& 0 \\
		\& \& 0 \arrow[r] \& q^{-1} R_\theta 1_{32^{(2)}1} \arrow[r] \& 0.
	\end{tikzcd}$$
	Thus
	\begin{equation}\label{Erst}
		M_2^{\rho,\sigma,\tau}(z_1^a z_2^b z_3^c \otimes z_1^m z_2^n) = z_1^{a+b+m}z_2^{c+n},
	\end{equation}
	and $(F_1^{\rho,\tau}M_2^{\rho,\sigma,\tau} - U_2^{\rho,\sigma,\tau})(z_1^a z_2^b z_3^c \otimes z_1^m z_2^n)$ is given by the above diagram, except the diagonal arrow is right multiplication with
	$$\left[\begin{matrix}
		-y_2^{b+m} y_4^{c+n} (y_2^a - y_1^a) 1_{32^{(2)}1}
	\end{matrix}\right].$$
	We may now take $F_2^{\rho,\sigma,\tau}(z_1^a z_2^b z_3^c \otimes z_1^m z_2^n)$ to be:
	$$\begin{tikzcd}[column sep=2em, row sep=4em, ampersand replacement=\&]
		0 \arrow[r] \& q^2 R_\theta 1_{3221} \arrow[r,"d_1^\rho"] \& q R_\theta 1_{2321} \oplus q R_\theta 1_{3212} \arrow[r,"d_0^\rho"] \arrow[dr,"{\left[\begin{matrix} \psi_1 y_2^{b+m}y_4^{c+n}\frac{y_2^a-y_1^a}{y_2-y_1} 1_{32^{(2)}1} \\ 0 \end{matrix}\right]}" right=2em] \& R_\theta 1_{2312} \arrow[r] \& 0 \\
		\& \& 0 \arrow[r] \& q^{-1} R_\theta 1_{32^{(2)}1} \arrow[r] \& 0.
	\end{tikzcd}$$

	Now $U_3^{\pi,\rho,\sigma,\tau}(z_1^a z_2^b z_3^c \otimes z_1^m z_2^n \otimes z_1^w)$ is given by:
	$$\begin{tikzcd}[column sep=2em, row sep=4em, ampersand replacement=\&]
		0 \arrow[r] \& q^2 R_\theta 1_{2321} \arrow[r,"d_1^\pi"] \arrow[drr,"{\left[\begin{matrix} -\psi_1 y_2^{b+m}y_4^{c+n+w}\frac{y_2^a-y_1^a}{y_2-y_1} 1_{32^{(2)}1} \end{matrix}\right]}" right=2em] \& q R_\theta 1_{2231} \oplus q R_\theta 1_{2132} \arrow[r,"d_0^\pi"] \& R_\theta 1_{2123} \arrow[r] \& 0, \\
		\& \& 0 \arrow[r] \& q^{-1} R_\theta 1_{32^{(2)}1} \arrow[r] \& 0,
	\end{tikzcd}$$
	so that $M_3^{\pi,\rho,\sigma,\tau} = [U_3^{\pi,\rho,\sigma,\tau}]$ is as in the theorem statement.

	For the second assertion, note that the existence of a second $A_\infty$-category structure $\CE_\theta = (\CE_\theta, N_1, N_2, \ldots)$ satisfying the conditions of Theorem~\ref{TKad} implies the existence of an isomorphism of $A_\infty$-categories $G: (\CE_\theta, M_1, M_2, \ldots) \to (\CE_\theta, N_1, N_2, \ldots)$ with $G_1$ being the identity on each morphism space.
	Assume, toward a contradiction, that such an isomorphism exists, and that $N_3 = 0$.

	Recall that we take $M_1 = N_1 = 0$, so ($**_2$) applied to $G$ reads $G_1 M_2 = N_2(G_1\otimes G_1)$, and since $G_1$ is the identity, we have $M_2 = N_2$.
	Now since $N_3 = 0$ by assumption, ($**_3$) reads
	\begin{equation}\label{EGRelation}
		G_2(M_2\otimes 1 - 1\otimes M_2) + M_3 = M_2(1\otimes G_2 - G_2\otimes 1).
	\end{equation}
	The restriction of (\ref{EGRelation}) to $\CE_\theta(\sigma,\tau)\otimes \CE_\theta(\rho,\sigma)\otimes\CE_\theta(\pi,\rho)$ is
	\begin{align}\label{EGRestricted}
		&{}G_2^{\pi,\rho,\tau}(M_2^{\rho,\sigma,\tau}\otimes 1^{\pi,\rho}) - G_2^{\pi,\sigma,\tau}(1^{\sigma,\tau}\otimes M_2^{\pi,\rho,\sigma}) + M_3^{\pi,\rho,\sigma,\tau} \nonumber \\
		&{}\qquad\qquad = M_2^{\pi,\sigma,\tau}(1^{\sigma,\tau}\otimes G_2^{\pi,\rho,\sigma}) - M_2^{\pi,\rho,\tau}(G_2^{\rho,\sigma,\tau}\otimes 1^{\pi,\rho}),
	\end{align}
	so according to (\ref{EZeroThings}), (\ref{EGRestricted}) becomes
	\begin{equation}\label{EGSimplified}
		G_2^{\pi,\rho,\tau}(M_2^{\rho,\sigma,\tau}\otimes 1^{\pi,\rho}) + M_3^{\pi,\rho,\sigma,\tau} = M_2^{\pi,\sigma,\tau}(1^{\sigma,\tau}\otimes G_2^{\pi,\rho,\sigma}).
	\end{equation}

	We have a formula for $M_2^{\rho,\sigma,\tau}$ given by (\ref{Erst}).
	We define
	$$F_1^{\pi,\sigma}(z_1^m z_2^n) :=
	\begin{tikzcd}[column sep=2em, row sep=3em, ampersand replacement=\&]
		0 \arrow[r] \& q^2 R_\theta 1_{2321} \arrow[r,"d_1^\pi"] \arrow[dr,"{\left[\begin{matrix} \psi_1 y_2^n y_4^m 1_{3221} \end{matrix}\right]}" right=1.2em]
		\& q R_\theta 1_{2231} \oplus q R_\theta 1_{2132} \arrow[r,"d_0^\pi"] \arrow[dr,"{\left[\begin{matrix} 0 \\ -\psi_2 \psi_1 y_2^n y_3^m 1_{3212} \end{matrix}\right]}" right=2em]
		\& R_\theta 1_{2123} \arrow[r] \& 0, \\
		\& 0 \arrow[r] \& q R_\theta 1_{3221} \arrow[r,"d_0^\sigma" below] \& R_\theta 1_{3212} \arrow[r] \& 0,
	\end{tikzcd}$$
	By our choices of $F_1^{\pi,\sigma}$ and $F_1^{\sigma,\tau}$, we have that $U_2^{\pi,\sigma,\tau}(z_1^a z_2^b z_3^c \otimes z_1^m z_2^n)$ is given by
	$$\begin{tikzcd}[column sep=2em, row sep=3em, ampersand replacement=\&]
		0 \arrow[r] \& q^2 R_\theta 1_{2321} \arrow[r,"d_1^\pi"] \arrow[drr,"{\left[\begin{matrix} -\psi_1 y_1^a y_2^{b+n}y_4^{c+m} 1_{32^{(2)}1} \end{matrix}\right]}" right=2em] \& q R_\theta 1_{2231} \oplus q R_\theta 1_{2132} \arrow[r,"d_0^\pi"] \& R_\theta 1_{2123} \arrow[r] \& 0, \\
		\& \& 0 \arrow[r] \& q^{-1} R_\theta 1_{32^{(2)}1} \arrow[r] \& 0
	\end{tikzcd}$$
	so that $M_2^{\pi,\sigma,\tau}(z_1^a z_2^b z_3^c \otimes z_1^m z_2^n) = z_1^{a+c+m} z_2^{b+n}$.

	Denote the left- and right-hand sides of (\ref{EGSimplified}) by $L$ and $R$, respectively.
	We apply $L$ and $R$ to two elements: $z_1\otimes 1 \otimes 1$ and $z_2\otimes 1 \otimes 1$.
	Note that the map
	$$G_2^{\pi,\rho,\sigma} : q^{-1}\CZ_2\otimes q^{-1}\CZ_1 \cong \CE_\theta^1(\rho,\sigma) \otimes \CE_\theta^1(\pi,\rho) \to \CE_\theta^1(\pi,\sigma) \cong \CZ_2,$$
	is a KLR degree $0$ map, so we must have $G_2^{\pi,\rho,\sigma}(1\otimes 1) = 0$.
	Thus, $R(z_1\otimes 1 \otimes 1) = R(z_2\otimes 1 \otimes 1) = 0$.
	On the other hand, we have $M_3^{\pi,\rho,\sigma,\tau}(z_1\otimes 1\otimes 1) = 1$, $M_3^{\pi,\rho,\sigma,\tau}(z_2\otimes 1\otimes 1) = 0$, and $M_2^{\rho,\sigma,\tau}(z_1\otimes 1) = M_2^{\rho,\sigma,\tau}(z_2\otimes 1) = z_1$, so by (\ref{EGSimplified}), we have
	\begin{align*}
		0 &= R(z_1 \otimes 1 \otimes 1) = L(z_1 \otimes 1 \otimes 1) = G_2(z_1\otimes 1) + 1, \\
		0 &= R(z_2 \otimes 1 \otimes 1) = L(z_2\otimes 1 \otimes 1) = G_2(z_1\otimes 1),
	\end{align*}
	a clear contradiction.
\end{proof}

\newcommand{\bib}[2]{#1}

\bib{
\bibliographystyle{abbrv}
\bibliography{divresTwo}
}{

}

\end{document}